\documentclass[10pt]{article}
\usepackage[sumlimits]{amsmath}
\usepackage{amsmath,amsfonts,amsthm,amssymb,amscd}
\newtheorem{theo}{Theorem}[section]
\newtheorem{prop}[theo]{Proposition}
\newtheorem{lemm}[theo]{Lemma}
\newtheorem{coro}[theo]{Corollary}

\newtheorem{defi}[theo]{Definition}

\newcommand\N{\mathbb{N}}
\newcommand\Z{\mathbb{Z}}

\newcommand\R{\mathbb{R}}

\title{On bounded pseudodifferential operators in Wiener
spaces }

\bigskip

\author{Laurent Amour, Lisette Jager and Jean Nourrigat}

\date{Universit\'e de Reims}

\begin{document}

\maketitle

%%%%%%%%%%%%%%%%%%%%%%%%%%%%%%%%%%%%%%%%%%%%%%%%%%%%%%%%ù
% A mettre 

%\subjclass[2010]{Primary  35S05 }
%\dedicatory{Dedicated to the memory of Bernard Lascar}

%%%%%%%%%%%%%%%%%%%%%%%%%%%%%%%%%%%%%%%%%%%%%%%%%%%%%%%%ù
\date{}

\begin{abstract}
We aim at extending the definition of the Weyl calculus to
an infinite
 dimensional setting, by replacing the phase space $ \R^{2n}$ by
 $ B^2$, where $ (i,H,B)$ is an abstract
 Wiener space. A first approach is to generalize the integral
 definition using the  Wigner function. The symbol 
 is then a function defined
 on $B^2$ and belonging to a $L^1$ space for a gaussian measure,
 the Weyl operator is defined as a quadratic form on a dense subspace of
$L^2(B)$. For example, the symbol can be the stochastic extension on $B^2$,
in the sense of L. Gross, of a function $F$ which is continuous and bounded on
 $H^2$.

In the second approach, this function $F$ defined on $H^2$  satisfies
 differentiability conditions analogous to the finite dimensional ones. One
 needs to introduce hybrid operators acting as  Weyl operators on the
 variables of
finite dimensional subset of $H$ and as  Anti-Wick operators on the rest
 of the variables. The final Weyl operator is then defined as a limit
and it is continuous on a $L^2$ space. Under rather weak conditions, it
 is an extension of the operator defined by the first approach.

We give examples of monomial symbols
linking this construction to the classical pseudodifferential operators theory
and other examples related to other fields or previous works on this subject.

\end{abstract}

\tableofcontents
%%%%%%%%%%%%%%%%%%%%%%%%%%%%%%%%%%%%%%%%%%%%%%%%%%%%%%%%%%%%%%%%%
%%%%%%%%%%%%%%%%%%%%%%%%%%%%%%%%%%%%%%%%%%%%%%%%%%%%%%%%%%%%%%%%%
%%%%%%%%%%%%%%%%%%%%%%%%%%%%%%%%%%%%%%%%%%%%%%%%%%%%%%%%%%%%%%%%%
%%%%%%%%%%%%%%%%%%%%%%%%%%%%%%%%%%%%%%%%%%%%%%%%%%%%%%%%%%%%%%%%%
%%%%%%%%%%%%%%%%%%%%%%%%%%%%%%%%%%%%%%%%%%%%%%%%%%%%%%%%%%%%%%%%%
\section{Introduction}\label{1}

\bigskip
The Weyl calculus associates with every function $F$, called symbol, bounded and
measurable on
 $\R^{2n}$ and every $h>0$ an operator,
 denoted by  $Op_h^{Weyl, Leb}(F)$,
from the  Schwartz space ${\cal S}(\R^n)$ into its dual space ${\cal
S}'(\R^n)$. For a function $f$ belonging to ${\cal S}(\R^n)$, this
 operator is, formally,  defined by

\begin{equation}
\label{(1.1)} 
(  Op_h^{Weyl, Leb}(F)f) (x) = (2 \pi h)^{-n} \int _{\R^{2n}} 
e^{\frac{i}{h}(x-y) \cdot \xi } F\left ( \frac{x+y}{2}, \xi \right ) f(y)
d\lambda ( y, \xi )
\end{equation}
where $\lambda$ is the Lebesgue measure. If $F$ is  $C^{\infty}$ and
bounded  on $\R^{2n}$ as well as all its derivatives, then
Calder\'on-Vaillancourt's theorem states that  $ Op_h^{Weyl,Leb}(F)$ extends to
 a bounded operator on $L^2(\R^n,\lambda)$ 
(see \cite{C-V}, \cite{HO}).

\bigskip

The definition of the Weyl operator as an application on ${\cal S}(\R^n)$
with values in  ${\cal S}'(\R^n)$ or, equivalently,  as a quadratic form 
on  ${\cal S}(\R^n)$, has already been extended to the infinite-dimensional case
for some specific symbols
by  Kree-R\c aczka \cite{K-R} and, up to a small
modification, by Bernard Lascar (\cite{LA-1}),
  (see \cite{LA-2} to \cite{LA-10} as well for applications).
In the present paper the hypotheses on the symbol of the operator 
(the function $F$ in (\ref{(1.1)})) are weaker than by these authors.
We also give a  Calder\'on-Vaillancourt type result in this context.

\bigskip

 The classical definition (\ref{(1.1)}) does not lend itself to an extension 
to an infinite dimensional case. We shall use instead the definition of
$Op_h^{Weyl, Leb}(F)$ in which the Wigner function appears. This operator
is the only one that satisfies, for all $f$ and $g$ in
${\cal S}(\R^n)$:
 \begin{equation}
  <  Op_h^{Weyl, Leb}(F)  f, g>_{L^2(\R^n,\lambda)}  =(2\pi h) ^{-n}
  \int _{\R^{2n}} F(Z)H_h ^{Leb}(f,g,Z) d\lambda (Z), \label{(1.2)} 
  \end{equation} 
where $H_h^{Leb}(f,g,.)$ is the Wigner function (for the Lebesgue
 measure):
\begin{equation}
H_h^{Leb} (f,g,Z) = \int _{\R^n} e^{-{\frac{i}{h} } t\cdot \zeta } f \left (
 z + \frac{t}{2} \right ) \overline {g \left (
 z - \frac{t}{2} \right )} d\lambda (t) \hskip 2cm Z =
 (z,\zeta) \in \R^{2n}  \label{(1.3)}
 \end{equation}
 (cf Unterberger \cite{U-2}, or Lerner  \cite{LE}, sections 2.1.1 and 2.1.2,
  or   Combescure  Robert  \cite{C-R}, section 2.2).

\bigskip

The expression (\ref{(1.2)}) is well defined for $f$  and $g$ in  
 ${\cal S}(\R^n)$ even if $F$ is only supposed to be measurable and bounded.
Indeed, the Wigner function  $H_h^{Leb} (f,g,\cdot )$  is in
${\cal S}(\R^{2n})$. We  first extend this definition  (\ref{(1.2)}) 
to the infinite dimensional case (using Wiener spaces) in the same spirit,
by defining the operator as a bilinear form on a convenient dense subspace,
with very weak assumptions on $F$ (weaker than those of 
\cite{LA-1} and \cite{K-R}) : see definition \ref{1.2} below.
To this aim we associate a Gaussian Wigner function $H_h^{Gauss} (f,g)$
 with every couple $(f,g)$, where $f$ and $g$ are convenient functions
 defined on a Wiener space. Then we prove that,
 under Calder\'on-Vaillancourt type conditions, 
this operator extends to a bounded 
operator on a $L^2$ space (Theorem \ref{1.4}, which is the main result).
 
\bigskip
We successively indicate what replaces the space $\R^n$,
the Lebesgue measure, the Schwartz's space  ${\cal S}(\R^n)$  and the Wigner
 function. Then we shall be ready to define the Weyl operator as a quadratic
 form, as in  (\ref{(1.2)}).

\bigskip

The space $\R^n$ is replaced by a real separable infinite-dimensional 
Hilbert space $H$, the configuration space. The symbol $F$ is a function 
on the phase space $H^2$. One denotes by  $a\cdot b$ the scalar product 
of two elements $a$ and $b$ of $H$ and by $|a|$ the norm of an element
$a$ of $H$.

\bigskip

The Lebesgue measure is replaced by a Gaussian measure associated with
 the norm of $H$. But, since $H$ is infinite-dimensional, this Gaussian 
measure will be a measure, not on $H$ (which is impossible) but on 
a convenient Banach space $B$ containing $H$. Abstract Wiener spaces are
 an adequate frame. An abstract Wiener space is a
 triple
$(i, H,B)$, where $H$ is a real, separable Hilbert space, $B$ a
Banach space and $i$ a continuous injection from $H$ into $B$ such
that $i(H)$ is dense in $B$, other conditions being satisfied (see
\cite{G-2,G-3,K} or Definition \ref{4.1}).
 The injection $i$ is generally
not mentioned and one usually identifies $H$ with its dual space, so
that the preceding hypotheses imply
\begin{equation}
 B' \subset H' = H \subset B, \label{(1.4)} 
 \end{equation}
where every space is continuously embedded as a dense subspace of the
 following one. For every $u$ in $B'$ and $x$ in $B$,
one denotes by $u(x)$ the duality between these elements and
one supposes that, if $x$ is in $H$, $u(x)= u\cdot x$,
This will be the case in the rest of the article.

\bigskip

One denotes by  ${\cal F}(H)$ (resp. ${\cal F}(B')$) the set of
all finite-dimensional subspaces $E$  of $H$ (resp. of  $B'$). We note that
 ${\cal F}(B') \subset {\cal F}(H) $ by (\ref{(1.4)}).
For every $E$ belonging to ${\cal F}(H)$ and every  positive $h$, 
one defines a probability measure $\mu_{E,h}$ on  $E$ by:
\begin{equation}
 d\mu _{E, h}(x) = (2 \pi h)^{-(1/2){\rm dim }(E)} e^{-\frac{|x|^2}{2h}}
d\lambda (x),\label{(1.5)} 
\end{equation}
where $\lambda$ is the Lebesgue measure on $E$,  the norm on $E$
being the restriction to $E$ of the norm of $H$.
 If $B$ is a Banach space satisfying (\ref{(1.4)}), the continuity and density
conditions being satisfied, one defines, for each  $E$ is ${\cal F}(B')$,
 an application $P_E : B \rightarrow E$ by:
\begin{equation}
 P_E(x) = \sum _{j= 1}^n u_j(x) u_j,
\label{(1.6)} 
\end{equation}
where  $\{ u_1,...,u_n\}$ is a basis of  $E$, orthonormal for
the restriction to $E$ of the scalar product of $H$. This
application does not depend on the choice of the basis. If the Banach space
 $B$ satisfies proper conditions, (see Definition \ref{4.1}
for precisions),
for each positive $h$,  one can derive a probability measure $\mu_{B ,h}$ 
on the Borel $\sigma-$algebra of $B$, with  the following property.
For every $E$ in  ${\cal F}(B')$  and every function
 $\varphi$ in $L^1(E,\mu _{E ,h})$, the function $\varphi \circ P_E$ is
 in $L^1(B,\mu _{B ,h})$, and one has
\begin{equation}
 \int _B (\varphi \circ P_E )(x) d\mu _{B,h} (x) = \int _E
\varphi (y)d\mu _{E,h} (y). \label{(1.7)}
\end{equation}

Other properties of the measure  $\mu_{B,h}$ are recalled in section
\ref{4.A}  and
examples are  given in Section \ref{8}.
In Section \ref{4.B}, the same notions will be seen for subspaces $E$ in
 ${\cal F}(H)$. The Banach space $B$ associated with $H$ is not unique but
for any choice of $B$ and any positive $h$, the space $L^2(B,\mu _{B ,h})$ is
 isomorphic to the symmetrized Fock space ${\cal F}_s (H)$,
which does not depend on the choice of $B$
 (cf \cite{F} or \cite{A-J-N-1}). If the Hilbert space $H$ is
 finite dimensional, then $B=H$. If not,  $B$ is sometimes derived from $H$
 thanks to a Hilbert-Schmidt
operator (cf \cite{G-3}, Example 2, p. 192)), but other constructions are
possible. In Section, \ref{8} examples related to Brownian motion,
field theory and interacting lattices  will be given.

\bigskip
Now let us introduce the space which will replace  ${\cal S}(\R^n)$ in 
 (\ref{(1.2)}) (\ref{(1.3)}). For every  $E$ in  ${\cal F}(H)$, we shall need
the isometric isomorphism $\gamma_{E, h/2}$ between
 $L^2(E,\mu _{E,h/2})$ and $L^2(E, \lambda)$ given by
\begin{equation}
 (\gamma_{E,h/2} \varphi ) (x) =  ( \pi h)^{-(1/4){\rm dim }(E)}
\varphi (x) e^{-\frac{|x|^2}{2h}} \hskip 2cm \varphi \in L^2(E,\mu _{E,h/2}) .
 \label{(1.8)}
\end{equation}

One denotes by ${\cal S}_E$ the space of all functions  $\varphi : E
\rightarrow {\bf C}$ such that  $\gamma_{E,h/2} \varphi$ belongs to the
Schwartz space of rapidly decreasing functions
 ${\cal S}(E)$.
For $E\in{\cal F}(B')$, let ${\cal D}_E$ be the set of applications $f: B \rightarrow {\bf
C}$ of the form
 $f = \varphi \circ P_E$, where $P_E : B \rightarrow E$ is defined by
(\ref{(1.6)}) and  $\varphi : E \rightarrow {\bf C}$ belongs to
 ${\cal S}_E$. We denote by   ${\cal D}$ the union of the spaces
 ${\cal D}_E$, taken over all  $E$ in ${\cal F}(B')$.
This space  ${\cal D}$ is dense in $L^2(B,\mu _{B,h/2})$.
Indeed, if $(e_j)_{(j\geq 1)}$ is a Hilbert basis of $H$, the
vectors $e_j$ belonging to $B'$, the set of functions on $B$ which
are polynomials of  a finite number of functions $x \rightarrow
e_j(x)$ is contained in  ${\cal D}$ and is dense in $L^2(B,\mu_{B
, h/2})$ (see for example \cite{A-J-N-1}).The constant functions belong to
 ${\cal D}_E$ for every  $E$ in  ${\cal F}(B')$.

\bigskip

We very often need the following classical result (Kuo, \cite{K},
Chapter 1 Section 4).

\begin{prop}\label{1.1} {\ }
\begin{enumerate}
\item
 The space $B'$ is contained in  $L^2(B ,\mu _{B,h})$ and  the norm of
 $u\in B'$, considered as an element of
 $L^2(B,\mu _{B,h})$, is equal to $\sqrt{h}| u |$.
\item
The inclusion of $B'$ into $L^2(B,\mu _{B,h})$
extends to a continuous linear map  from $H$ into $L^2(B,\mu _{B ,h})$,
 with norm $\sqrt {h}$, denoted by $u\rightarrow \ell_u$.
\end{enumerate}
\end{prop}

The first point can be seen by applying  (\ref{(1.7)}) to a one dimensional
 space $E$.
The map $\ell $ turns $H$ into a Gaussian space in the sense of \cite{J} or into
a  "Gaussian random process", in the sense of  \cite{SI}.
 In the case of Example 8.2,
for every $u$ in the Cameron-Martin space,  $\ell_u$  is the It\^o
integral of the function $u'$.

\bigskip

Let us now define the Wigner-Gauss function, which will replace the  usual
Wigner function. Let  $(i, H, B)$  be a Wiener space satisfying (\ref{(1.4)}).
For every subspace $E$ in  ${\cal F}(H)$, for all $\varphi $ and $\psi$ in
 ${\cal S}_E$, one defines a function 
 $\widehat  H^{Gauss}_h (\varphi,\psi) $ on  $E^2$, setting:
\begin{equation} 
 \widehat  H^{Gauss}_h (\varphi,\psi) (z, \zeta) =
 e^{\frac{|\zeta|^2 }{ h}} \int _E e^{-{2\frac{i}{h}}   \zeta \cdot t  } \varphi (z+t)
 \overline {\psi (z-t)} d\mu _{E,h/2} (t)
\hskip 2cm (z,\zeta) \in E^2. \label{(1.9)}
\end{equation}
 One notices that, for $Z$ in $E^2$:
 \begin{equation} 
\widehat H^{Gauss}_h (\varphi ,\psi ) (Z) = 2^{-{\rm dim}(E)} e^{\frac{|Z|^2}{ h}}
H_h^{Leb} ( \gamma_{E,h/2} \varphi  ,\gamma_{E,h/2} \psi ) (Z)
\hskip 2cm Z \in E^2 . \label{(1.10)} 
 \end{equation}
One will see (Proposition \ref{4.8}) that this function belongs to
 $L^1(E^2,\mu_{E^2,h/2})$ and to $L^2(E^2,\mu_{E^2,h/4})$ as well.
For all functions  $f$ and $g$ in
 ${\cal D}_E$, ($E$ in ${\cal F}(B')$), of the form
 $f = \widehat f \circ P_E$  et $g = \widehat g \circ P_E$,
 with $\widehat f $ and $\widehat g$ in ${\cal S}_E$,
one defines the Wigner-Gauss transform
 $ H^{Gauss}_h (f,g) $, which is the function defined on $E^2$ by:
\begin{equation}
 H^{Gauss}_h (f,g) (Z) = \widehat H^{Gauss}_h (\widehat f,\widehat g) (P_EZ)
\hskip 2cm  a.e Z\in B^2. \label{(1.11)} 
\end{equation}
One writes $P_E$ instead of $P_{E^2}$ for the sake of simplicity.
 One will see in Section \ref{4.B} how to modify this definition if 
 $E \in {\cal F}(H)$. According to (\ref{(1.7)}), it follows from Proposition
\ref{4.8} that this function is in
 $L^1(B^2,\mu_{B^2,h/2})$ and  $L^2(B^2,\mu_{B^2,h/4})$. If $f$ and $g$ are
in  ${\cal D}$, the subspace $E$ such that $f$ and $g$ are in
 ${\cal D}_E$ is not unique, but the function defined above does not depend on
 $E$. Proposition \ref{4.8} states that this Wigner-Gauss transformation
extends, by density, from  $L^2(B,\mu_{B,h/2})\times L^2(B,\mu_{B,h/2}) $ 
to  $L^2(B^2, \mu_{B^2,h/4})$ and to the space of continuous functions defined
 on $H^2$. One will see, in  (\ref{(4.E.6)}),  another expression of
 $H^{Gauss}_h (f,g) $, using  Segal-Bargmann transforms of $f$ and $g$.

\bigskip

Now we are almost ready to define the Weyl operator associated with a symbol
 $F$. If $(i, H, B)$ is a Wiener space, we have the choice of two 
phase spaces:  $H^2$ an $B^2$. The first one is equipped with the
symplectic form
 $\sigma ((x,\xi), (y, \eta)) = y\cdot \xi - x \cdot \eta $, but not with a
 measure adapted to the integrations we want to conduct.
On the contrary, the space $B^2$ is equipped with the Gaussian measure
 $\mu_{B^2,h}$, but it is not a symplectic space.

\bigskip

This difficulty is overcome in the following way. The symbols will be 
initially defined as functions on $H^2$. There exists an operation, 
introduced by L. Gross (\cite{G-1}, see Ramer \cite{RA}) and usually called
stochastic extension,
which associates with a Borel function $F$ on $H$ or $H^2$  a Borel function
 $\widetilde F$ on $B$ or on $B^2$.
This stochastic extension, which  is not a genuine extension,  will be recalled
in Definition \ref{4.4}. What will appear in the initial definition formula
(\ref{(1.13)})
of the Weyl operators is the stochastic extension $\widetilde F$ of the
initial symbol $F$.

\bigskip

In the second step, concerned with the bounded extension in
$L^2(B,\mu _{B, h/2})$ of this initial operator, we shall restrict ourselves to
bounded symbols.  But in the initial definition of the operator as a quadratic
form on  ${\cal D}$, a polynomial growth will be enough.
For the initial definition (\ref{(1.13)}), the function $\widetilde{F}$ 
(stochastic
extension of the  symbol $F$) will be in $L^1 ( B^2, \mu_{B^2 ,h/2})$
and the polynomial growth  will be expressed in terms of the existence of a
nonnegative integer   $m$ such that the following norm is defined.
\begin{equation}
 N_m ( \widetilde{F})=  \sup _{Y\in H^2 }
\frac{\Vert \widetilde{F}(\cdot + Y)\Vert_{L^1( B ^2,\mu_{B^2,h/2})}} {(1+|Y|)^m }.
  \label{(1.12)} 
\end{equation}
This norm is finite if the function $F$ is bounded or if it is a polynomial
 expression of degree $m$ with respect to functions
 $(x,\xi) \rightarrow \ell _a (x) + \ell _b (\xi)$, with  $a$ and $b$
in $H$.

\begin{defi}\label{1.2}
  Let $(i, H, B)$ be an abstract Wiener
space satisfying (\ref{(1.4)}) and $h>0$. Let $\widetilde{F}$ be a  function
in $L^1 ( B ^2,\mu_{B^2,h/2})$. Suppose there exists $m\geq 0$ such that
 the norm
$N_m (\widetilde{F})$ is finite.  We define a quadratic form
 $Q^{Weyl}(\widetilde{F}) $ on  ${\cal D}
\times {\cal D}$ in the following way. For  all $f$ and $g$ in
 ${\cal D}$, one sets:
\begin{equation}
Q_h^{Weyl}(\widetilde{F})  (f, g)
= \int _{B^2} \widetilde{F} (Z) H^{Gauss }_h (f,g) ( Z)
d\mu_{B^2,h/2} (Z) \label{(1.13)} 
\end{equation}
where  $H^{Gauss }_h (f,g)$ is defined in (\ref{(1.11)}).
\end{defi}

The convergence will be proved in Proposition \ref{4.10}. If $\widetilde{F}$
is bounded, the convergence is a consequence  of Proposition \ref{4.8},
since $H^{Gauss }_h (f,g) $ is in $L^1 (B^2, \mu_{B^2, h/2})$. If we only
suppose  that $N_m(\widetilde{F})$ is finite, other arguments are necessary.

 \bigskip

One sees the relationship with the classical definition   (\ref{(1.2)}),
 (\ref{(1.3)}). If $F$ has the form $F= \widehat{F} \circ P_E$, where
 $E$ is in  ${\cal F}(B')$ and $\widehat{F}$ is a measurable bounded function
on $E ^2$, if $f = \widehat{f}\circ P_E$ and  $g = \widehat{g} \circ P_E$,
where $\widehat{f}$ et $\widehat{g}$ are in  ${\cal S}_E$, one has
\begin{equation}
Q_h^{Weyl}(\widehat{F}  \circ P_E) (\widehat{f} \circ P_E, \widehat{g} \circ P_E  )
= < Op_h^{Weyl, Leb} (\widehat{F}  )  \gamma_{E,h/2} \widehat{f},  \gamma_{E,h/2} \widehat{g}
>_{L^2 (E,\lambda)}  \label{(1.14)} 
\end{equation}
where $Op_h^{Weyl, Leb} (\widehat{F}  ) $ is defined on $E$
as it is in (\ref{(1.2)}), (\ref{(1.3)})
on $\R^n$. This equality comes from (\ref{(1.10)}).

\bigskip

One can compare this definition with the definition in \cite{K-R} or 
\cite{LA-1}.
The authors define an anti-Wick operator associated with a symbol,
a function defined on $H^2$ which is, for example, the Fourier transform
of a complex measure bounded on $H^2$. In \cite{K-R} they associate, too, with
 such a symbol, a Weyl operator defined as a quadratic form on a dense
subspace. When the symbol $F$ is  the Fourier transform
of a complex measure bounded on $H^2$, we prove, in Proposition \ref{8.6}, 
that $F$ admits a stochastic extension $\widetilde F$, with which
Definition \ref{1.2} above associates a Weyl operator, defined as a quadratic
 form. For this kind of symbols, the Weyl operator (as a quadratic form)
is explicitly written in (\ref{(8.C.7)}) and (\ref{(8.C.3)}) or, equivalently,
in (\ref{(8.C.7)}) and (\ref{(8.C.9)}). This last form can be found in
\cite{K-R} as well. The anti-Wick operator can be found in \cite {K-R, LA-1}
in the form (\ref{(8.C.8)})-(\ref{(8.C.3)}). Our Definition \ref{1.2}
is more general, insofar as admitting a stochastic extension is
more general than being the Fourier transform of a bounded measure.
Section \ref{8.B} gives other examples of classes of functions
admitting stochastic extensions.

\bigskip

However, all functions do not have a stochastic extension. If $H$ is infinite
 dimensional, the norm function, which associates with every $x\in H$ its norm
$|x|$, admits no stochastic extension (see  Kuo, \cite{K}, Chapter 1, 
section 4).
In the same way, the function $x \rightarrow e^{-|x|^2}$, defined on $H$, has 
no stochastic extension. We shall not be able to define a pseudodifferential
operator, whose symbol would be  $F(x,\xi) =e^{-|x|^2}$. This operator would
 have to be the multiplication by $e^{-|x|^2}$, but this function, lacking 
an extension, is only defined on a negligible set (the Hilbert space $H$
is negligible in $B$, cf  \cite{K}) and this does not make sense.
On the contrary, the function $t \rightarrow e^{- (At)\cdot t}$, where
the operator $A$ is 
positive, symmetrical and trace-class on $H$, has a stochastic extension
(see Proposition \ref{8.10}).

\bigskip

It now remains to extend the bilinear form defined above on the
subspace  ${\cal D}$ to get a linear operator bounded on $L^2(B ,\mu
_{B,h/2})$. When the symbol is the Fourier transform of a measure
bounded on $H^2$ (case treated in \cite{K-R}), the upper bound on the  norm
is a consequence of Proposition \ref{8.13}. For other cases
 we must specify the hypotheses on  $F$,
which will strongly depend on the choice of a Hilbert basis of $H$.
 We can now state the hypothesis on the function $F$, which will be
the symbol in our  version of the  Calder\'on-Vaillancourt theorem.

\begin{defi}\label{1.3}
 Let  $(i, H,B)$  be an abstract Wiener space
 satisfying (\ref{(1.4)}). The norm of $H$ will hereafter be denoted by
 $|\cdot |$ and the scalar product of two elements $a$ and  $b$ of $H$
will be denoted by $a \cdot b$. The norm of an element of $H^2$ is denoted by
  $|\cdot |$ as well. For all $X = (x ,\xi)$ and $Y
= (y ,\eta)$ in $H^2$, we set
$\sigma (X,Y) = y\cdot \xi - x\cdot \eta$. We choose a Hilbert
basis $(e_j )_{(j\in \Gamma)}$ of $H$, each vector belonging to
$B'$, indexed by a countable set $\Gamma$. Set $u_j = (e_j,0)$ and
$v_j = (0, e_j)$ $(j\in \Gamma)$. A multi-index is a map
 $(\alpha,\beta )$ from $\Gamma $ into
$\N \times \N$ such that $\alpha_j = \beta _j = 0$ except for a
finite number
 of indices. Let $M$ be a nonnegative real number, $m$  a nonnegative integer
 and  $\varepsilon = (\varepsilon_j )_{(j \in \Gamma)}$ a family of
nonnegative real numbers.  One denotes by $ S_m(M, \varepsilon)$ the
set of bounded continuous functions $ F:H^2\rightarrow {\bf
C}$ satisfying the following condition. For every multi-index
 $(\alpha,\beta)$ such that $0 \leq \alpha_j \leq
m$ and $0 \leq \beta_j \leq m$ for all $j\in \Gamma$, the following
derivative
\begin{equation}
\partial_u^{\alpha}\partial_v^{\beta}  F =  \left [\prod _{j\in \Gamma }
\partial _{u_j} ^{\alpha_j} \partial _{v_j} ^{\beta_j}\right ]  F  \label{(1.15)}
\end{equation}
is well defined, continuous on
 $H^2$ and satisfies, for every $(x,\xi)$ in  $H^2$
\begin{equation}
\left | \left [\prod _{j\in \Gamma } \partial_{u_j} ^{\alpha_j}
 \partial _{v_j} ^{\beta_j}\right ]  F(x,\xi)
\right |    \leq M \prod _{j\in \Gamma } \varepsilon_j ^{\alpha_j +
\beta_j}\ . \label{(1.16)} 
\end{equation}
\end{defi}

The choice of indexing the basis  $(e_j)$ by an arbitrary countable
set $\Gamma$  is motivated  by possible applications in lattice
theory.
 Cordes \cite{C}, Coifman Meyer \cite{C-M} (for the standard quantization,
not for the Weyl one)
 and Hwang \cite{HW} remarked that, in the
Calder\'on-Vaillancourt bounds, it is enough to consider multi-indices
 $(\alpha,\beta )$ such that
 $0\leq \alpha_j \leq 1$ et $0\leq \beta_j \leq 1$ for every  $j$, which
 inspired
the definition of $S_m(M, \varepsilon)$. One finds in Section \ref{8}
examples of functions of  $S_m(M,\varepsilon)$ coming from
interacting lattices models.

\bigskip

The main result can be stated as follows.

\begin{theo}\label{1.4}
 Let  $(i, H,B)$  be an abstract Wiener
space satisfying (\ref{(1.4)}) and $h$ be a positive number. Let 
 $(e_j )_{(j\in \Gamma)}$ be a Hilbert space basis of $H$, each vector belonging
 to
$B'$, indexed by a countable set $\Gamma$. Let $F$  be a function on  $H^2$
satisfying the following two hypotheses.
\begin{description}
\item[(H1)]
 The function  $F$ belongs to the class $ S_2(M, \varepsilon)$
 of  Definition  \ref{1.3}, where
  $M$ is a nonnegative real number and $\varepsilon =
(\varepsilon_j )_{(j \in \Gamma)}$ a square summable family of
nonnegative real numbers.
\item[(H2)]
  We assume that $F$ has a stochastic extension
$\widetilde F$ with   respect to both
 measures  $\mu_{B^2,h}$ and  $\mu_{B^2 ,h/2}$ (see Definition  \ref{4.4}).
\end{description}
 Then there exists an operator,
denoted by  $Op_h^{Weyl}(F)$, bounded in  $L^2(B, \mu_{B, h/2})$,
such that, for all $f$ et $g$ in ${\cal D}$
\begin{equation}
 <Op_h^{Weyl}(F) f,g> = Q_h^{Weyl}(\widetilde F) (f,g),\label{(1.17)} 
\end{equation}
where the right hand side is defined by Definition
 \ref{1.2}. Moreover,  if $h$ is in  $(0, 1]$:
\begin{equation} \Vert Op_h^{Weyl}(F)\Vert _{{\cal L}(
L^2(B,\mu_{B,h/2}) )}\leq M  \prod_{j\in \Gamma} (1 + 81 \pi h
S_{\varepsilon} \varepsilon _j^2),\label{(1.18)} 
\end{equation}
where
\begin{equation}
 S_{\varepsilon} = \sup _{j\in \Gamma} \max (1 ,\varepsilon_j^2). \label{(1.19)}
\end{equation}
\end{theo}

One will see in Proposition \ref{8.7} that, if a function $F$ is in
 $S_1(M, \varepsilon)$ and if the sequence $\varepsilon =
(\varepsilon_j )_{(j \in \Gamma)}$ is summable, then 
 $F$ satisfies  the hypothesis (H2).
 Proposition \ref{8.14} gives an example of function in $S_m (M, \varepsilon)$,
inspired by the lattice theory.

\bigskip

The operator  $Op_h^{Weyl}(F)$ associated with $F\in S_2(M, \varepsilon)$ will 
not be defined by an integral expression, but as the limit, in
 ${\cal L}( L^2(B,\mu_{B,h/2}) )$, of a sequence of operators.
In order to define this sequence  we shall associate with each subspace  $E$ in
 ${\cal F}(B')$, an operator denoted by $Op_h^{hyb,E}(\widetilde{F})$ and
 bounded on $L^2(B,\mu _{B,h/2})$. This operator behaves as a Weyl operator 
on a set of variables and as an anti-Wick operator on the
other variables. It will have the
form $Q_h^{Weyl }(G)$, where $G$ is obtained from the extension
 $\widetilde F$  by applying a partial heat operator concerning the
 variables in $E^{\perp}$ (see Section \ref{2}).
 Then  we shall prove that, if one replaces $E$ by a sequence
 $E (\Lambda_n) = Vect (e_j) _{(j\in \Lambda_n)}$ where $(\Lambda _n)$ is an
increasing sequence of finite subsets of $\Gamma$ whose union is
$\Gamma$, then the sequence of operators $(Op_h ^{hyb, E(\Lambda_n)}(G))_n$ is
 a Cauchy sequence in ${\cal L}(L^2(B,\mu _{B ,h/2}))$. Its limit will be
 denoted by $Op_h^{Weyl}(F)$.
  We shall see that, if one restricts this operator
to get a bilinear form on  ${\cal D}$, it coincides with the one
defined by  Definition \ref{1.2}.
In particular it does not depend on the Hilbert
basis $(e_j)$ chosen to construct it, nor on the sequence  $(\Lambda_n)$.

 \bigskip

The  hybrid  operator associated with each finite
dimensional subspace is defined in  Section \ref{2}.
 Section \ref{3} presents more precisely the main steps of the proof of
 Theorem \ref{1.4}. Section \ref{4} first recalls classical facts 
about Wiener spaces, then gives a precise definition of the stochastic
extensions, recalls the Segal Bargmann transformation, gives the
necessary upper bounds about the Gaussian Wigner function and gives another
 expression of the Weyl and hybrid operators, more convenient to get
norm estimates. Section \ref{5} establishes the norm estimates which 
prove that the sequence constructed to approach $Op_h^{Weyl}(F)$ is indeed a
Cauchy sequence. In Section \ref{6} we prove that  $Op_h^{Weyl}(F)$ 
is really an extension of the initially defined operator (Definition
 \ref{1.2}). Section \ref{7} studies the Wick symbol of this operator.
Section \ref{8} gives examples of Wiener spaces, stochastic extensions
and symbols belonging to $S_m(M, \varepsilon)$.

\bigskip
This paper would not have been written without a workgroup organized by
 F. H\'erau, who induced us strongly to get acquainted with
Wiener spaces.

%%%%%%%%%%%%%%%%%%%%%%%%%%%%%%%%%%%%%%%%%%%%%%%%%%%%%%%%%%%%%%%%%
%%%%%%%%%%%%%%%%%%%%%%%%%%%%%%%%%%%%%%%%%%%%%%%%%%%%%%%%%%%%%%%%%
%%%%%%%%%%%%%%%%%%%%%%%%%%%%%%%%%%%%%%%%%%%%%%%%%%%%%%%%%%%%%%%%%
%%%%%%%%%%%%%%%%%%%%%%%%%%%%%%%%%%%%%%%%%%%%%%%%%%%%%%%%%%%%%%%%%
%%%%%%%%%%%%%%%%%%%%%%%%%%%%%%%%%%%%%%%%%%%%%%%%%%%%%%%%%%%%%%%%%
\section{Anti-Wick, Weyl and hybrid operators. }\label{2}

\subsection{Heat semigroup and anti-Wick operator}
\label{2.A}
%%%%%%%%%%%%%%%%%%%%%%%%%%%%%%%%%%%%%%%%%%%%%%%%%%%%%%%%%%%%%%%%%
%%%%%%%%%%%%%%%%%%%%%%%%%%%%%%%%%%%%%%%%%%%%%%%%%%%%%%%%%%%%%%%%%

In an Euclidean finite dimensional space $E$
one may associate with any  Borel bounded function
$F$  on $E ^2$, and for all $h>0$, an operator $ {Op}_h^{AW,Leb}(F)$ 
called the anti-Wick or Berezin-Wick operator. Instead of referring
to the usual definition, let us say that it is the only operator satisfying
\begin{equation}
 < {Op}_h^{AW, Leb}(F) f,g  >_{Leb} = < {Op}_h^{Weyl, Leb}
(e^{{\frac{h}{4}}\Delta }F) f,g  > _{Leb} \label{(2.1)}
\end{equation}
for all $f$ and $g$ in the  Schwartz space ${\cal S}(E)$.  This
operator has a bounded extension in $L^2(E,\lambda)$,
 also denoted by $ {Op}_h^{AW, Leb}(F)$, which verifies
\begin{equation}
 \Vert  {Op}_h^{AW, Leb}(F) \Vert _{{\cal L} (L^2(E ,
\lambda))} \leq \Vert F \Vert _{L^{\infty}(E ^2)}. \label{(2.2)} 
\end{equation}
This property is easier to extend and hence  is more convenient as a
 starting point than the usual definition. See \cite{F} or \cite{LE},
 Chapter 2.

\bigskip

To extend Definition (\ref{(2.1)}) to the finite dimensional setting,
we need an analogue of the heat semigroup on $B^2$, that is the family of
 operators $\widetilde H_t$ defined by:
\begin{equation}
(\widetilde H_{ t}F) (X) = \int _{B^2} F(X +Y) d\mu _{B^2, t} (Y)
\hskip 2cm X\in B^2 . \label{(2.3)}
\end{equation}
The operator  $\widetilde H_{ t}$ corresponds, in the finite dimensional case, 
to the operator $e^{\frac{t}{2}\Delta }$. For every $t>0$, the operator 
$\widetilde H_{ t} $  is bounded in the space of all bounded Borel functions
 on $B^2$ (see Hall \cite{HA}). Moreover, it is bounded and with a norm smaller
 than $1$ from  $L^p ( B^2,\mu _{B^2,t+h})$ in
 $L^p ( B^2,\mu _{B^2,h})$  ($1 \leq p < \infty $, $t>0$, $h>0$),
see Proposition \ref{4.5} below. 

\bigskip

For every positive $h$, one associates with each bounded Borel function $F$
on $B^2$ an operator called anti-Wick operator. It will first be defined as a
 quadratic form on ${\cal D}\times {\cal D}$, denoted by $Q_h^{AW}(F)$ and such
that
\begin{equation}
 Q_h^{AW}(F)(f,g) = Q_h^{Weyl}(\widetilde H_{h/2}  F) (f, g), \label{(2.4)}  
\end{equation}
where the  Weyl quadratic form is the one of Definition \ref{1.2}.
One will see in Corollary \ref{4.12} an expression which is closer to
the usual definition, as well as the fact that this quadratic form is
linked with a bounded operator, whose norm will be specified.

\bigskip

The following step is associating, with every subspace $E$ in ${\cal F}(B')$,
a hybrid operator, with the help of a partial heat operator.
This hybrid operator  behaves as a
Weyl operator associated with $F$ regarding some  functions,  and
as an anti-Wick operator associated with the same symbol
regarding some other functions.

%%%%%%%%%%%%%%%%%%%%%%%%%%%%%%%%%%%%%%%%%%%%%%%%%%%%%%%%%%%%%%%%%
%%%%%%%%%%%%%%%%%%%%%%%%%%%%%%%%%%%%%%%%%%%%%%%%%%%%%%%%%%%%%%%%%
\subsection{Wiener measure decomposition}\label{2.B}

The following proposition proved in Gross \cite{G-4} or in \cite{RA} allows to
split the variables for our  hybrid operators.

\begin{prop}\label{2.1}
 Let $(i, H,B)$  be an abstract Wiener space. In the following, the injection
 $i$ is implicitly understood, $H$ is identified with its dual space in such
 a way that (\ref{(1.4)}) is verified and we fix a subspace $E$  in
 ${\cal F}(B')$.  Let $E^{\perp}$ be the following space:
\begin{equation}
 E^{\perp} = \{  x\in B, \hskip 1cm u(x) =0 \ \ \ \ \ \forall u\in E\}.
\label{(2.XX)}
\end{equation}
Let $i_1$ denote the injection of $E^{\perp}\cap H$ in $E^{\perp}$. Then:
\begin{enumerate}
\item
 The system $(i_1, E^{\perp}\cap H,E^{\perp})$ is an abstract Wiener space.
 We denote by $\mu _{E^{\perp}, h}$ the Gaussian measure of
 parameter $h$ on $E^{\perp}$.
\item
 We have, for all $h>0$:
 \begin{equation}
 \mu_{B,h}= \mu _{E,h} \otimes \mu _{E^{\perp},h}. \label{(2.5)} 
 \end{equation}
\end{enumerate}
\end{prop}

The map $x \rightarrow (P_E (x), x- P_E(x))$ (where $P_E$ is defined in 
(\ref{(1.6)})) is a bijection between  $B$ and $E \times E^{\perp}$. We denote
 by $x = (x_E ,x_{E^{\perp}})$ the variable in $B$ and by
 $X = (X_E ,X_{E^{\perp}})$ the variable in $B^2$. Sometimes, we shall write
$X = X_E  + X_{E^{\perp}}$.

%%%%%%%%%%%%%%%%%%%%%%%%%%%%%%%%%%%%%%%%%%%%%%%%%%%%%%%%%%%%%%%%%
%%%%%%%%%%%%%%%%%%%%%%%%%%%%%%%%%%%%%%%%%%%%%%%%%%%%%%%%%%%%%%%%%
\subsection{Partial heat semigroup and hybrid operators}
\label{2.C}
If $E$ is in  ${\cal F}(B')$, one can define, as in (\ref{(2.3)}), on the one
 hand a heat semigroup  acting only on the variables in $E$ and, on the other
 hand, another one acting on the variables of $E^{\perp}$
 (defined in Proposition \ref{2.1}).
There are two kinds of operators acting on the variables in $E$: 
the first kind acts on a space of functions defined on $H$, the second
kind, on a space of functions defined on $B$.
For the operator acting on the variables of  $E^{\perp}$, the second version
 only is available.
 
\bigskip

Let $E$ be in ${\cal F}(B')$. For all $t>0$, we
define an operator $ H_{E,t}$ on the space of bounded continuous
functions on $H^2$, setting, when $ F$ is such a function:
\begin{equation}
( H_{E,t} F) (X) = \int _{E^2}  F( X +Y_{E} ) d\mu _{E^2, t} (Y_{E})\hskip 2cm
 X \in H^2 . \label{(2.6)}
\end{equation}
We define likewise an operator $\widetilde H_{E,t}$, defined by the same
 formula, but acting on the space of bounded Borel functions on $B^2$.

\bigskip

We may also define a partial heat semi-group, only acting on the
variables lying in the  subspace $E^{\perp}$. For all bounded Borel
functions $ F$   on $B^2$ and for all positive $t$, one may define a
function $\widetilde H_{E^{\perp},t}F$ on $B^2$ by setting, for
all $X$ in $B^2$ :
\begin{equation}
(\widetilde H_{E^{\perp},t}F) (X) = \int _{(E^{\perp})^2} F(X
+Y_{E^{\perp}}) d\mu _{(E^{\perp})^2, t} (Y_{E^{\perp}}). \label{(2.7)}
\end{equation}
\bigskip

The operators $H_{E,t}$ and $\widetilde H_{E,t}$ will mainly appear in
 Section \ref{3}.
The  $\widetilde H_{E ^{\perp}, t}$ will be used now and in Section \ref{6}.
\bigskip

 We are now ready to define the hybrid operator.

\begin{defi}\label{2.2}
 If $F$ is a bounded Borel function on $B^2$ and if $E$ is in ${\cal F}(B')$,
 we denote by $Q_h^{hyb, E}(F)$ the quadratic form on ${\cal D}\times{\cal D}$ 
such that, for all $f$ and $g$ in ${\cal D}$, we have:
\begin{equation}
 Q_h^{hyb, E}(F)( f,g) =  Q_h^{Weyl}(\widetilde H_{E^{\perp},h/2}F)( f,g).
\label{(2.8)}  
\end{equation}
\end{defi}

The following result underlines the relationship between
hybrid operators when associated with the same symbol, but with two
different subspaces, one being included into the other one.

 Let $F$ be a Borel function bounded on $B^2$. Let $E_1$ and $E_2$ be
 in ${\cal F}(B')$,  such that $E_1 \subset E_2$. Let $S$ be
the orthogonal complement of $E_1$ in $E_2$ (for the scalar product
of $H$). Then, we have:
\begin{equation}
 Q_h^{hyb, E_1}(F)  =   Q_h^{hyb, E_2}(\widetilde H_{S,h/2}F)\label{(2.9)}
\end{equation}
and we have:
\begin{equation}
 \widetilde H_{S,h/2}F (X) =(\pi h)^{-{\rm dim}(S)} \int _{S^2} 
 e^{\frac{|Y|^2}{ h}} F(X +Y ) d\lambda (Y)  \hskip 2cm X \in B^2 .\label{(2.10)}
\end{equation}
If $\{ e_1,...,e_n \}$ is an orthonormal basis of the orthogonal
complement $S$ of $E_1$ in $E_2$, denoting by $D_j$  the subspace of $H$
spanned by $e_j$, we have:
\begin{equation}
\widetilde H_{S,h/2}F = \left [\prod_{j\leq n} \widetilde H_{D_j,h/2} \right ]
 F \label{(2.11)}  
\end{equation}
where $\widetilde H_{D_j, h/2}$ is defined in (\ref{(2.6)}).

\bigskip

The proof of  (\ref{(2.9)}) and (\ref{(2.10)}) relies on  the Definition 2.2:
%---
$$Q_h^{hyb, E_j}(F) = Q_h^{Weyl }(\widetilde H_{E_j^{\perp},h/2}F) \hskip
2cm j=1, 2$$
%---
and on the following equality:
%---
$$\widetilde H_{E_1^{\perp},h/2} = \widetilde H_{E_2^{\perp},h/2} \widetilde H_{S,h/2}.$$
%---
%%%%%%%%%%%%%%%%%%%%%%%%%%%%%%%%%%%%%%%%%%%%%%%%%%%%%%%%%%%%%%%%%
%%%%%%%%%%%%%%%%%%%%%%%%%%%%%%%%%%%%%%%%%%%%%%%%%%%%%%%%%%%%%%%%%
\section{Plan of the proof of Theorem \ref{1.4} }\label{3}

Let $(e_j)_{(j\in \Gamma)}$ be a Hilbert basis of $H$, as in Definition \ref{1.4}.
For every $j$ in  $\Gamma$, we denote  by $D_j$ the subspace of $H$
spanned by  $e_j$. For every positive $t$, let $H_{D_j,t} $ and 
$\widetilde H_{D_j,t} $  be the operators defined in (\ref{(2.6)}),
the first one acting on functions on $H^2$, the second one, on functions
on $B^2$. The integration domain in  (\ref{(2.6)}) is $D_j^2$,
which is the subspace of $H^2$ spanned by  $(e_j, 0)$ and $(0, e_j)$.
 For every finite subset  $I$ of  $\Gamma$, set:
\begin{equation}
 \widetilde T_{I, h} = \prod _{j\in I} (I - \widetilde H_{D_j, h/2}), \hskip 2cm 
\widetilde S_{I, h}= \prod_{j\in I}  \widetilde H_{D_j, h/2}. \label{(3.1)}
\end{equation}
These operators act in the space of bounded Borel functions on $B^2$. We
 denote by $ T_{I, h}$ the operator defined as in (\ref{(3.1)}), but acting in
 the space of bounded continuous functions on $H^2$.

\bigskip

 For every finite subset  $I$ of  $\Gamma$, let $E(I)$ be the subspace of $H$
spanned by the $e_j$ ($j\in I$).
\begin{equation}
 E(I) = {\rm Vect } \ \{ ( e_j),\qquad j\in \Gamma \}. \label{(3.2)}
\end{equation}

\bigskip

The main Theorem   \ref{1.4} is a consequence of the Propositions
\ref{3.1} et \ref{3.3} below, which will be proved respectively in Sections
\ref{5} and \ref{6}.

\begin{prop}\label{3.1}
  Let $ F$ be in $ S_2(M, \varepsilon)$. We assume that $F$ has a stochastic 
extension  $\widetilde F$  for the measure $\mu_{B^2,h}$.
 For every finite subset $I$ in $\Gamma$,  for all $h$ in $(0, 1]$, there
 exists a bounded operator, which will be denoted by
  $Op_h ^{hyb, E(I )} (\widetilde T_{I, h} \widetilde F)$
 such that, for all $f$ and $g$ in ${\cal D}$, with the notations
 (\ref{(2.8)}), (\ref{(3.1)}) and (\ref{(3.2)}):
\begin{equation}
Q_h ^{hyb, E(I )} (\widetilde T_{I, h} \widetilde F)(f,g)  =
  < Op_h ^{hyb, E(I )} (\widetilde T_{I, h} \widetilde F)f,g>. \label{(3.3)}
\end{equation}
 Moreover its norm satisfies:
\begin{equation}\Vert Op_h ^{hyb, E(I )} (\widetilde T_{I, h} \widetilde F) 
\Vert_{{\cal L}(L^2(B,\mu_{B,h/2}) )} \leq M (81 \pi  hS_{\varepsilon}) ^{|I|}
\prod_{j\in I} \varepsilon_j^2. \label{(3.4)}
\end{equation}
\end{prop}

If one admits Proposition \ref{3.1} (which is a consequence of the
combined Propositions \ref{5.2} and \ref{5.3}),  one can prove 
the following result.

\begin{prop}\label{3.2}
Let $ F$ be a function  in $ S_2(M, \varepsilon)$, where the family 
$(\varepsilon_j)_{(j\in \Gamma)}$ is square summable. Set $h>0$. We assume that
$F$ has a stochastic extension  $\widetilde F$ for the measure $\mu_{B^2,h}$.
Then, for every increasing
sequence $(\Lambda_n)$ of finite subsets in $\Gamma$, whose union is
$\Gamma$,  there exists a sequence of operators, denoted by  
$(Op_h^{hyb,E(\Lambda_n)}(\widetilde F))$, such that, with the notations 
(\ref{(2.8)}) and (\ref{(3.2)}), for all $f,g$ in ${\cal D}$,:
\begin{equation}
Q_h^{hyb, E(\Lambda_n)}(\widetilde F)(f,g)   =
<  (Op_h^{hyb, E(\Lambda_n)}(\widetilde F)) f, g> \label{(3.5)}
\end{equation}
Moreover, the sequence of operators
 $(Op_h^{hyb,E(\Lambda_n)}(\widetilde F))_{(n\geq 1)}$ is a Cauchy sequence in
${\cal L}(L^2(B,\mu_{B,h/2}))$. Its limit, denoted by
 $Op_h^{Weyl}(F)$, satisfies  (\ref{(1.18)}).
\end{prop}

{\it Proof.}
We have, for every  continuous and bounded function $G$  on $B^2$,
 for any finite subset $\Lambda $ in $\Gamma$:
\begin{equation} G = \sum _{I \subseteq \Lambda } \widetilde T_{I,h}
 \widetilde S_{\Lambda \setminus I, h} G \label{(3.6)} 
\end{equation}

The sum runs over all the subsets in $\Lambda$, including the empty
set and $\Lambda$ itself. As a consequence, the equality (\ref{(3.5)}) will
 be satisfied if we set:

\begin{equation}
  Q_h ^{hyb, E(\Lambda)}  (\widetilde{F}) = \sum _{I \subseteq \Lambda } 
  Q_h ^{hyb, E(\Lambda)} (\widetilde T_{I,h} 
\widetilde S_{\Lambda \setminus I, h} \widetilde{F}). \label{(3.7)}
\end{equation}

From  (\ref{(2.9)}) and (\ref{(2.10)}), applied with subspaces $E(I)$ and
$E(\Lambda)$:

\begin{equation}
 Q_h ^{hyb, E(\Lambda)} (\widetilde{F}) = 
\sum _{I \subseteq \Lambda }   Q_h ^{hyb, E(I)}
 (\widetilde T_{I,h}  \widetilde{F}). \label{(3.8)}
\end{equation}

The term corresponding to $I= \emptyset$ is the anti-Wick  quadratic form
associated with $\widetilde{F}$ and adapted to the Gaussian measure.
 Therefore, according to Proposition \ref{3.1}, there exists 
a bounded operator, denoted by  $Op_h ^{hyb, E(\Lambda)} $
 such that \ref{(3.5)} is satisfied. If
$(\Lambda _n )$ is an increasing sequence of finite subsets of
$\Gamma$, we then have, if $m< n$:
\begin{equation} 
Op_h ^{hyb, E(\Lambda_n)} (\widetilde{F}) - 
Op_h ^{hyb, E(\Lambda_m)} (\widetilde{F})=
\sum _{I\in {\cal P}(m, n)}  Op_h ^{hyb, E(I)}
 (\widetilde T_{I,h}  \widetilde{F}) \label{(3.9)}
\end{equation}
 where ${\cal P}(m, n)$ is the family of subsets  $I$ in $\Gamma $,
included in  $\Lambda_n$, but with at least one element not
belonging to $\Lambda_m$.  From (\ref{(3.9)}) and from Proposition 
 \ref{3.1}, we have, when $m<n$:
$$
\begin{array}{lll}
\displaystyle
 \Vert Op_h ^{hyb, E(\Lambda_n)} (\widetilde F)-Op_h^{hyb, E(\Lambda_m)} (\widetilde F)
\Vert  _{{\cal L}( L^2(B,\mu_{B,h/2}) )}&
\displaystyle \leq  \sum _{I\in {\cal
P}(m, n)} \Vert  Op_h ^{hyb, E(I)}
 (\widetilde T_{I,h}  \widetilde F)\Vert  _{{\cal L}( L^2(B,\mu_{B,h/2}) )} \\&
\displaystyle  \leq M \sum _{I\in {\cal P}(m, n)}(81 \pi hS_{\varepsilon} )^{|I|}
 \prod _{j\in I}\varepsilon_j^2.
\end{array}
$$
As a consequence, if $m<n$:
%----
$$\Vert Op_h^{hyb, E(\Lambda _n )}(\widetilde F) - Op_h^{hyb, E(\Lambda_m)}
(\widetilde F) \Vert _{{\cal L}( L^2(B,\mu_{B,h/2}) )}\leq M 81
\pi h S_{\varepsilon} \left [ \sum _{j\in \Lambda _n \setminus
\Lambda_m} \varepsilon_j^2\right ] \prod _{k\in \Lambda_n} (1 + 81
\pi hS_{\varepsilon} \varepsilon _k^2). $$
%--
If the family $(\varepsilon_j^2)_{(j\in \Gamma)}$ is summable, the
above right hand-side product stays bounded independently of $n$,
whereas the sum tends to $0$ when $m\rightarrow +\infty$. As a
consequence, the sequence $(Op_h ^{hyb, E(\Lambda_n)} (\widetilde
F))$ is a Cauchy sequence in ${\cal L}( L^2(B,\mu_{B,h/2}) )$.
Likewise:
%-----
$$\Vert Op_h^{hyb, E(\Lambda _n )}(\widetilde F)  \Vert _{{\cal L}(
L^2(B,\mu_{B,h/2}) )}\leq M  \prod _{k\in \Lambda_n} (1 + 81 \pi
h S_{\varepsilon} \varepsilon _k^2). $$
%--
Therefore, the limit of this sequence of operators, denoted by
$Op_h^{Weyl}(F)$, verifies  (\ref{(1.18)}).\hfill $\square$ 

\bigskip

One could think that the operator $Op_h^{Weyl}(F)$ depends on the sequence 
$(\Lambda _n)$, but the following proposition proves that it is not the case.
It will be proved in Section \ref{6}.

\begin{prop}\label{3.3}
 Let $ F$ belong to $ S_2(M,\varepsilon)$, where the sequence 
$(\varepsilon_j)_{(j\in \Gamma )}$
is square summable. Set $h>0$. Define the function $\widetilde F$ on $B^2$  as
 the stochastic extension of $ F$ both for the measure $\mu_{B^2,h}$ and for
 the measure $\mu _{B^2,h/2}$.
 Let $(\Lambda_n)$ be an increasing sequence  of finite subsets
  of $\Gamma$, whose union is $\Gamma$. Then, we have, for every $f$
  and $g$ in ${\cal D}$, setting $E_n = E(\Lambda_n)$
\begin{equation}
 \lim _{n\rightarrow + \infty} Q_h^{Weyl}(\widetilde H_{E_n^{\perp} ,h/2}\widetilde F)
 (f,g) =  Q_h^{Weyl}(\widetilde F) (f,g)   \label{(3.10)} 
\end{equation}
\end{prop}

One can notice that, for a  given function $F$ in $ S_2(M,\varepsilon)$, 
the stochastic extension (for the measure $\mu_{B^2,s}$)
 $\widetilde F$ is unique, but only up to a $\mu_{B^2,s}$- negligible set.
Nevertheless, if $\widetilde F$ and $\widetilde G$ are two stochastic
 extensions of the same function $F$, for the measures
 $\mu_{B^2,h}$ and  $\mu _{B^2,h/2}$, one can check that 
 $\widetilde H_{E_n^{\perp} ,h/2}\widetilde F$ and $\widetilde H_{E_n^{\perp} ,
h/2}\widetilde G$ are equal almost everywhere for  the measure
  $\mu _{B^2,h/2}$.

\bigskip

{\it End of the proof of Theorem \ref{1.4}.}
Once Proposition \ref{3.2} has been established, it only remains to prove
 (\ref{(1.17)}).
 Now for all $f$ and $g$ in ${\cal D}$ one has, setting $E_n = E(\Lambda_n)$:
%---
$$ < Op_h^{Weyl}(F) f,g> = \lim _{n\rightarrow + \infty}<(Op_h^{hyb, E_n } 
(\widetilde F)) f,g> = \lim _{n\rightarrow + \infty}  Q_h^{Weyl}
(\widetilde H_{E_n^{\perp} ,h/2}\widetilde F)( f,g).  $$
%---
Hence the equality (\ref{(1.17)}) of Theorem \ref{1.4} comes from the above
equality and from (\ref{(3.10)}).

%%%%%%%%%%%%%%%%%%%%%%%%%%%%%%%%%%%%%%%%%%%%%%%%%%%%%%%%%%%%%%%%%
%%%%%%%%%%%%%%%%%%%%%%%%%%%%%%%%%%%%%%%%%%%%%%%%%%%%%%%%%%%%%%%%%
\section{Some useful operators in Wiener spaces}\label{4}

We first recall the precise definition of Wiener spaces as well as some
of their properties. Then we  adapt some classical notions for Wiener spaces :
coherent states, Segal-Bargmann transformation. Next we  give properties
of the Gaussian Wigner function of Section \ref{1}. This will allow us to 
write Definition  \ref{1.2} of the Weyl operators and Definition \ref{2.2} 
of the hybrid operators in a way more suitable for norm estimates.
This will yield Proposition \ref{4.11}, which will be used in Section \ref{5}
to prove Proposition  \ref{3.1}.

%%%%%%%%%%%%%%%%%%%%%%%%%%%%%%%%%%%%%%%%%%%%%%%%%%%%%%%%%%%%%%%%%
%%%%%%%%%%%%%%%%%%%%%%%%%%%%%%%%%%%%%%%%%%%%%%%%%%%%%%%%%%%%%%%%%
\subsection{Abstract Wiener spaces}\label{4.A}

If $H$ is a real separable infinite-dimensional Hilbert space, it is impossible
to define on its Borel $\sigma$-algebra a measure $\mu_{H, h}$
such that  (\ref{(1.7)}) holds with $B$ instead of $H$.

\bigskip

Nevertheless, one can define a promeasure (or cylindrical
probability measure in the sense of \cite{K-R}) $\mu_{H, h}$ on the cylinder sets
of $H$. For every $E$ in ${\cal F}(H)$ and every positive $h$,
one can define a Gaussian measure $\mu_{E,h}$ on $E$ by (\ref{(1.5)}),
where $\lambda_E$ is the Lebesgue measure on $E$ (normalized in a natural way)
and    $|\cdot|$ is the norm of  $H$.
A cylinder set of $H$ is any set of the form $C = \pi_E^{-1}(\Omega)$, where
 $E\in {\cal F}(H)$,  $\pi_E: H \rightarrow E$ is the orthogonal projection 
and $\Omega$  is a Borel set of $E$. If $C$ is such a cylinder set, 
one sets $\mu_{H,h}(C)= \mu_{E,h}(\Omega)$. In other words, for every
Borel set $\Omega$ of $E$:
\begin{equation}
\mu _{H,h} (  \pi_E^{-1}(\Omega))  = (2\pi h)^{-{\rm dim} (E)/2}
 \int _{\Omega}   e^{-\frac{|y|^2}{ 2h}} d\lambda (y). \label{(4.1)} 
\end{equation}
One defines this way an additive set function on the cylinder sets of
$H$, but if $H$ is infinite-dimensional, this function is not $\sigma$-additive 
and  $\mu_{H,h}$ does not extend as a measure on the $\sigma$-algebra generated
by the cylinder sets (which is the Borel $\sigma$-algebra).

\bigskip

If the Hilbert space $H$ is included into a Banach space $H$, the canonical
injection being continuous and having a dense range, so that (\ref{(1.4)})
is satisfied, one can define the cylinder sets of $B$ as the sets
 $C=P_E^{-1}(\Omega)$, where $E\in {\cal F}(B')$, where $P_E: B \rightarrow E$
is defined by (\ref{(1.6)}) and where $\Omega$ is a Borel set of $E$. One
then defines an additive set function  $\mu_{B,h}$ on the cylinder sets of
$B$ as in  (\ref{(4.1)}):
%------
$$ \mu _{B,h} (  P_E^{-1}(\Omega))  = (2\pi h)^{-{\rm dim} (E)/2}
 \int _{\Omega}   e^{-\frac{|y|^2}{ 2h}} d\lambda (y).
     $$
But this time, if $B$ is well chosen, the additive sets function $ \mu _{B,h}$
extends as a measure on the  Borel $\sigma$-algebra of $B$. The following definition specifies the conditions which $B$ must fulfill.

\begin{defi}\label{4.1}\cite{G-2, G-3,K}
 An abstract Wiener space is a triple  $(i, H,B)$ where $H$ is a real
separable Hilbert space, $B$ a Banach space and $i$ a continuous injection
from $H$ into $B$, such that $i(H)$ is dense in $B$ and satisfying, moreover,
the following condition. For all positive $\varepsilon $ and $h$, there exists
a subspace  $F$ in ${\cal F}(H)$ such that, for all  $E$ in ${\cal F}(H)$,
 orthogonal to $F$,
$$ \mu _{H,h} \left ( \left \{ x\in H,\ \ \Vert i(\pi_E(x)) \Vert _B >
\varepsilon \right \} \right ) < \varepsilon. $$
%---
The norm on $B$ is said to be ``measurable''.
\end{defi}

If  $(i, H,B)$ is an abstract Wiener space (in other words if 
the norm of $B$ is measurable in the sense above) and if (\ref{(1.4)}) holds
 one proves that, for every
positive  $h$, the additive set function  $\mu_{B,h}$, defined on the
 cylinder set functions of $B$, extends as a measure on the  Borel 
$\sigma$-algebra of $B$  and has the following property. For
every finite system  $\{ u_1,...,u_n\}$ of $B'$, which is orthonormal with
 respect to the scalar product of  $H$, the functions $x \rightarrow u_j
(x)$ (defined on  $B$) are independent Gaussian random variables and the system
 $\{ u_1, ...u_n \}$ has the normal distribution  $ \mu _{\R^n,h}$. 
See \cite{G-2}, \cite{G-3} and \cite{K} (consequence of the Theorems 4.1 and
 4.2, Chapter 1).  For every $E$ in  ${\cal F}(B')$ and every 
$\varphi $ in $L^1(E, \mu _{E,h})$, the equality  (\ref{(1.7)}) is satisfied,
 according to the transfer Theorem.

\bigskip

Let us recall a classical example of the explicit computation of an integral,
where $a$ is in the complexified space  $H_{\bf C}$:

\begin{equation}
 \int _B e^{ \ell_a (x)}  d\mu _{B,h} (x) = e^{ h \frac{a^2}{ 2}}.
\label{(4.2)} 
\end{equation}
One has set $a^2 = |u|^2 - |v|^2 + 2i u\cdot v $ if $a = u+i v$, with
$u$ and $v$ in $H$. Let us recall, too, that for all $a$ in $H$and for all 
 $p\geq 1$:
\begin{equation} \int _B |\ell_a (x)|^p   d\mu _{B,h} (x) = 
 \frac{(2h)^{p/2}}{ \sqrt { \pi} } |a|^p \
\Gamma \left (\frac{p+1}{ 2} \right ). \label{(4.3)}
\end{equation}
One sees, too, that for all $a$ and $B$ in $H$,
 \begin{equation}
\int_B e^{\ell_b(u)} |\ell _a (u)|^p d\mu_{B,h} (u) =
e^{h\frac{|b|^2}{ 2}} \int_{\R} |\sqrt {h} |a|v + h a\cdot b|^p d\mu_{\R,1}(v).
\label{(4.3bis)}
 \end{equation}
The following proposition allows to deal with translations by a vector
$a$ belonging to $H$. There is no such result for a  translation by a vector
$a$ belonging to $B$.
\begin{prop}\label{4.2}\cite{K}, p 113,114
Let $(i,H,B)$ be an abstract Wiener space and  $\mu _{B,h}$ its measure.
For all $g\in L^1(B, \mu_{B,h}) $, one has, for all $a$
in $H$:
\begin{equation}
\int_B g(x) d\mu_{B,h} (x)= e^{-\frac{1}{2h}|a|^2}
\int_B g(x+a) e^{ - \frac{1}{h}\ell_a (x) }  d\mu _{B,h}(x).\label{(4.4)} 
\end{equation}
 For every $a\in H$ and every $f$ in $ L^2(B, \mu _{B,h}) $, one has
\begin{equation}
\int_B |f(x)|^2 d\mu_{B,h} (x)= e^{-\frac{1}{2h}|a|^2}
\int_B |f(x+a)|^2 e^{ - \frac{1}{h}\ell_a (x) }  d\mu _{B,h}(x).\label{(4.5)}  
\end{equation}
\end{prop}

Let us recall the theorem of Wick :
\begin{theo}\label{4.3} {\bf Wick}
 Let $u_1, ... u_{2p}$ be vectors of $H$ ($p\geq 1$). Let  $h>0$. Then
one has
\begin{equation}
 \int _B \ell _{u_1} (x) ... \ell _{u_{2p}} (x) d\mu _{B,h} (x) = h^p
\sum _{(\varphi, \psi)\in S_p} \prod _{j=1}^p < u_{\varphi (j)}, u_{\psi (j)}>
\label{(4.6)}
\end{equation}
where $S_p$ is the set of all couples  $(\varphi, \psi)$ of injections from
$\{ 1,...,p\}$ into $\{1,...., 2p\}$ such that:
\begin{enumerate}
\item
 For all  $j\leq p$,  $\varphi (j) < \psi (j)$.
\item
 The sequence $(\varphi (j) ) _{(1\leq j \leq k)}$ is an increasing sequence.
\end{enumerate}
\end{theo}

One deduces from (\ref{(4.2)}) the following inequalities, which hold for all
$a$ and $b$ in the complexified of $H$ :
\begin{equation}
\begin{array}{lll}
\displaystyle
 \int _B \left| e^{\ell _a (x) } - e^{\ell _b (x) } \right |^2 d\mu _{B,h} (x)
\leq  4h|a -b |\ (|a| + |b|)\ e^{ 2h \max ( |{\rm Re}\ a |^2 ,|{\rm Re} \  b  |^2) }\\ 
\displaystyle  \int_B |e^{\ell_a(x)}-e^{\ell_b(x)}|^2 \ d\mu_{B,h}(x) \leq
e^{2h \max(|{\rm Re}(a)|,|{\rm Re}(b)|)^2} h|a-b|^2 
(1+4h \max(|{\rm Re}(a)|,|{\rm Re}(b)|)^2)
\end{array}  \label{(4.7)} 
\end{equation}

%%%%%%%%%%%%%%%%%%%%%%%%%%%%%%%%%%%%%%%%%%%%%%%%%%%%%%%%%%%%%%%%%
%%%%%%%%%%%%%%%%%%%%%%%%%%%%%%%%%%%%%%%%%%%%%%%%%%%%%%%%%%%%%%%%%
\subsection{Stochastic extensions}\label{4.B}

In order to define the stochastic extension of a function
  $f:H\rightarrow {\bf C}$, we first define the extension of the
orthogonal projection  $\pi_E: H \rightarrow E$, where $E\in{\cal F}(H)$.

\bigskip

 With each $E$ in ${\cal F}(H)$, one can associate a map
 $\widetilde \pi_E: B\rightarrow E$ defined almost everywhere by
\begin{equation}
 \widetilde \pi_E(x) = \sum _{j= 1}^{{\rm dim}(E)} \ell_{u_j}(x) u_j,
 \label{(4.B.1)}
\end{equation}
where the $u_j$ $(1\leq j \leq {\rm dim}(E))$ form an orthonormal
basis of $E$ and $\ell_{u_j}$ is defined in Proposition \ref{1.1}. This map does
 not depend on the choice of the orthonormal basis and, when $E\subset B'$, it
 coincides with the map $P_E$ already defined by (\ref{(1.6)}).

\bigskip

Remark that (\ref{(1.7)}) still holds for the subspaces $E$ in ${\cal F}(H)$ 
(and not necessarily in  ${\cal F}(B'))$, provided  $P_E$ is replaced by
$\widetilde \pi_E$. See Lemma 4.7 in \cite{K}.

\bigskip

Notice that, for every subspace $E$ in  ${\cal F}(H)$:
\begin{equation}
 a \cdot (\widetilde \pi _E (x) ) = \ell _{\pi_E(a)} (x)
 \hskip 2cm a\in H \ \ \ \ \
a.e. x\in B \label{(4.B.2)}
\end{equation}
\bigskip

The following notion has been introduced by  Gross \cite{G-1},
 who gives conditions for the existence of the extension. Other conditions
 or examples will be found in Section \ref{8} or in \cite{A-J-N-1}.

\begin{defi}\label{4.4} \cite{G-1,G-2,G-3} \cite{RA} \cite{K}
 Let $(i, H,B)$  be an abstract Wiener space satisfying (\ref{(1.4)}).
 (The inclusion $i$ will be omitted). Let $h$ be a positive real number.
\begin{enumerate}
\item
 A  Borel function $ f$, defined on $H$, is said to admit a stochastic
 extension $\widetilde f$ with respect to the measure  $\mu_{B,h}$ if, for
every increasing sequence $(E_n)$ in ${\cal F}(H)$, whose union is dense in $H$,
the sequence of functions $ f \circ \widetilde \pi _{E_n}$ (where 
 $\widetilde \pi _{E_n}$ is defined by  (\ref{(4.B.1)})) converges in
probability with respect to the measure $\mu_{B,h}$  to
$\widetilde  f$.  In other words, if, for every  $\delta>0$,
\begin{equation}
 \lim _{n\rightarrow + \infty } \mu_{B,h} \left(  \left \{
x \in B,\ \ \ \  | f \circ \widetilde \pi_{E_n } (x) -  \widetilde f (x) |
> \delta \right \} \right ) = 0. \label{(4.B.3)}
\end{equation}
\item 
A function $f$ is said to admit a stochastic extension 
 $\widetilde f\in L^p (B,\mu _{B,h})$ in the sense of $L^p (B,\mu _{B,h})$ 
($1\leq p< \infty$) if, for every increasing sequence $(E_n)$ in ${\cal F}(H)$,
 whose union is dense in $H$, the functions  $f \circ\widetilde \pi_{E_n}$
are in  $L^p (B,\mu_{B ,h})$ and if the sequence  $f \circ \widetilde \pi_{E_n}$
converges in $L^p (B,\mu _{B,h})$ to $\widetilde f$.
\end{enumerate}
 One defines likewise 
the stochastic extension of a function on $H^2$ to a function on  $B^2$.
\end{defi}

If  $\widetilde f$ is the stochastic extension  of a function 
$f : H\rightarrow H$, one cannot say that $f$ is the restriction of
 $\widetilde f$ to $H$. Since $H$ is negligible (see Kuo \cite{K}), this is 
irrelevant. For every $a$ in $H$ one sees that the application
 $u \rightarrow u\cdot a$, defined on  $H$, admits a stochastic extension
which is the function  $\ell_a$. This is a consequence of the equality
 (\ref{(4.B.2)}). In Definition \ref{4.4}, the functions can take their
values in a Banach space. Hence one can say that the application
 $\widetilde \pi_E$ of  (\ref{(4.B.1)}) is the  stochastic extension of the
 orthogonal projection  $\pi_E : H \rightarrow E$. One will find in Section 
\ref{8} examples of functions admitting stochastic  extensions. In particular,
if a function $f$ is bounded on $H$ and uniformly continuous with respect to
the restriction to $H$ of the norm of $B$, then it admits a stochastic extension
$\widetilde f$, which coincides with its density extension (see Kuo \cite{K},
Chapter 1, Theorem 6.3).
 
\bigskip

 If a Borel function $f$ is bounded on $H$, does not depend on $h$ and admits,
 for every positive $h$, a stochastic  extension with respect to $\mu_{B,h}$, 
 this extension may depend on $h$. It is not the case in the situation of
 Proposition \ref{8.5}.
 In the other cases, we may consider that the stochastic
 extension is independent of $h$ if $h$ varies in a  countable subset
 $Q$ of  $(0, +\infty)$.

%%%%%%%%%%%%%%%%%%%%%%%%%%%%%%%%%%%%%%%%%%%%%%%%%%%%%%%%%%%%%%%%%
%%%%%%%%%%%%%%%%%%%%%%%%%%%%%%%%%%%%%%%%%%%%%%%%%%%%%%%%%%%%%%%%%
\subsection{The Heat operator (continued)}\label{4.C}

In the rest of this work, $(i, H, B)$ represents an  abstract Wiener space
satisfying  (\ref{(1.4)}) and the injection $i$ will be omitted. We complete
the investigation begun in Sections  \ref{2.A} and \ref{2.C}.

\bigskip

Let $E$ be in  ${\cal F}(B')$, let  $E^{\bot}$ be its orthogonal space, defined 
in (\ref{(2.XX)}) and let $t>0$. For every  Borel function  $F$ bounded on
 $B^2$, let $\widetilde H_{E^{\bot},t}F$  be the function defined in (\ref{(2.7)}):
$$
(\widetilde H_{E^{\bot},t} F)(X)= \int_{(E^{\bot})^2} F(X+Y)\ d\mu_{(E^{\bot})^2,t}(Y).
$$
For all real numbers  $h_1>0$ and $h_2>0$, one defines a probability measure
 on   $B^2$, using Proposition \ref{2.1},  by:
\begin{equation}
 \nu _{B^2,E^2, h_1, h_2 } = \mu_{E^2,h_1}\otimes \mu_{(E^{\bot})^2,h_2} .
\label{(4.C.2)}
\end{equation}

\begin{prop}\label{4.5}
 Let  $h_1, h_2$ and $t$ be positive real numbers. The operator 
 $\widetilde H_{E^{\bot},t}$ is bounded and its norm is at most $1$ from
$L^p(B^2, \nu _{B^2,E^2, h_1, h_2 +t })$ into
$L^p(B^2, \nu _{B^2,E^2, h_1, h_2  })$  for every $p\in [1,+\infty]$. Let
%------
\begin{equation}
(\widetilde M_{E^{\bot},t,h_1,h_2} G)(X)= \int_{(E^{\bot})^2}
G(X_E, Y_{E^{\bot}} + \frac{h_2}{ t+h_2} X_{E^{\bot}} ) \
d\mu_{(E^{\bot})^2,\frac{th_2}{ t+h_2  }}( Y_{E^{\bot}}).  \label{(4.C.3)}
\end{equation}
The operator  $\widetilde M_{E^{\bot},t,h_1,h_2}$ is bounded and its norm is at
 most $1$ from  $L^p(B^2,\nu _{B^2,E^2, h_1, h_2  } )$
into  $L^p(B^2,\nu _{B^2,E^2, h_1, h_2 +t } )$ for every $p\in [1,+\infty]$.
Moreover, for all functions  $F\in L^p(B^2,\nu _{B^2,E^2, h_1, h_2 +t })$
and  $G \in L^q(B^2,\nu _{B^2,E^2, h_1, h_2  })$,
 with $p$ and $q$ conjugate exponents, one has:
\begin{equation}
\int_{B^2}  \widetilde H_{E^{\bot},t} F(X) G(X) \ d \nu _{B^2,E^2, h_1, h_2 }(X)=
\displaystyle  \int_{B^2} F(X)
 \widetilde M_{E^{\bot},t,h_1,h_2}  G(X) \ d \nu _{B^2,E^2, h_1, h_2+t }(X).
\label{(4.C.4)} 
\end{equation}
\end{prop}

{\it  Proof. } Suppose first that $p$ is finite. Let $F$ be a bounded
cylindrical Borel function : there exist a subspace  $D$  in ${\cal F } (B')$
and a Borel function  $\widehat{F}$, bounded on $D^2$, such that, for every
$X\in B^2$, $\ F(X)= \widehat{F}(P_D(X)).$ One can assume, without loss of
 generality, that $E\subset D$. Let $(e_1,\dots,e_s,e_{s+1},\dots,e_d)$ be a
basis of $D$, orthonormal with respect to the scalar product of $H$, 
 $(e_1,\dots,e_s)$ being a basis of $E$. Let us use the following notation :
%---
$$\widehat{F}(\sum_{k=1}^d x_ke_k,\sum_{k=1}^d \xi_ke_k)
=f(\widetilde{x},x',\widetilde{\xi},\xi' )$$
%---
with $\widetilde{x}=(x_1,\dots,x_s)$, $x'=(x_{s+1},\dots,x_d)$ and the same
convention for $\xi$. Let  $p\in [1,+\infty[$. 
According to Jensen's inequality, 
$$
|| \widetilde H_{E^{\bot},t}F||^p_{L^p(B^2,\nu_{B^2, E^2,h_1, h_2})}
\leq \int_{B^2}\int_{(E^{\bot})^2}
|  \widehat{F}(P_D(X+Y) )|^p d\mu_{(E^{\bot})^2,t}(Y)
d \mu_{E^2,h_1}(X_E) d \mu_{(E^{\bot})^2,h_2 }(X_{E^{\bot}}),
$$
therefore, by (\ref{(1.7)}),
$$\begin{array}{lll}
\displaystyle
|| \widetilde H_{E^{\bot},t} F||^p_{L^p(B^2,\nu_{B^2, E^2,h_1, h_2})} \leq
 (2\pi t )^{-d+s} (2\pi h_1 )^{-s}(2\pi h_2 )^{-d+s}\\
\displaystyle  \int_{\R^{4d-2s}}
\left|f(\widetilde{x}, x'+y',\widetilde{\xi},\xi'+\eta')\right|^p \
e^{- \frac{1}{ 2t} (| y'|^2+ |\eta'|^2)}
e^{- \frac{1}{ 2h_1}(|\widetilde{x}|^2+|\widetilde{\xi}|^2)}e^{- \frac{1}{ 2h_2}
(|x'|^2+|\xi'|^2)} d\lambda(\widetilde{x},x',y',\widetilde{\xi},\xi',\eta'). 
\end{array}
$$
The transition from  $x'+y',\xi'+\eta'$ to  $z',\zeta'$ and a classical
computation on the convolution of Gaussian functions give
$$
\begin{array}{lll}
\displaystyle
 || \widetilde H_{E^{\bot},t}F||^p_{L^p(B^2,\nu_{B^2, E^2,h_1, h_2})}\\
\leq \displaystyle (2\pi h_1)^{-s}(2\pi(h_2+t))^{s-d}
\int_{\R^{2d}} |f(z,\zeta)|^p e^{- \frac{1}{ 2(h_2+t)}( | z'|^2+|\zeta'|^2)}
e^{- \frac{1}{ 2h_1} (| \widetilde{z}|^2+| \widetilde{\zeta}|^2) }
d\lambda(\widetilde{z},z',\widetilde{\zeta},\zeta' )
\\\displaystyle
 = \int_{B^2} |F(X)|^p  d \mu_{E^2,h_1}(X_E)d \mu_{(E^{\bot})^2,h_2 +t}(X_{E^{\bot}}).
\end{array}
$$
One extends this result, using the density of the cylindrical functions
in $L^p(B^2,\nu_{B^2, E^2,h_1, h_2+ t})$.

If $F\in L^{\infty}(B^2,\nu_{B^2, E^2,h_1, h_2+ t}) $, let  $G$ be the function equal
to $F$ on the set
$$
\{ X\in B^2 : |F(X)|\leq ||F||_{L^{\infty}(B^2,\nu_{B^2, E^2,h_1, h_2+ t}) } \}
$$
and to $0$ elsewhere. The function  $F-G$ is equal to $0$
 $\nu_{B^2, E^2,h_1, h_2+ t}$-almost everywhere, thence  in
$L^1(B^2,\nu_{B^2, E^2,h_1, h_2+ t}))$, which  implies that
 $\widetilde H_{E^{\bot},t} F= \widetilde H_{E^{\bot},t} G$
$\nu_{B^2, E^2,h_1, h_2}$ almost everywhere.
Now, for all $X$ in $B^2$, $ \widetilde H_{E^{\bot},t} G(X)\leq
 ||F||_{L^{\infty}(B,\nu_{B^2, E^2,h_1, h_2+ t})} $, then
$ \widetilde H_{E^{\bot},t} F(X)\leq ||F||_{L^{\infty}(B,\nu_{B^2, E^2,h_1, h_2+ t})}$
for $\nu_{B^2, E^2,h_1, h_2}$- almost every $X$.

The result about $\widetilde M_{E^{\bot},t,h_1,h_2}$ can be proved in a similar way.
One can also introduce the partial dilation
$\Delta_{ E^{\bot}, \lambda} $, defined by
$$
\Delta_{ E^{\bot}, \lambda}G(X_E,X_{E^{\bot}}) = G(X_E, \lambda X_{E^{\bot}}),
$$
which is isometrical from $L^p(B^2,\nu_{B^2, E^2,h_1, h_2})$ in 
$L^p(B^2,\nu_{B^2, E^2,h_1, h_2/\lambda ^2})$ and remark that
$$
\widetilde M_{E^{\bot},t,h_1,h_2} =  \Delta_{ E^{\bot}, \frac{h_2}{t+h_2}}
\circ \widetilde H_{E^{\bot},\frac{th_2}{ t+h_2}}.
$$
The  equality  (\ref{(4.C.4)}) is the result of computations similar
to the precedent ones.

\bigskip

When  $E=\{ 0_B\}$, this yields results for the global heat
 operator, the product measures being ``global'' on $B^2$.

Gross and Kuo gave results about these operators in the case when they act 
on functions more regular than $L^p$.

%%%%%%%%%%%%%%%%%%%%%%%%%%%%%%%%%%%%%%%%%%%%%%%%%%%%%%%%%%%%%%%%%
%%%%%%%%%%%%%%%%%%%%%%%%%%%%%%%%%%%%%%%%%%%%%%%%%%%%%%%%%%%%%%%%%
\subsection{Coherent states and  Segal Bargmann transformation }\label{4.D}

For every  $E $ in ${\cal F}(H)$, the coherent states on $E$ are functions
 ${\Psi}^{E, Leb}_{X,h}  $ defined,  for all $X =(a,b) \in E ^2$ and
for any  $h>0$, by:

\begin{equation} \Psi^{E, Leb}_{X,h} (u)  =    (\pi h)^{ -{\rm dim}(E)/4}
e^{-\frac{| u-a|^2 }{ 2h}} e^{\frac{i}{ h} u .b - \frac{i}{ 2h} a. b}
 \qquad X = (a,b) \in E^2 \qquad u\in E . \label{(4.D.1)} 
\end{equation}
\bigskip

It is well known that, for all $u$ and $v$ in  $L^2(E,\lambda)$, one has
\begin{equation}
< u,v> _{Leb} = (2 \pi h)^{-{\rm dim}(E)} \int _{E^2}  < u ,\Psi^{E, Leb}_{X ,
h}>_{Leb} < \Psi^{E, Leb}_{X,h},v>_{Leb} d\lambda (X) \label{(4.D.2)}
\end{equation}
where $< \cdot,\cdot >_{Leb}$  is the scalar product of  $L^2(E ,\lambda)$.
\bigskip

One can also define Gaussian coherent states, which are functions on the Hilbert
 space $H$. For every $X = (a,b)$ in  $H^2$, one defines a function
$\Psi_{X,h} $ on $H$, setting:
\begin{equation}
\Psi_{X,h} (u) = e^{\frac{1}{h} u\cdot (a+ib) -\frac{1}{2h}|a|^2 -
\frac{i}{2h} a\cdot b}. \label{(4.D.3)}
\end{equation}
This function is only defined for $u$ in $H$. It cannot be used as it is
in any integral, but it admits a stochastic extension in the sense of
 $L^2(B,\mu_{B,h/2})$ (see Definition \ref{4.4}). This extension is defined by
\begin{equation}
\widetilde \Psi_{X,h}(u)
 = e^{ \frac{1}{h} \ell_{a+ib}(u)-\frac{1}{2h}|a|^2 -\frac{i}{2h} a\cdot b} \quad X =
 (a,b) \in H^2 \  {\rm a.e.}  u\in B .  \label{(4.D.4)}
\end{equation}
  
The norm of the extension is equal to $1$ (see (\ref{(4.3)})). According  to
 (\ref{(4.2)}), one has, for all $U$ and $V$ in $H$
\begin{equation}
 < \widetilde \Psi_{U,h} ,\widetilde  \Psi_{V,h} >
_{L^2(B,\mu_{B,h/2})} = e^{-\frac{1}{4h} |U-V|^2 + 
\frac{i}{ 2h} \sigma (U,V)}. \label{(4.D.5)}
\end{equation}

\bigskip

For every $E$ in  ${\cal F}(H)$ and every $X$ in $E^2$ one can write:
\begin{equation} 
\Psi_{X,h}= \Big(\gamma_{E,h/2}^{-1}\Psi^{E,Leb}_{X,h}\Big)\circ P_E. \label{(4.D.6)}
\end{equation}

\bigskip

To define the Segal-Bargmann transformation one starts from its well-known
analog in a finite dimensional case (see \cite{F} or \cite{SJ})). For every
$E$ in ${\cal F}(H)$ and every function  $\varphi $ in ${\cal S}_E$, one
defines a function   $\widehat T_h\varphi $ on $E^2$ setting:
\begin{equation}
 (\widehat T_h  \varphi ) (x,\xi) = e^{-\frac{1}{4h}  (x-i\xi)^2  } \
\int _{E} \varphi (y) e^{ \frac{1}{ h} y \cdot (x - i \xi)  } d\mu
_{E, h/2} (y) \hskip 2cm  (x,\xi) \in  E^2 \label{(4.D.7)}
\end{equation}
where $(x-i\xi)^2 = |x|^2 -|\xi|^2 - 2i x\cdot \xi$. One knows (\cite{F})
 that $\widehat T_h$ is a partial isometry from  $L^2(E,\mu _{E h/2})$ in
  $L^2( E^2,\mu _{E^2, h})$ and that, for all $Z = (z,\zeta)$ in $E^2$:
\begin{equation}
 (\widehat T_h  \varphi ) (z,\zeta) = \int _{E^2} e^{\frac{1}{2h}
 (x+i \xi)\cdot  (z-i\zeta) } (\widehat T_h  \varphi ) (x,\xi)
 d\mu _{E^2,h} (x,\xi). \label{(4.D.8)} 
\end{equation}
One can see, too, that:
\begin{equation}
 \widehat T _h \varphi (X) = e^{\frac{1}{4h} |X|^2} < \gamma _{E,h/2}
\varphi ,\Psi_{X,h}^{E, Leb} > _{Leb}. \label{(4.D.9)} 
\end{equation}

\bigskip

For every function $f$ in ${\cal D}$, of the form $f =\widehat f\circ P_E$,
with $E$ in  ${\cal F} (B')$ and  $\widehat f $ in ${\cal S}_E$, one defines
a function $T_hf$ on $B^2$ by:
\begin{equation}
(T_hf)(X)=(\widehat T_h \widehat f)(P_E(X))\hskip 1cm X\in B^2. \label{(4.D.10)}
\end{equation}
If the same function $f$ has two expressions
 $f = \widehat{f_1} \circ P_{E_1} = \widehat{f_2} \circ P_{E_2}  $, with $E_1$ 
and $E_2$ in ${\cal F}(B')$, $\widehat{f_1} $ in ${\cal S}_{E_1}$ and
$\widehat{f_2} $ in ${\cal S}_{E_2}$, one checks that the right term of
 (\ref{(4.D.10)}) is the same, be it computed with $E_1$ or $E_2$, which 
justifies the notation  $T_hf$. According to the preceding, one has
\begin{equation}
\Vert T_hf \Vert_{L^2(B^2,\mu _{B^2 ,h})}=
\Vert \widehat T_h \widehat f \Vert _{L^2(E^2,\mu _{E^2 ,h})}=
\Vert f\Vert_{L^2(B,\mu_{B ,h/2})}= \Vert \widehat f \Vert_{L^2(E,\mu_{E,h/2})}.
 \label{(4.D.11)}
\end{equation}
By density, $T_h$ extends as a partial isometry from $L^2( B,\mu _{B,h/2})$
in $L^2( B^2,\mu _{B^2 ,h})$, denoted by $T_h$ as well. One sees that, for
$f\in {\cal D}_E$,
 $T_h f$ is anti-holomorphic if one identifies $(x,\xi)$
 and $x+ i \xi$. This application  $T_h$  is not surjective. One can find in
 Driver Hall \cite{D-H} and Hall \cite{HA} a study of its image.

\bigskip

We shall need norms for the elements of  ${\cal D}$. For all $E$ in
 ${\cal F} (B')$, for all $f$ in ${\cal D}_E$, written as
 $f = \widehat  f \circ P_E$, with $\widehat f $ in ${\cal S}_E$, for 
every integer  $m\geq 0$ and every positive $h$, one sets:
\begin{equation}
I_{E,m,h}(f) =(2\pi h)^{-{\rm dim}(E)} \int_{E^2}|(\widehat T_h \widehat f) (X)|
 (1+ |X|)^m e^{-\frac{1}{4h} |X|^2} d\lambda (X). \label{(4.D.12)} 
\end{equation}
This integral is convergent, for, according to (\ref{(4.D.9)}):
%---
$$ I_{E ,m,  h} (f) =(2\pi h)^{-{\rm dim}(E)}  \int _{E^2}  |< \gamma _{E,h/2}
\widehat f ,\Psi_{X,h}^{E, Leb} > _{Leb}|  (1+ |X|)^m d\lambda (X) .$$
%--
Since the function $\gamma _{E,h/2} \widehat f $ is in the Schwartz space
${\cal S}(E)$, it is well known that this integral is convergent.

\bigskip

For all $\widehat f$ in  ${\cal S}_E$ one can write, in the sense of integrals
of functions valued in $L^2(B,\mu _{B,h/2})$,
\begin{equation} 
\widehat  f = (2 \pi h) ^{-{\rm dim}(E)} \int _{E^2} e^{-\frac{1}{4h} |X|^2}
 (\widehat T_h \widehat f) (X)\Psi _{X,h} d\lambda (X) .\label{(4.D.13)} 
\end{equation}
The equality (\ref{(4.D.13)}) follows from (\ref{(4.D.6)}), (\ref{(4.D.9)})
and from the classical property, which makes sense since 
$\gamma _{E,h/2} \widehat f$ is in the Schwartz space ${\cal S}(E)$:
%-----
$$ \gamma_{E,h/2}\widehat f = (2\pi h)^{-{\rm dim}(E)} \int _{E^2}
< \gamma_{E,h/2} \widehat f,{\Psi}^{E, Leb}_{X,h}>_{Leb} {\Psi}^{E, Leb}_{X,h}
d\lambda (X) .$$
%------

\bigskip

We will need the following function, for all $X = (x,\xi)$, $Y= (y, \eta)$
and $Z = (z, \zeta)$ in $H^2$:
\begin{equation}
K_h^{AW} ( X, Y, Z)= e^{\frac{1}{2h} ( (x+i \xi)\cdot      (z-i\zeta) +
   (y-i\eta) \cdot  (z+i \zeta)  )}  . \label{(4.D.14)}   
\end{equation}

\begin{prop}\label{4.6}
Let $E$ be in ${\cal F}(B')$. Suppose that $E$ decomposes as $E=E_1\oplus E_2$,
 where $E_1$ and $E_2$ in ${\cal F}(B')$ are mutually orthogonal. One considers
 two functions   $\widehat f$ and $\widehat g$ in ${\cal S}_E$. One identifies 
$\widehat T_h \widehat f$ with a function $(\widehat T_h\widehat f)(X_1,X_2)$
on $(E_1)^2\times (E_2)^2$ and proceeds similarly for $\widehat T_h \widehat g$.
Then, for all  $Z_1$ in  $E_1$ and all $X_2$ and $Y_2$ in $E_2$:
\begin{equation}
\label{(4.D.15)}  
\int_{(E_1)^4} K_h^{AW} (X_1,Y_1,Z_1) (\widehat T_h \widehat f) (X_1,X_2)
\overline {(\widehat T_h \widehat g) (Y_1,Y_2)} d\mu _{(E_1)^4, h} (X_1, Y_1) =% \\
  (\widehat T_h\widehat f) (Z_1,X_2)
 \overline {(\widehat T_h \widehat g) (Z_1,Y_2)} 
\end{equation}
\end{prop}

\bigskip

{\it  Proof.} For all $X_2=(x_2,\xi_2)$ in $(E_2)^2$ and all $t$ in $E_1$ set:
%---
$$ \varphi _{X_2} (t) =  e^{-\frac{1}{4h}  (x_2-i\xi_2)^2  } \
\int _{E_2} \widehat f(t, y) e^{ \frac{1}{h} y \cdot (x_2 - i \xi_2) } d\mu_{E_2,h/2} (y).
 $$
%----
One can apply (\ref{(4.D.8)}), replacing  $\varphi $ with $\varphi _{X_2} $
and  $E$ with $E_1$. One gets:
$$ 
(\widehat T_h\widehat f) (Z_1,X_2) = \int _{(E_1)^2} e^{\frac{1}{2h} 
 (x_1+i \xi_1)\cdot      (z_1-i\zeta_1)} (\widehat T_h\widehat f) (X_1,X_2)
d\mu _{(E_1)^2, h} (X_1).
 $$
By a similar treatment  of $g$, one deduces (\ref{(4.D.15)}) 

%%%%%%%%%%%%%%%%%%%%%%%%%%%%%%%%%%%%%%%%%%%%%%%%%%%%%%%%%%%%%%%%%
%%%%%%%%%%%%%%%%%%%%%%%%%%%%%%%%%%%%%%%%%%%%%%%%%%%%%%%%%%%%%%%%%
\subsection{Wigner Gaussian function }\label{4.E}

The function $\widehat H_h^{Gauss} (\varphi ,\psi )$ has been defined in
 (\ref{(1.9)}) for functions $\varphi $ and $ \psi $ in ${\cal S}_E$, 
($E\in {\cal F} (H)$). We well see below that this definition  extends to
 $\varphi $ and $\psi $ in $L^2(E,\mu_{E,h/2})$. This is the case for the 
function $\Psi_{X,h}$ ($X = (x ,\xi)$)  defined in  (\ref{(4.D.3)}), which
is in  $L^2(E,\mu _{E,h/2})$ if $E$  contains $x$ and $\xi$. We shall 
compute $\widehat H_h^{Gauss}(\varphi,\psi)$ (for $\varphi$ and $\psi$ in
 ${\cal S}_E$)  using the Segal-Bargmann transforms of $\varphi$ and $\psi$.
\begin{prop}\label{4.7}{\ } 
\begin{enumerate}
\item
 For every subspace $E$  in ${\cal F} (H)$, for all $\varphi $ and $\psi $ in
 $L^2(E,\mu_{E,h/2})$, the equality (\ref{(1.9)}) defines a  function 
 $\widehat H^{Gauss}_h (\varphi ,\psi ) $ continuous  on $E^2 $ and satisfying,
 for all   $Z$ in $E^2$:
\begin{equation}
\left| \widehat H_h^{Gauss}(\varphi,\psi ) (Z)\right|\leq
e^{\frac{1}{h}|Z|^2}|| \varphi ||_{L^2(E,\mu_{E,h/2})}  \ ||\psi ||_{L^2(E,\mu_{E,h/2})} .
\label{(4.E.1)} 
\end{equation}
\item
 For all $X$, $Y$ and $Z$ in $H$ one can write:
\begin{equation}
\widehat H^{Gauss}_h (\Psi_{X,h},  \Psi_{Y,h} ) (Z)=
e^{-\frac{1}{4h} (|X|^2+ |Y|^2)} K_h^{Weyl} (X,Y, Z)\label{(4.E.2)}
\end{equation}
where one sets, for all $X=(x,\xi)$, $Y=(y,\eta)$ and $Z=(z,\zeta)$ in $H^2$:
\begin{equation}
 K_h^{Weyl} (X,Y, Z)= e^{\frac{1}{h} ( (x+i \xi)\cdot  (z-i\zeta) +(y -i \eta)
 \cdot (z+i\zeta) - \frac{1}{2} (x+i \xi) \cdot  (y -i \eta) )}\label{(4.E.3)}  
\end{equation}
\item
 Let $\varphi $ and $\psi $ be in  ${\cal S}_E$ ($E\in {\cal F}(H)$). Then the
 function $\widehat H_h^{Gauss}(\varphi,\psi)$ satisfies, for all $Z$ in $E^2$:
\begin{equation}
\widehat H^{Gauss}_h (\varphi , \psi ) (Z) =  \int _{E^4}
 K_h ^{Weyl} (X,Y, Z) (\widehat T_h  \varphi  ) (X) \overline {(
\widehat T_h \psi ) (Y) }  d\mu_{E^4, h} (X, Y)\label{(4.E.4)} 
\end{equation}
\end{enumerate}
\end{prop}

{\it  Proof. Point 1)}  We begin by proving (\ref{(4.E.1)}) for  $\varphi $
 and $\psi $ in  ${\cal S}_E$. For all $Z = (z,\zeta)$ in  $E^2$, by 
Cauchy Schwarz,
$$
\begin{array}{llll}
\displaystyle  \left| \widehat H_h^{Gauss}(\varphi,\psi) (z,\zeta)\right|
 & \displaystyle \leq  e^{\frac{1}{h}|\zeta|^2}
\int_E |\varphi (z+t)e^{-z\cdot t/h} | |\psi (z-t)
e^{ z\cdot t/h} | \  d\mu_{E,h/2}(t)\\
 & \displaystyle \leq e^{\frac{1}{h}|\zeta|^2}
 \left(\int_E |\varphi(z+t)|^2e^{-2 z\cdot t/h} \ d\mu_{E,h/2}(t)
\int_E |\psi (z+t)|^2  e^{-2 z\cdot t /h}  \  d\mu_{E,h/2}(t)\right)^{1/2}.
\end{array}
$$
Now,
$$ \int_E | \varphi(z+t)|^2 e^{-2 z\cdot t/h}  \  d\mu_{E,h/2}(t)=
 e^{|z|^2/h} \int_E | \varphi(t) |^2  \  d\mu_{E,h/2}(t).
$$
\smallskip

{\it Point 2) }  For all $p$ and $q$ in the complexified  $H_{\bf C}$ of $H$,
we remark that the Wigner Gauss transform of
%---
$$ f(t) = e^{p\cdot t } \hskip 2cm g(t) = e^{q\cdot t} \hskip 2cm t\in H$$
%----
satisfies, for all $Z = (z,\zeta)$ in $H^2$
\begin{equation}
\widehat H_h^{Gauss}(f,g) (z,\zeta) = e^{ {\frac{h}{4} } (p- \overline q)^2
+ z\cdot (p + \overline q)- i\zeta \cdot (p -\overline q)}. \label{(4.E.5)} 
\end{equation}
 To establish point 2), setting $X = (x,\xi)$,
 $Y = (y,\eta)$, it suffices to apply (\ref{(4.E.5)}) to
 $p = (x + i \xi)/h $, $q = (y + i \eta)/h$.

\smallskip

{\it Point 3)} Remark that, according to (\ref{(4.E.1)}), for all $Z$ in $E^2$,
 the application
$(\varphi,\psi)\rightarrow \widehat H^{Gauss}_h (\varphi,\psi) (Z)$ is a
continuous bilinear form on $L^2(E,\mu _{E, h/2}) \times L^2(E, \mu _{E,h/2}) $.
Moreover, for all $\varphi $ in ${\cal S}_E$, the application associating
 $ e^{-\frac{1}{4h} |X|^2} (\widehat T_h \varphi ) (X)\widetilde \Psi _{X,h}$
with  $X\in E^2$ is integrable, in the Bochner sense, with values in 
$L^2(E, \mu _{B,h/2})$, and its integral is given by (\ref{(4.D.13)}). The
same holds for $\psi $. It yields, as in Corollary 2 in  Section V.5 
in Yosida \cite{Y}:
%-------
$$\widehat  H^{Gauss}_h (\varphi, \psi) (Z) = (2\pi h)^{-2{\rm dim}(E)} \int _{E^4}
 e^{-\frac{1}{4h} (|X|^2+|Y|^2) }\widehat H^{Gauss}_h (\Psi_{X,h} ,\Psi_{Y,h} )(Z)
(\widehat T_h\varphi)(X) \overline {(\widehat T_h\psi)(Y)} d\lambda (X,Y). $$
%-----
The equality  (\ref{(4.E.4)}) follows from that, according to (\ref{(4.E.2)}).

\bigskip

The Wigner-Gauss transform  $H^{Gauss}_h (f,g)$ of two functions $f$ and $g$ 
in ${\cal D}_E$ ($E$ in ${\cal F}(B')$) is defined by (\ref{(1.11)}) if
 $f = \widehat f \circ P_E$  and $g = \widehat g\circ P_E$, with $\widehat f$ 
and $\widehat g$ in  ${\cal S}_E$. One has:
\begin{equation}
 H^{Gauss}_h (f,g) (Z) =  \int _{E^4}\widetilde K_h ^{Weyl} (X,Y, Z)
(\widehat T_h \widehat f)(X) \overline {(\widehat T_h \widehat g ) (Y) }
  d\mu_{E^4, h} (X, Y)  \hskip 2cm a.e. Z\in B^2 \label{(4.E.6)} 
\end{equation}
 with:
\begin{equation}
\widetilde K_h^{Weyl}(X,Y,Z)= e^{\frac{1}{h} (\ell_{x+i \xi}(z-i\zeta)
 +\ell_{y -i \eta} (z+i\zeta)  - \frac{1}{2} (x+i \xi) \cdot  (y -i \eta) )}
  \hskip 2cm (X,Y)\in H^2,\ a.e. Z\in B^2. \label{(4.E.7)} 
\end{equation}
Indeed, according to  (\ref{(4.E.3)}), (\ref{(4.E.7)}) and (\ref{(4.B.2)}),
 one has, for all $X$ and $Y$ in $E^2$, for almost every $Z$ in $B^2$
\begin{equation}
 K_h ^{Weyl} (X,Y, P_E(Z))=  \widetilde K_h^{Weyl} (X,Y, Z).\label{(4.E.8)}
\end{equation}
 Note that
%---
$$ H_h^{Gauss}  (f,1) (Z) = (T_hf) (2Z). $$
%------
Let us now investigate the bounds on the norm and the density extensions
of the Wigner-Gauss transformation.

\begin{prop}\label{4.8}
 For all  $f$ and  $g$ in ${\cal D}$, the Gaussian Wigner function
$H_h^{Gauss} (f, g)$ is in $L^1(B^2, \mu_{B^2, h/2})$.
The operator  which associates, with all functions $f$ and $g$ in ${\cal D}$,
 $H_h^{Gauss} (f,g)$, extends uniquely as a continuous bilinear application
with norm  $\leq 1$, from $L^2(B,\mu_{B, h/2})\times L^2(B ,\mu_{B, h/2})$ to
 $L^2(B^2,\mu_{B^2, h/4})$.
The operator  which associates, with all functions $f$ and $g$ in ${\cal D}$,
 $H_h^{Gauss} (f,g)$, extends uniquely as a continuous bilinear application
 from $L^2(B,\mu_{B, h/2})\times L^2(B ,\mu_{B, h/2})$ to the space
of continuous functions on  $H^2$ and for all $Z$ in $H^2$:
\begin{equation}
\left| H_h^{Gauss}(f,g) (Z)\right|  \leq
e^{\frac{1}{h}|Z|^2}||f||_{L^2(B,\mu_{B,h/2})} \ ||g||_{L^2(B,\mu_{B,h/2})} .\label{(4.E.9)} 
\end{equation}
\end{prop}

{\it  Proof.} Let  $f$ and $g$ be in  ${\cal D}$, let $E$ be such that $f$ and
 $g$ are in ${\cal D}_E$ and let $\widehat f$ and $\widehat g$ be functions  
defined on $E$ such that $f=\widehat f\circ P_E$ and $g =\widehat g\circ P_E$.
One knows that 
$H_h^{Gauss}(f,g)= \widehat H_h^{Gauss}(\widehat{f},\widehat{g})\circ P_E$, where 
$\widehat H_h^{Gauss} (\widehat{f},\widehat{g}) $ satisfies (\ref{(1.10)}). If
the functions $f$ and $g$ are in ${\cal D}_E$, the functions
 $\gamma_{E,h/2} \widehat f $ and $\gamma_{E,h/2} \widehat g$ are in
${\cal S}(E)$. According to Folland \cite{F} (Proposition 1.92), the function
$ H_h^{Leb} ( \gamma_{E,h/2} \widehat f,\gamma_{E,h/2} \widehat g) $ is in 
 ${\cal S} (E^2)$, hence in $L^1(E^2, \lambda)$. We deduce from that, according
 to (\ref{(1.10)}), that $\widehat  H_h^{Gauss} (\widehat{f},\widehat{g})$ 
is in $L^1(E^2, \mu _{E^2, h/2})$. The first point of the proposition
follows as a consequence, according to (\ref{(1.7)}).
The functions $\gamma_{E,h/2} \widehat f $ and $\gamma_{E,h/2} \widehat g $
  are in  $L^2(E,\lambda )$, where $\lambda $ is the Lebesgue measure on 
 $E$. According to  Folland \cite{F} (Proposition 1.92), one has
$$
\begin{array}{lll}
\displaystyle
\Vert H_h^{Leb} ( \gamma_{E,h/2} \widehat f ,\gamma_{E,h/2} \widehat g)
\Vert_{L^2 (E^2,\lambda)}
& \displaystyle \leq (2\pi h)^{ {\rm dim}(E)/2 }
\Vert \gamma_{E,h/2} \widehat f \Vert_{L^2 (E ,\lambda)}
\ \Vert \gamma_{E,h/2} \widehat f \Vert_{L^2 (E,\lambda)}\\
& \displaystyle=  (2\pi h)^{ {\rm dim}(E)/2 }
\Vert \widehat f\Vert_{L^2(E,\mu_{E, h/2})}\ \Vert \widehat g\Vert_{L^2(E,\mu_{E, h/2})}.
\end{array}
$$
It follows that, according to (\ref{(1.10)}):
%-----
$$\Vert \widehat  H_h^{Gauss}(\widehat f ,\widehat g) \Vert _{L^2(E^2, \mu _{E^2, h/4})}
\leq \Vert \widehat f \Vert _{L^2(E, \mu _{E, h/2})} \
 \Vert \widehat g \Vert _{L^2(E, \mu _{E, h/2})}. $$
%----
Hence, according to (\ref{(1.7)}):
%----
$$\Vert    H_h^{Gauss} (f,g) \Vert_{L^2(B^2, \mu _{B^2, h/4})}
\leq \Vert  f \Vert_{L^2(B, \mu _{B, h/2})} \ \Vert g \Vert_{L^2(B, \mu _{B, h/2})}. 
$$
%-----
The extension from  $ L^2(B,\mu_{B, h/2})\times L^2(B,\mu_{B, h/2}) $ in the
space of continuous functions on  $H^2$ is a consequence of (\ref{(4.E.1)}).

The extensions of the Wigner-Gauss transformation will also be denoted by
$H_h^{Gauss}$. For all $X$ and $Y$ in $H^2$, for almost
 $\mu_{B^2, h/4}$  all $Z$ in $B^2$ one has:
\begin{equation} 
H^{Gauss}_h ( \widetilde \Psi_{X,h}, \widetilde  \Psi_{Y,h} ) (Z) =
e^{-\frac{1}{4h} (|X|^2+ |Y|^2)} \widetilde K_h^{Weyl} (X,Y, Z)\label{(4.E.10)} 
\end{equation}
One verifies that this function (of $Z$) has a $L^2 (B^2 ,\mu _{B^2,h/4})$
norm equal to $1$.

 The following proposition  will be useful in Section  \ref{7}.

\begin{prop}\label{4.9}
If $(E_n)$ is an increasing sequence of subspaces in  ${\cal F}(B')$, whose
union is dense in $H$, then for all  $X$ and $Y$ in $H^2$:
\begin{equation}
\lim _{n\rightarrow +\infty}
 H^{Gauss}_h ( \Psi_{X,h}\circ P_{E_n},\Psi_{Y,h}\circ P_{E_n} )(Z) =
e^{-\frac{1}{4h} (|X|^2+ |Y|^2)} \widetilde K_h^{Weyl} (X,Y, Z). \label{(4.E.11)} 
\end{equation}
The convergence (as a function of $Z$) is in the sense of
 $L^2(B^2,\mu_{B^2,h/2})$.
\end{prop}

{\it  Proof.}
 According to (\ref{(1.11)}) and (\ref{(4.E.2)}),
$$
 H^{Gauss}_h ( \Psi_{X,h}\circ  P_{E_n},  \Psi_{Y,h}\circ  P_{E_n}  )(Z)
= e^{-\frac{1}{4h}(|X|^2+ |Y|^2)}K_h^{Weyl}(X,Y, P_{E_n}(Z)).
$$
Thanks to the definition (\ref{(4.E.3)}) and to (\ref{(4.7)}),
one proves that 
$$
\lim_{n\rightarrow \infty} K_h^{Weyl}(X,Y, P_{E_n}(Z))
= \widetilde K_h^{Weyl}(X,Y,Z),
$$
the limit being taken in $L^2(B^2, \mu_{B^2,h/2})$.
This proves the proposition.

%%%%%%%%%%%%%%%%%%%%%%%%%%%%%%%%%%%%%%%%%%%%%%%%%%%%%%%%%%%%%%%%%
%%%%%%%%%%%%%%%%%%%%%%%%%%%%%%%%%%%%%%%%%%%%%%%%%%%%%%%%%%%%%%%%%
\subsection{Convergence  in Definition \ref{1.2}.}\label{4.F}

If the function $\widetilde F$ is bounded and if $f$ and $g$  are in
 ${\cal D}$, the expression $Q_h^{Weyl}(\widetilde F)(f,g) $ from Definition
\ref{1.2} is well defined, according to Proposition \ref{4.8}. We show
that the same property is true under the weaker condition that, for a
given $m\geq 0$, the norm $ N_m (\widetilde F)$ in (\ref{(1.12)})  is finite.
This condition is satisfied, for example,  for functions of the form
$\widetilde F(x,\xi) =\prod \ell_{a_j}(x)^{\alpha_j}\prod \ell_{b_k}(\xi)^{\beta_k} $,
with $a_j$ and $b_j$ in $H$, 
using (\ref{(4.3bis)}). One can  show that $ N_m(\widetilde F)<\infty$,
with $m =\sum \alpha_j +\sum \beta_k$. 

\begin{prop}\label{4.10}
Let $\widetilde F$ be a Borel function on $B^2$, such that the norm 
$ N_m (\widetilde F)$ defined in (\ref{(1.12)}) is finite for a nonnegative 
integer $m$. For each subspace $E\in {\cal F}(B')$, for all $f$ and $g$ 
in ${\cal D}_E$, we have:
\begin{equation}
| Q_h^{Weyl}(\widetilde F) (f,g) | \leq  \ I_{E ,m,  h}(f) \ I_{E ,m,  h} (g)
\  N_m (\widetilde F),  \label{(4.F.1)} 
\end{equation}
where  $I_{E ,m,  h}(f) $ is defined in (\ref{(4.D.12)}).
\end{prop}

{\it Proof.} Let $f$ and $g$  be in  ${\cal D}_E$. According to Definition 
\ref{1.2}, to (\ref{(1.11)})  and to  Proposition \ref{4.7}, we only have
 to check that the function:
\begin{equation}
(X,Y, Z) \rightarrow
|\widetilde K_h^{Weyl }(X,Y, Z) |\ | (\widehat {T_h}\widehat{f}) (X)| \  |(
\widehat {T_h}\widehat{g} ) (Y) |\ |\widetilde F (Z)|  \label{(4.F.2)}
\end{equation}
is in $L^1(E^4\times B^2)$ for the measure  $d\mu_{E^4,h}(X,Y)d\mu_{B^2, h/2}(Z)$.
To this aim, one applies the change of variables  (\ref{(4.4)})
 with $B$ replaced with $B^2$, $h$ with
 $h/2$, $g(Z) =|\widetilde F(Z)
\widetilde K_h^{Weyl} (X,Y, Z)|$ and $a = (X+Y)/2$.
Note that 
\begin{equation}
 \widetilde K_h^{Weyl}\left(X,Y,Z+\frac{X+Y}{2} \right) 
e^{-\frac{1}{h} \ell _{X+Y} (Z)-\frac{1}{4h}  |X+Y|^2 } =  e^{\frac{1}{4h}  (|X|^2 +|Y|^2)} 
e^{\frac{i}{h} \varphi (X,Y, Z)}\label{(4.F.3)}  
\end{equation}

with
\begin{equation} \varphi (X,Y,Z) =\frac{1}{2} \sigma (X,Y)
+ \ell_{\xi-\eta}(z)-\ell_{x-y}(\zeta),
\label{(4.F.4)}
\end{equation}
where  $\sigma$  is the symplectic form. Consequently, for all  $X$ and $Y$
 in $H$:
%-----
$$ 
\int_{B^2} |\widetilde K_h^{Weyl }(X ,Y,Z)|\ |\widetilde F (Z)| d\mu _{B^2,h/2}(Z)
\leq e^{\frac{1}{4h}  (|X|^2 +|Y|^2)}
 \int_{B^2}  |\widetilde F\left( Z +\frac{X+Y}{2} \right ) |d\mu_{B^2,h/2}(Z).$$
%----
Hence, for every integer $m$ such that the norm  $N_m (\widetilde F)$ 
is finite, 
%-----
$$ \int_{B^2} |\widetilde K_h^{Weyl }(X,Y,Z)|\  |\widetilde F(Z)| d\mu_{B^2,h/2}(Z)
\leq e^{\frac{1}{4h}  (|X|^2 +|Y|^2)}  \left ( 1 + \frac{|X+Y|}{2} \right )^m
N_m (\widetilde F).$$
%------
The function (\ref{(4.F.2)}) is indeed in $L^ 1(E^4 \times B^2)$ for the 
announced measure.  Therefore
%-----
$$|Q_h ^{Weyl }( \widetilde F) (f,g)|\leq N_m (\widetilde F)  
\int _{E^4} e^{\frac{1}{4h}  (|X|^2 +|Y|^2)}  \left ( 1 + \frac{|X+Y|}{2} \right )^m |
 (\widehat {T_hf}) (X)| \  |(\widehat {T_hg} ) (Y) | d\mu_{E^4, h} (X, Y).$$
%-------
This proves the claim.

\bigskip

This implies that the quadratic form in (\ref{(1.13)}) is well defined 
if  $\widetilde F$ is bounded or if there exists an integer $m\geq 0$ such that
 $N_m (\widetilde F)$ is finite.

\bigskip

One can write, for all $f$ and $g$ in ${\cal D}_E$:
\begin{equation}
 Q_h ^{Weyl }( \widetilde F) (f,g) = \int _{E^4 \times B^2}
\widetilde K_h^{Weyl }(X,Y, Z)\ (\widehat {T_h}\widehat{f}) (X) \ \overline {(
\widehat {T_h}\widehat{g}) (Y)} \widetilde F (Z) d\mu_{E^4, h} (X, Y)
d\mu _{B^2,h/2}(Z). \label{(4.F.5)}
\end{equation}

%%%%%%%%%%%%%%%%%%%%%%%%%%%%%%%%%%%%%%%%%%%%%%%%%%%%%%%%%%%%%%%%%
%%%%%%%%%%%%%%%%%%%%%%%%%%%%%%%%%%%%%%%%%%%%%%%%%%%%%%%%%%%%%%%%%
\subsection{Anti-Wick and hybrid operators}\label{4.H}

We now rewrite the hybrid operator of Definition \ref{2.2} in a form more
suitable for norm estimates. At the same time we shall rewrite the anti-Wick 
operator defined by (\ref{(2.4)}) and give an expression nearer to
the usual definition. One will use the measure $\nu_{B^2,E^2,h_1, h_2}$ 
 defined in (\ref{(4.C.2)}).

\begin{prop}\label{4.11}
Let $f$ and $g$ be in  ${\cal D}$. Let $D$ be in ${\cal F} (B')$, $\widehat f$ 
and  $\widehat g$ in  ${\cal S}_D$ be such that $f = \widehat f \circ P_D$ and
 $g = \widehat g \circ P_D$. For every Borel function $\widetilde F$ bounded on 
 $B^2$, for every subspace $E$ in  ${\cal F}(B')$:
\begin{equation}
 Q_h^{hyb,{E}}(\widetilde F)(f,g)
 = \int_{ D^2\times D^2\times B^2} K_h^{hyb, E} (X,Y,Z)
(\widehat T_h \widehat f) (X) \overline {(\widehat T_h \widehat g ) (Y) }
 \widetilde F(Z)  d\mu _{D^4, h} (X,Y) d\nu _{B^2, E^2, h/2, h} (Z)
\label{(4.G.1)}  
\end{equation}
where  $ K_h^{hyb, E}$ is the function defined for all $(X,Y)$ in $H^4$ and
almost every $Z$ in $B^2$, by:
\begin{equation}
 K_h^{hyb, E} (X, Y, Z)  = \widetilde K_h^{Weyl} (X_E, Y_E, Z_E)
\widetilde K_h^{AW} (X_{E^{\perp}},Y_{E^{\perp}}, Z_{E^{\perp}}). \label{(4.G.2)} 
\end{equation}
Here  we set $X = (X_E,X_{E^{\perp}})$, ... according to Proposition \ref{2.1},
$\widetilde  K_h^{Weyl} $  has been defined by  (\ref{(4.E.7)}) and
\begin{equation}
\widetilde K_h^{AW} ( X, Y, Z)= e^{\frac{1}{2h} (\ell _{x+i \xi} ( (z-i\zeta) +
   \ell _{y-i\eta} (z+i \zeta)  )} . \label{(4.G.3)} 
\end{equation}
The above integral does not depend on the subspace $D$ such that $f$ and $g$
are in ${\cal D}_D$ (for $f$ and $g$ in ${\cal D}$.)
\end{prop}

{\it  Proof.} According to Definitions  \ref{2.2} and \ref{1.2}, one has
%---
$$ Q_h^{hyb,{E}}(\widetilde F)(f,g) = 
\int_{B^2}( \widetilde H_{E^{\perp}, h/2} \widetilde F)(Z) G(Z)
 d \mu _{B^2,h/2}(Z)
$$
%--
with
%---
$$ G(Z) = H_h^{Gauss} (f,g) (Z). $$
%---
This function  $G$ belongs to $L^{1} (B^2, \mu _{B^2,h/2})$ thanks to
Proposition \ref{4.8}. According to Proposition \ref{4.5},  applied with
 $p= +\infty$, $q=1$, $h_1 =h_2 =h/2$ and $t= h/2$,  we obtain:
%------
$$
 Q_h^{hyb, {E}}(\widetilde F) (f,g) =
\int_{B^2}\widetilde F(Z)
 ( \widetilde M_{E^{\perp} ,h/2, h/2, h/2}G)(Z) d\nu_{B^2, E^2, h/2, h}(Z).
$$
%-----
According to (\ref{(4.E.6)}), the form of $f$ and $g$ imply that the function
 $G= H_h^{Gauss} (f,g)$ has the following expression:
%----
$$ G(Z) = \int _{ D^4 } \widetilde K_h^{Weyl} (X,Y, Z)(\widehat T_h \widehat f ) (X) \overline {(
\widehat T_h \widehat g ) (Y) }
 d\mu _{D^4, h} (X,Y). $$
%----
Using the expression (\ref{(4.C.3)}) of $\widetilde M_{E^{\perp},h/2,h/2,h/2}$
 and Fubini's theorem, we obtain:
%----
$$ (\widetilde M_{E^{\perp} ,h/2, h/2, h/2}G) (Z)= \int _{ D^4 } \Phi (X, Y, Z)
(\widehat {T_h}\widehat{f}) (X) \overline {(
\widehat {T_h}\widehat{g}) (Y) }  d\mu _{D^4, h} (X,Y) $$
where
$$
\begin{array}{lll}
\displaystyle  \Phi(X,Y,Z) &  \displaystyle =\int_{(E^{\perp})^2} \widetilde
  K_h^{Weyl}\left(X,Y,\left (Z_E, V+\frac{1}{2}Z_{E^{\perp}}  \right ) \right )
d\mu_{(E^{\perp})^2, {\frac{h}{4}}}(V) \\ & \displaystyle
=\widetilde  K_h^{Weyl} (X_E,Y_E, Z_E)  \int _{(E^{\perp})^2}
\widetilde  K_h^{Weyl} \left (X_{E^{\perp}},Y_{E^{\perp}},V+\frac{1}{2}Z_{E^{\perp}}
 \right  ) d\mu_{(E^{\perp})^2, {\frac{h}{4}}}(V) .
\end{array}
$$
Here the functions $\ell_x,\dots,\ell_{\eta}$ appearing in the 
kernels are linear, since $x,\dots, \eta$
belong to $D\subset B'$.
One checks, using (\ref{(4.2)}), (with $B$ replaced with
$(E^{\perp})^2$ and $h$ with $h/4$), that:
%----
$$  \int _{(E^{\perp})^2}
\widetilde K_h^{Weyl} (X_{E^{\perp}},Y_{E^{\perp}},V+\frac{1}{2}Z_{E^{\perp}}  ) d\mu_{(E^{\perp})^2, {\frac{h}{4}}}(V)
=  K_h^{AW} (X_{E^{\perp}},Y_{E^{\perp}}, Z_{E^{\perp}}). $$
%---
The proposition follows from that.

\begin{coro}\label{4.12}
Let $f$ and $g$ be in  ${\cal D}$. Let $D$ be in ${\cal F} (B')$, $\widehat f$
 and $\widehat g$ in  ${\cal S}_D$ such that $f = \widehat f \circ P_D$ and
 $g = \widehat g \circ P_D$.  For every Borel function $\widetilde F$, bounded
 on  $B^2$, one has:
%---
$$ Q_h^{AW}(\widetilde F) (f,g)  = \int _{  B^2 } \widetilde F(Z)
(\widehat T_h\widehat f) (P_D(Z ))
 \overline {(\widehat T_h \widehat g) (P_D(Z ))}  d\mu _{B^2,  h} (Z). $$
 %-----
 There exists an operator  $Op_h^{AW}( \widetilde F)$, bounded in
 $L^2 ( B,\mu _{B,h/2})$, such that
 %---
$$ Q_h^{AW}(\widetilde F) (f,g)  = < Op_h^{AW}( \widetilde F) f, g>.$$
%---
The norm of this operator is smaller than  $\Vert \widetilde F\Vert _{\infty}$.
\end{coro}

{\it  Proof.} One applies Proposition \ref{4.11} with $E =\{ 0\}$,$E^{\perp}= B$.
This yields
%---
$$ Q_h^{AW}(\widetilde F)(f,g) = \int_{D^4\times B^2 } \widetilde K_h^{AW}(X,Y,Z)
(\widehat T_h \widehat f) (X) \overline {(\widehat T_h \widehat g ) (Y) }
 \widetilde F(Z)  d\mu _{D^4, h} (X,Y) d\mu _{B^2,  h} (Z) .  $$
%----
According to Proposition \ref{4.6}, with $E_1$ replaced with $D$ and $E_2$
 with $\{ 0 \}$, one has:
%-----
$$ \int _{D^4} \widetilde K_h^{AW} ( X, Y, Z) (\widehat T_h \widehat f) (X)
\overline {(\widehat T_h \widehat g) (Y)}
d\mu _{D^4, h} (X, Y) = (\widehat T_h\widehat f) (P_D(Z ))
 \overline {(\widehat T_h \widehat g) (P_D(Z ))}, $$
%------
since $\widetilde K_h^{AW}(X,Y,Z)= K_h^{AW}(X,Y,P_D(Z))$.
The corollary then follows, according to (\ref{(4.D.11)}).

%%%%%%%%%%%%%%%%%%%%%%%%%%%%%%%%%%%%%%%%%%%%%%%%%%%%%%%%%%%%%%%%%
%%%%%%%%%%%%%%%%%%%%%%%%%%%%%%%%%%%%%%%%%%%%%%%%%%%%%%%%%%%%%%%%%
\subsection{Partial Heat semigroups and stochastic extensions. }\label{4.I}

 We shall first be concerned with establishing a property for the stochastic
 extension given by Definition \ref{1.2}. Namely, we shall prove the fact
 that it commutes with heat operators acting on subspaces of
 finite dimension, provided  that these operators act on
 bounded Lipschitz functions on $H$.

\bigskip

For all  $E$  in ${\cal F}(B')$, and for all $t>0$, we shall use the operators
 $ H_{E,t}$ and $ \widetilde H_{E,t}$ both defined by (\ref{(2.6)}), but
 the first one acting in spaces of functions on $H^2$, whereas the second
 one acts in spaces of functions on $B^2$.

\begin{prop}\label{4.13}
 Let $ F$ be a bounded Lipschitz function  on $H^2$. Fix $h>0$ and $t>0$. 
Assume that $F$ has a stochastic extension $\widetilde F$ for the measure
$\mu_{B,h}$. Then, for each  $E$ in ${\cal F}(B')$,  the function $H_{E,t}F$
 has a stochastic extension for the measure $\mu_{B,h}$, and this stochastic
extension is equal to $\widetilde H_{E,t} \widetilde F$.
\end{prop}

{\it Proof.}  Let $(E_n)$ be an increasing sequence in ${\cal F}(H)$, whose
 union is dense in $H$. In order to obtain that $( H_{E,t}{F})\circ P_{E_n}$
converges to $\widetilde H_{E,t}\widetilde F$ in $\mu_{B,h}-$ probability, it is
necessary and sufficient to verify that, for all subsequences
 $(( H_{E,t}{F})\circ P_{E_\varphi(n)})_n$, one may  extract a further subsequence
 which converges $\mu_{B,h}$-almost everywhere to
 $\widetilde H_{E,t}\widetilde F$. Set $\psi : \N\rightarrow \N$ such that 
$({F}\circ P_{E_{\varphi(\psi(n))}})_n$ tends to $\widetilde F$ $\mu_{B,h}$-almost
 everywhere. Set $G_n=E_{\varphi(\psi(n))}$. For $X \in B$,we have
$$
\begin{array}{lll}
\displaystyle  (H_{E,t}{F})\circ P_{G_n}(X)- \widetilde H_{E,t}\widetilde F(X)
\\
\displaystyle  = (2 \pi t)^{-\dim(E)}\int_{E^2}  e^{-\frac{|Y|^2}{ 2t}} 
\left( {F}( P_{G_n}(X)+Y) - {F}( P_{G_n}(X+Y))\right) \ d\lambda(Y)\\
\displaystyle + (2 \pi t)^{-\dim(E)}\int_{E^2} e^{-\frac{|Y|^2}{ 2t}} 
\left({F}(P_{G_n}(X+Y))-\widetilde F(X+Y)\right) \ d\lambda(Y).
\end{array}
$$
%-----
By assumption, there exists $K>0$ such that for all $ (Z,V) \in H^2$:
%---
$$|{F}(Z+V)-{F}(Z)|\leq K |V|,
$$
%----
which shows, by the dominated convergence theorem, that the first
term converges to 0 for all $X$.

We shall prove that the second term converges to $0$ in $L^1(B,\mu_{B,h})$, 
which will give the convergence $\mu_{B,h}$-almost everywhere for
  a subsequence. Set
%---
$$
A_n= \int_{B^2} \left|\int_{E^2}  e^{-\frac{|Y|^2}{ 2t}} 
\left({F}(P_{G_n}(X+Y)) -\widetilde F(X+Y)\right)\ d\lambda(Y)\right|\
 d\mu_{B^2 ,h}(X).
$$
%----
Fubini's theorem and the decomposition (\ref{(2.1)}),
 $ \mu_{B,h}(X)=\mu_{E,h}(X_E)\otimes \mu_{(E^{\bot}),h}(X_{E^\bot})$, give
%---
$$
A_n\leq \int_{E^4\times (E^{\bot})^2 } e^{-\frac{|Y|^2}{ 2t}}
 \left|{F}( P_{G_n}(X_E+Y,X_{E^\bot})) -\widetilde F(X_E+Y,X_{E^\bot}) \right|
 \ d\lambda(Y) d\mu_{E^2,h }(X_E) d\mu_{(E^{\bot})^2,h}(X_{E^\bot}).
$$
%-----------------
Then the change of variables $X_E+Y \rightarrow Z$ and an explicit 
computation on gaussian functions yield:
%---
$$
A_n\leq C(h,t) \int_{E^2\times (E^{\bot})^2 }
 \left|{F}( P_{G_n}(Z,X_{E^\bot})) -\widetilde F(Z,X_{E^\bot})\right|
 e^{|Z|^2 (\frac{1}{ 2h} -\frac{1}{ 2(t+h)})}  
  d\mu_{E^2,h}(Z)  d\mu_{(E^{\bot})^2,h}(X_{E^\bot}).
$$
%---
We return to $B$ setting $W=(Z,X_{E^\bot})$. We have
%---
$$
A_n\leq C(h,t) \int_{B^2 } \left|{F}( P_{G_n}(W)) -\widetilde F(W) \right|
 e^{|P_{E}(W)|^2 (\frac{1}{ 2h} -\frac{1}{ 2(t+h)})}    d\mu_{B^2,h}(W).
$$
%---
The term in the above integral tends to $0$ $\mu_{B,h}$ almost
everywhere and is bounded by
%---
$$
\Psi (W)= 2 e^{|P_{E}(W)|^2 (\frac{1}{2h} -\frac{1}{ 2(t+h)})} . \sup_{X\in H^2}| F(X)|.$$
%---
Using (\ref{(1.7)}), an explicit computation shows that the bound $\Psi$ is
in $L^1(B,\mu_{B,h})$. The second term tends indeed to $0$ in
$L^1(B,\mu_{B,h})$, thus almost everywhere for $\mu_{B,h}$, when extracting again
 a subsequence, indexed by $\zeta$. The term 
$(H_{E,t}\widetilde{F})\circ P_{E_{\varphi(\psi(\zeta(n)))}}-H_{E,t}F$ is then a sum of
 two terms tending  $\mu_{B,h}$-almost everywhere to $0$.\hfill $\square$ 

\bigskip

Using  the Hilbertian basis $(e_j)_{(j\in \Gamma )}$ chosen for our construction,
 we defined in (\ref{(3.1)}), for each finite subset $I$ in $\Gamma$,
 operators denoted by $ \widetilde T_{I, h}$ and $ \widetilde S_{I, h}$ .
These operators act in the space of bounded Borel functions on $B^2$. We denote
 by $ T_{I, h}$ the operator defined similarly to (\ref{(3.1)}), but acting in
 the space of bounded continuous functions on $H^2$.

\begin{prop}\label{4.14}
 If $\widetilde F$ is a Borel function  bounded on $B ^2$ and 
 is the stochastic extension
 for the measure $\mu _{B^2, h}$ of a function $ F$ in $ S_1(M, \varepsilon)$,
 the family $(\varepsilon_j)_{(j\in \Gamma)}$ being square summable, then
$\widetilde T_{I,h}\widetilde F$ is the stochastic extension of
the function $ T_{I, h} F$ for the measure $\mu _{B^2,h}$.
\end{prop}

{\it Proof.} Since $T_{I,h}= \prod_{j\in I}(I-H_{e_j})$, products like
 $\prod_{j\in J}H_{e_j}= H_{E,h/2}$ appear,  with 
$E=E(J) ={\rm Vect}(e_j,j\in J), J\subset I$. We then apply Proposition
 \ref{4.13} for each of these terms, replacing
 $t$ with $h/2$ and $E$ with  $E(J)$, $J\subseteq I$. One may use this
proposition since every  function $ F$ in $ S_1(M, \varepsilon)$
verifies, for all $X$ and $V$ in $H^2$:
%---
$$ | F(X+ V) -  F(X)| \leq M |V|\sqrt {2}  \left [
\sum_{j\in \Gamma} \varepsilon_j^2 \right ]^{1/2}. $$
%---
Indeed, let $\eta>0$ and let $J$ be a finite subset of $\Gamma$ such
that $ | {F}(X+V)- {F}(X+ P_{E(J)} (V))|<\eta.$ We obtain, using the
 hypothesis of differentiability along the $u_j$ and $v_j$,
$$
|{F}(X+P_{E(J)} (V)) - {F}(X)|\leq \sum_{j\in J} M\varepsilon_j(
|V\cdot u_j| + |V\cdot v_j|).
$$
Then
$$
| {F}(X+V)- {F}(X)| \leq  \eta+ M
\left(2\sum_{j\in \Gamma}\varepsilon_j^2\right)^{1/2}|V|,
$$
 which gives the result.\hfill $\square$ 

%%%%%%%%%%%%%%%%%%%%%%%%%%%%%%%%%%%%%%%%%%%%%%%%%%%%%%%%%%%%%%%%%
%%%%%%%%%%%%%%%%%%%%%%%%%%%%%%%%%%%%%%%%%%%%%%%%%%%%%%%%%%%%%%%%%
\section{Proof of Proposition  \ref{3.1}.}\label{5}

Proposition  \ref{3.1}  is a direct consequence of Propositions \ref{5.2} and
 \ref{5.3} below.

\bigskip

For every integer $m$ and for any subset $I$ in $\Gamma$, ${\cal M}_m(I)$ 
designates the set of mappings $(\alpha,\beta)$ from $I$ into 
$\{ 0, 1,..., m \}\times \{ 0, 1,..., m \} $. When $(e_j)_{(j\in \Gamma)}$ 
is the Hilbertian basis of $H$ chosen for our construction, we set
 $u_j = (e_j,0)$ and $v_j = (0, v_j)$ for all $j\in \Gamma$. For each
 multi-index $(\alpha,\beta)$ in ${\cal M}_m(I)$, we set:
\begin{equation}
\partial_u^{\alpha} \partial_v^{\beta} = \prod_{j\in I} \partial_{u_j} ^{\alpha_j}
 \partial _{v_j} ^{\beta_j} .\label{(5.1)} 
\end{equation}
For any function $ F$ in $ S_m(M, \varepsilon)$, $(m \geq 0)$, let us define :
\begin{equation}
N_{I ,h}^{(m)}( F)  = \sum _{(\alpha,\beta)\in {\cal M}_m(I)}
h^{(|\alpha|+|\beta|)/2}  \left \Vert \partial _u^{\alpha} \partial
_v^{\beta}   F \right \Vert _{\infty } \label{(5.2)} 
\end{equation}
where $\Vert \cdot \Vert _{\infty}$ denotes the supremum  norm for
bounded functions on $H^2$.  For a function defined
on a subspace $E$ of $H$, let $N_{I ,h}^{(m)}( F)$
denote the same expression, the supremum being taken on $E$. 
 Adapting the proof of Unterberger
\cite{U-2} or \cite{A-J-N-1}, we obtain the following proposition.

\begin{prop}\label{5.1}
Let   $I$ be a finite subset of $\Gamma$ and $E = E(I)$ be the subspace 
spanned by the $e_j$ $(j\in I)$. Let $H$ be a continuous function on $E^2$, 
such that $N_{I ,h}^{(2)}( H)$ is well defined. We set, for all $U$ and $V$
 in $L^2(E^2, \mu_{E^2,h})$:
\begin{equation}
 Q(U, V) = \int_{E^6} K_h^{Weyl} (X,Y, Z) U(X) V(Y) H(Z) d\mu _{E^4, h} (X,Y)
d\mu _{E^2,h/2} (Z) \label{(5.3)}
\end{equation}
Then, we have:
\begin{equation}
|Q(U,V)|\leq \left( \frac{9 \pi}{2}\right)^{|I|}N_{I,h}^{(2)}(H) 
\Vert U \Vert_{ L^2(E^2,\mu_{E^2,h})} \Vert V \Vert_{ L^2(E^2,\mu_{E^2,h})} \label{(5.4)} 
\end{equation}
\end{prop}

\bigskip

{\it Proof.}  We shall study the following integral kernel:
%----
$$ \Phi (X,Y) = \int_{E^2} K_h^{Weyl} (X,Y, Z) H(Z) d\mu _{E^2,h/2} (Z).$$
%-----
 By the change of variables  $Z \rightarrow Z + (X+Y)/2$, using (\ref{(4.E.3)})
 one obtains:
%------
$$ \Phi (X,Y) = e^{\frac{1}{4h} (|X|^2 + |Y|^2)} \int _{E^2}
e^{\frac{i}{h}\psi (X,Y,Z) } H \left (Z +\frac{X+Y}{2}\right) d\mu _{E^2,h/2}(Z)$$
%----
where $\psi (X,Y, Z)$ is defined by
$$
\psi (X,Y, Z)=\frac{1}{2}\sigma(X,Y) -\sigma(Z,X-Y).
$$
 Integrating by parts,
we  find an upper bound for this integral.
 For every $X= (x,\xi)$ in $E^2$, set:
%---
$$ K_E(X) = \prod _{j\in I} \left ( 1 + \frac{x_j^2}{ h} \right )
 \left ( 1 + \frac{x_j^2}{ h} \right ) \hskip 2cm x_j =  e_j (x) = u_j(X)
 \hskip 2cm \xi _j = e_j(\xi) = u_j(X).$$
 %--
Integrations by parts then show that, for all $X$ and $Y$ in $E^2$:
 %---
 $$K_E(X - Y )  | \Phi ( X,Y )| \leq e^{\frac{1}{4h} (|X|^2 + |Y|^2)}
\int_{E^2}\left |(L H\left (Z+\frac{X+Y}{ 2}\right)\right |d\mu_{E^2 ,h/2} (Z) $$
%---
where $L$ is the differential operator defined by:
%-----
$$ L =\prod_{j\in I} L_{z_j} L_{\zeta_j} \qquad L_{z_j}=
 \sum_{k=0}^2 a_k \left (\frac{z_j}{\sqrt {h}}\right )h^{k/2}\partial _{u_j}^k $$
%---
where
$$  a_0 (z) = 3 - 4 z^2 \qquad a_1(z) = 4z \qquad a_2(z) = -1.
$$
%----
We then deduce that:
%----
$$K_E(X-Y)|\Phi (X,Y)|\leq e^{\frac{1}{4h}(|X|^2+|Y|^2)}
\sum_{(\alpha,\beta)\in {\cal M}_2(I)}h^{(|\alpha|+|\beta|)/2}
\left \Vert \partial_u^{\alpha} \partial_v^{\beta} H \right
\Vert_{\infty} \prod_{j\in I} C_{\alpha_j}C_{\beta_j} $$
%---
with:
%---
$$ C_k = \pi ^{-1/2} \int_{\R} |a_k(z)|e^{-z^2} dz \hskip 2cm k=0, 1,2.$$
%-----
Notice that $\max ( C_0,C_1,C_2) \leq 3$. Therefore:
%----
$$ | \Phi (  X,Y)| \leq e^{\frac{1}{4h} (|X|^2 + |Y|^2)}
 9^{|I|} N_{I ,h}^{(2)}( H) \ K_E(X - Y )  ^{-1} .$$
%------
We have:
\begin{equation}
 \int _{E^2} K_E(X)^{-1} d\lambda (X) = (h\pi^2)^{{\rm dim}(E)}. \label{(5.5)} 
\end{equation}
Hence,
%-----
$$ |Q(U, V) |\leq (2\pi h)^{-2{\rm dim}(E)} 9^{|I|} N_{I ,h}^{(2)}( H)
\int_{E^4}  K_E(X-Y )^{-1}e^{-\frac{1}{4h} (|X|^2 + |Y|^2)} U(X) V(Y) d\lambda (X,Y).
 $$
%-------
Using Schur's lemma and (\ref{(5.5)}), one obtains (\ref{(5.3)}) and 
(\ref{(5.4)}). The proof of the proposition is then completed.\hfill $\square$ 

\bigskip

The first step of the proof of Proposition \ref{3.1} is the following result:

\bigskip

\begin{prop}\label{5.2}
Let $G$ be a function  on $H^2$ and $I$ be a finite subset of $\Gamma$,
 such that $N_{I ,h}^{(2)}( G)$ is well defined. We suppose that $G$ has a  
stochastic extension $\widetilde G$ (for the measure $\mu _{B^2,h}$). Let
 $E = E(I)$ be the subspace spanned by the $e_j$ $(j\in I)$. Then, there 
exists a bounded
operator $Op_h ^{hyb, E(I )} (\widetilde G)$ in $L^2(B,\mu_{B,h/2}) $
such that, for all $f$ and $g$ in ${\cal D}$:
\begin{equation}
 Q_h^{hyb, E(I)}(\widetilde G)( f,g) = < Op_h ^{hyb, E(I )} (\widetilde G) f,g>.
\label{(5.6)}
\end{equation}
Moreover, we have:
\begin{equation} \Vert Op_h^{hyb,E(I)}(\widetilde G)\Vert_{{\cal L}(L^2(B,\mu_{B h/2}))}
 \leq \left ( \frac{9\pi}{2} \right)^{|I|}N_{I ,h}^{(2)}( G) \label{(5.7)} 
\end{equation}
\end{prop}

{\it Proof.} We use the expression for the quadratic form
given by Proposition \ref{4.11}. Let $D\in {\cal F}(B')$ be a subspace
such that  $f$ and $g$ are  in ${\cal D}_D$. Since the quadratic form defined
 in  (\ref{(4.G.1)}) does not depend on $D$, we may assume that $D$ contains
 $E =E(I)$. Let $S$ be the orthogonal of $E$ in $D$. The variable of $B^2$ may
 be written $Z = (Z_E,Z_S,Z_{D^{\perp}})$ with $Z_{E^{\perp}}=(Z_S,Z_{D^{\perp}})$. 
Let $\widehat f$ and $\widehat g$ be in ${\cal S}_D$ such that 
$f = \widehat f \circ P_D$ and $g = \widehat g \circ P_D$. The Segal Bargmann 
transforms $\widehat T_h \widehat f $ and $\widehat T_h \widehat g $  are
 functions of $Z_D = (Z_E,Z_S)$. With these notations, the equality
 (\ref{(4.G.1)}) may be written:
 %---
$$ 
\begin{array}{lll}
\displaystyle
Q_h^{hyb, {E(I)}}(\widetilde G) (f,g ) 
 = &\displaystyle  \int_{ E^4 \times S^4 \times B^2 }
K_h^{Weyl} (X_E, Y_E, Z_E)K_h^{AW} (X_S, Y_S, Z_S)\\
 &\displaystyle (\widehat T_h \widehat f)(X_E,X_S) 
\overline {(\widehat T_h\widehat g)(Y_E,Y_S)}
 \widetilde G(Z_E, Z_S,Z_{D^{\perp}} ) \\
& \displaystyle  d\mu _{E^4, h} (X_E,Y_E) d\mu _{S^4, h} (X_S,Y_S)
 d\mu_{E^2,h/2}(Z_E) d\mu_{S^2,h}(Z_S) d\mu_{(D^{\perp})^2,h} (Z_{D^{\perp}}).
\end{array}
$$
 %----
One applies Proposition \ref{4.6} with $E_1$ replaced with $S$ and $E_2$ with
 $E$. This allows to write:
%-----
$$
\begin{array}{lll}
\displaystyle Q_h^{hyb, {E(I)}}(\widetilde G) (f,g) &\displaystyle
 = \int _{ E^4 \times B^2 }
K_h^{Weyl} (X_E, Y_E, Z_E)(\widehat T_h \widehat f) (X_E, Z_S) \overline {(
\widehat T_h \widehat g ) (Y_E, Z_S) }  \widetilde G(Z_E, Z_S,Z_{D^{\perp}} ) 
\\ & \displaystyle d\mu _{E^4, h} (X_E,Y_E)
 d\mu _{E^2,  h/2} (Z_E) d\mu _{S^2,  h} (Z_S)d\mu _{(D^{\perp})^2,  h} (Z_{D^{\perp}} )
\end{array}
 $$
 %----
 Let $(\Lambda _n)$ be an increasing sequence of finite subsets of $\Gamma$,
containing $I$ and whose union is $\Gamma$. Since $\widetilde G$ is a
 stochastic extension of $G$  in the sense of Definition \ref{4.4},
 one remarks that
$$
\lim_{n\rightarrow \infty}||  G\circ P_{E(\Lambda_n)}-
 \widetilde G||_{L^2( E^2\times (E^{\bot})^2,\mu_{E^2, h/2}\otimes \mu_{(E^{\bot})^2, h})} =0.
$$
 Hence, for a subsequence, $ G\circ P_{E(\Lambda_{\varphi(n)})}$ converges to $
 \widetilde G$  $\mu_{E^2, h/2}\otimes \mu_{(E^{\bot})^2, h}$-almost everywhere.
 Since $X\mapsto < \gamma_{D,h/2}\widehat f, \Psi_{X,h}^{D,Leb}>_{Leb}  $
and $X\mapsto < \gamma_{D,h/2}\widehat g, \Psi_{X,h}^{D,Leb}>_{Leb}  $
are rapidly decreasing  one can write, using the dominated convergence
 Theorem, 
 $$
\begin{array}{lll}
\displaystyle Q_h^{hyb,{E(I)}}(\widetilde G)(f,g)
&\displaystyle =\lim_{n\rightarrow \infty } \int_{ E^4 \times B^2 }
K_h^{Weyl} (X_E, Y_E, Z_E)(\widehat T_h \widehat f) (X_E, Z_S) \overline {(
\widehat T_h \widehat g ) (Y_E, Z_S) }\\
 &\displaystyle (G\circ P_{E(\Lambda_{\varphi(n)})} (Z_E,Z_S ,Z_{D^{\perp}}) d\mu_{E^4,h}(X_E,Y_E)
 d\mu _{E^2,  h/2} (Z_E) d\mu _{S^2,  h} (Z_S)d\mu _{(D^{\perp})^2,  h} (Z_{D^{\perp}} ).
\end{array}
 $$
For all $n$, all  $Z_S$ in  $S^2$ and $Z_{D^{\perp}}$ in $(D^{\perp})^2$,
 one defines a continuous function on $E^2$ by
  %---
$$H_{n, Z_S ,Z_{D^{\perp}} } (Z_E) = (G\circ P_{E(\Lambda _n)} (Z_E, Z_S,Z_{D^{\perp}} ) .$$
   %----
The norm  $N_{I ,h}^{(2)}$ of this function is well defined by (\ref{(5.2)}),
  the supremum
being taken on $E$ and satisfies
\begin{equation}
 N_{I ,h}^{(2)}( H_{n, Z_S ,Z_{D^{\perp}} }) \leq  N_{I ,h}^{(2)}( G). \label{(5.XXX)}
\end{equation}
One can, then, apply Proposition \ref{5.1}, $H$ being replaced with
 $H_{n, Z_S,Z_{D^{\perp}} }$ and $U( X_E) = (\widehat T_h \widehat f) (X_E, Z_S) $,
 $V( X_E) = (\widehat T_h \widehat g) (Y_E, Z_S) $.
Using Proposition \ref{5.1}, Cauchy Schwarz's inequality and (\ref{(5.XXX)}),
 one obtains, taking the limit on $n$,
 %-------
 $$ Q_h^{hyb, {E(I)}}(\widetilde G) (f,g) \leq  \left ( \frac{9\pi}{2} \right
)^{|I|}N_{I ,h}^{(2)}( G)A(f) A(g)  $$
%------
 where:
 %---
 $$ A(f)^2 = \int_{E^2\times S^2}  |(\widehat T_h \widehat f) ( X_E, Z_S)|^2
 d\mu _{E^2, h} (X_E ) d\mu _{S^2,  h} (Z_S).$$
 %----
 According to (\ref{(4.D.11)}), one has $A(f) =\Vert f \Vert_{ L^2(B,\mu_{B,h/2})}$.
 Hence,
  %-------
 $$ Q_h^{hyb, {E(I)}}(\widetilde G) (f,g) \leq \left(\frac{9\pi}{2} \right)^{|I|}
N_{I ,h}^{(2)}( G)\Vert f \Vert _{ L^2(B,\mu_{B,h/2})}
\Vert g \Vert _{ L^2(B,\mu_{B,h/2})}. $$
%------
The proposition is thus easily deduced.\hfill $\square$ 

\bigskip

We next involve the operator $T_{I,h}$ defined in (\ref{(3.1)}). Proposition
 \ref{3.1} will then follow from Proposition \ref{5.2}, applied with
$ G = T_{I, h} F$ combined with Proposition \ref{5.3} below, applied to $ F$.
 Indeed, Proposition \ref{4.14} shows that the stochastic extension of $ G$ for
 the measure $\mu_{B^2,h}$ is $\widetilde G = \widetilde T_{I,h}\widetilde F$.

\begin{prop}\label{5.3}
For any  $F$ in $S_2(M,\varepsilon)$, for every finite subset $I$ of $\Gamma$
and for all $h>0$, the function $T_{I,h}F$ defined after (\ref{(3.1)})
 satisfies:
\begin{equation}
 N_{I ,h}^{(2)}(T_{I, h} F)\leq M (18 S_{\varepsilon} h)^{|I|}
\prod _{j\in I}\varepsilon_j^2  \label{(5.7N)}
\end{equation}
\end{prop}

Setting $u_j = (e_j,0)$ and $v_j=(0,e_j)$, we write
 $\Delta_j =\partial_{u_j}^2 + \partial_{v_j}^2 $,  and we consider the operator
 $H_{D_j, h/2}$ defined as in (\ref{(2.8)}), but acting on the functions on
 $H^2$. The proof of Proposition \ref{5.3} will rely on the next lemma.

\begin{lemm}\label{5.4}
 If $H_{D_j,h/2}$ is defined in (\ref{(2.6)}), with $E$ replaced
by $D_j= {\rm Vect}(e_j)$, we may write:
\begin{equation}
 I - H_{D_j,h/2} = A_j   = B_j \partial _{ u_j} + C_j \partial _{ v_j}
 =\frac{h}{4} V_j  \Delta_j \label{(5.8)} 
\end{equation}
 where the operators $B_j$, $C_j$ and $V_j$ are bounded in the space $C_b$ of
 continuous bounded functions  on $H^2$.  Moreover:
\begin{equation}
\Vert A_j \Vert _{{\cal L}(C_b)} \leq 2 \hskip 1.5cm
 \Vert B_j \Vert _{{\cal L}(C_b)} \leq (h/2)^{1/2}
 \hskip 1cm
 \Vert C_j \Vert _{{\cal L}(C_b)} \leq  (h/2)^{1/2}
 \hskip 1.5cm
 \Vert V_j \Vert _{{\cal L}(C_b)} \leq 1. \label{(5.9)} 
\end{equation}
\end{lemm}

{\it Proof.} The first estimate in (\ref{(5.9)}) is standard. From the heat
 kernel's explicit expression, we get the first equality in (\ref{(5.8)}) with:
\begin{equation}
( B_j \varphi ) (x,\xi) = - (\pi h) ^{-1} \int_{\R^2\times [0, 1]}
e^{- \frac{1}{h} (u^2 + v^2)}  u \varphi ( x + \theta u e_j,\xi +
\theta v e_j )dudv d\theta \label{(5.10)}
\end{equation}
\begin{equation}
 ( C_j \varphi ) (x,\xi) = - ( \pi h)^{-1} \int_{\R^2\times [0,
1]} e^{- \frac{1}{h} (u^2 + v^2)}v \varphi ( x + \theta u e_j,\xi +
\theta v e_j )dudv d\theta \label{(5.11)} 
\end{equation}
$$
 (V_j \varphi) (x,\xi) = - ( \pi h)^{-1} \int_{\R^2\times [0,
1]} e^{- \frac{1}{h} (u^2 + v^2)} 2 \theta \varphi ( x + \theta u e_j,\xi +
\theta v e_j )dudv d\theta 
$$
We then deduce the norm estimates of $C_j$ and $V_j$ in (\ref{(5.9)}). The last
equality in (\ref{(5.8)}) and the bound of $V_j$ in (\ref{(5.9)}) are obtained
using integrations by parts in (\ref{(5.10)}) and (\ref{(5.11)}).
\hfill $\square$ 

\bigskip

{\it Proof of  Proposition \ref{5.3}.} For any multi-index $(\alpha ,\beta )$ 
in ${\cal M}_2(I)$, one may rewrite the operator
 $\partial _u^{\alpha} \partial _{v}^{\beta} T_{I, h}$
under the following form :
%---
$$\partial _u^{\alpha} \partial _{v}^{\beta} T_{I, h} = \prod _{j\in I}
U_j \ \partial_{u_j}^{\alpha_j}\partial_{v_j}^{\beta_j}$$
%---
with:
%---
$$ U_j =  \left \{ 
\begin{array}{lllll}
A_j &{\rm if}& \alpha _j + \beta _j \geq 2\\
B_j \partial _{u_j} + C_j \partial _{v_j} &{\rm if}& \alpha _j +\beta _j = 1\\
\frac{h}{4} V_j  \Delta _j &{\rm if}& \alpha _j + \beta _j = 0\\
 \end{array} \right. . $$
%-----
Given the norm estimates for the  operators $A_j,\dots,D_j$ obtained in Lemma
 \ref{5.4}, we then deduce that, if $F$ belongs to $ S_2(M, \varepsilon)$:
%---
$$ h^{(|\alpha |+ |\beta |)/2}  \Vert \partial_u^{\alpha}\partial_{v}^{\beta}  T_{I, h}F
\Vert _{\infty}  \leq M (2h S_{\varepsilon})^{|I|}\prod _{j\in I}
\varepsilon_j^2.
$$
%---
Since the number of elements of ${\cal M}_2(I)$ is $9^{|I|}$, we
have indeed proved  Proposition \ref{5.3}. \hfill $\square$ 

%%%%%%%%%%%%%%%%%%%%%%%%%%%%%%%%%%%%%%%%%%%%%%%%%%%%%%%%%%%%%%%%%
%%%%%%%%%%%%%%%%%%%%%%%%%%%%%%%%%%%%%%%%%%%%%%%%%%%%%%%%%%%%%%%%%
\section{Proof of Proposition \ref{3.3}.}\label{6}

Proposition \ref{3.3} will  follow from Lemma \ref{6.1}  stated
below, as will be proved at the end of the section.
The proof of Lemma \ref{6.1} will use  Lemma \ref{6.2}.

\bigskip

\begin{lemm}\label{6.1}
 Let $ F$  be in $ S_2(M, \varepsilon)$, where the sequence
 $(\varepsilon_j)_{(j\in \Gamma )}$ is square summable. Let $h$ be positive. Let
 $\widetilde F$ be a function on $B^2$ which is the stochastic extension of
 $ F$ both for the measure $\mu _{B^2,h}$ and for the measure $\mu _{B^2,h/2}$.
 Then, setting $E_n =E(\Lambda_n)$, and using the operators
 $ \widetilde H_{E_n^{\perp} ,h/2}$ defined in (\ref{(2.7)}):
\begin{equation}
 \Vert \widetilde H_{E_n^{\perp} ,h/2}\widetilde F -
 \widetilde F \Vert _{L^1 (B^2,\mu _{B^2 ,h/2})} \leq K M
 \left [ \sum _{\lambda _j \notin \Lambda_n} \varepsilon_j^2 \right ]^{1/2} + KM 
\sum _{\lambda _j \notin \Lambda_n} \varepsilon_j^2  . \label{(6.1)} 
\end{equation}
\end{lemm}

{\it Proof.} 
 We have $\widetilde H_{E_n^{\perp},h/2}( F\circ P_{E_n}) =  F\circ P_{E_n}$
 and consequently:
\begin{equation}
\Vert \widetilde H_{E_n^{\perp},h/2}\widetilde F-\widetilde F\Vert_{L^1(B^2,\mu_{B^2,h/2})}
\leq \Vert \widetilde F-F\circ P_{E_n}\Vert_{L^1(B^2,\mu_{B^2, h/2})}+
 \Vert \widetilde H_{E_n^{\perp},h/2}
(\widetilde F-F\circ P_{E_n})\Vert_{L^1(B^2,\mu_{B^2,h/2})} \label{(6.2)}
\end{equation}
One denotes by  $\nu_{n h}$ the measure $\nu_{B^2, E_n^2, h/2, h}$  defined by
 (\ref{(4.C.2)}). According to Proposition \ref{4.5}, the operator 
$\widetilde H_{E_n^{\bot},h/2}$ is bounded from
 $L^1(B^2,\nu_{n,h})$ into $L^1(B^2, \mu_{B^2,h/2})$, with a norm $\leq 1$.
Consequently, since the injection of
$L^2 (B^2,\nu_{n,h})$ in $L^1 (B^2,\nu_{n,h})$ has norm 1:
\begin{equation}
\Vert \widetilde
 H_{E_n^{\perp},h/2}\widetilde F-\widetilde F\Vert_{L^1(B^2,\mu_{B^2,h/2})} \leq 
 \Vert\widetilde F - F\circ P_{E_n}\Vert_{L^2 (B^2,\mu _{B^2,h/2})} +
\Vert  \widetilde F -  F\circ P_{E_n}\Vert _{L^2 (B^2,\nu_{n ,h})} .\label{(6.3)} 
\end{equation}
From Lemma \ref{6.2}, we have, for $m>n$:
\begin{equation}
\begin{array}{lll}
\displaystyle
\Vert (F\circ P_{E_m} ) -(F \circ P_{E_n} )\Vert_{L^2(B^2,\mu_{B^2,h})} + 
 \Vert (F\circ P_{E_m})-( F\circ P_{E_n})\Vert_{L^2(B^2,\nu_{n ,h})} \\
\displaystyle
 \leq C M \left [ \sum _{j\notin \Lambda_n} \varepsilon_j^2\right ]^{1/2} +
 C M  \sum _{j\notin \Lambda_n} \varepsilon_j^2  
\end{array}
\label{(6.4)} 
\end{equation}
where $C$ is a constant which depends only on $h$. Since $\widetilde F$ is
 the stochastic extension of $ F$ for the measure $\mu_{B^2 ,h}$ and since  $F$
 is bounded, we deduce that the sequence $ F \circ P_{E_m}$ tends to
 $\widetilde F$ in $L^2 (B^2,\mu_{B^2,h})$. For all $n$ fixed, the measures
 $\nu_{n,h}= \mu_{E_n^2,h/2}\otimes \mu_{(E_n^{\bot})^2,h}$ and
 $\mu_{B^2,h}= \mu_{E_n^2,h}\otimes \mu_{(E_n^{\bot})^2,h}$ are equivalent as a product
of equivalent measures (since $E_n$ is finite dimensional). Then, the 
convergence in probability for one of them is equivalent to  the convergence
 in probability for the other one. Since $\widetilde F$ is the stochastic
 extension of $F$ for the measure $\mu_{B^2,h/2}$, we thus deduce that, for all
$n$ fixed, the sequence $ F \circ P_{E_m}$ tends $\widetilde F$ in probability 
for the measure $\nu _{n, h}$ when $m$ tends to infinity.  Consequently, this
 sequence tends to $\widetilde F$ in $L^2 (B^2,\nu_{n, h})$ for $n$ fixed,
 $m$ going to infinity. In particular:
\begin{equation}  \label{(6.5)}
\begin{array}{lll}
\displaystyle
 \Vert  \widetilde F - F\circ P_{E_n}\Vert _{L^2 (B^2,\mu _{B^2,h/2})} +
\Vert  \widetilde F -  F\circ P_{E_n}\Vert _{L^2 (B^2,\nu_{n ,h})} \leq \\ 
\displaystyle \leq C M \left [ \sum _{j\notin \Lambda_n} \varepsilon_j^2\right ]^
{1/2}+ C M \sum _{j\notin \Lambda_n} \varepsilon_j^2 .
\end{array}
\end{equation}
The Lemma is then deduced from (\ref{(6.3)}) and (\ref{(6.5)}).\hfill $\square$ 

\begin{lemm}\label{6.2}
 With the notations and assumptions of Lemma \ref{6.1}, Inequality (\ref{(6.4)})
 holds true, with a constant $C>0$ depending only on $h$.
\end{lemm}

{\it Proof.} Let $\Lambda_m \setminus \Lambda_n =\{\alpha_1 ,\dots,\alpha_p \}$,
 with $\alpha_j \in \Gamma$.  Set $I_0 = \Lambda_n$ and
 $I_j = \Lambda_n \cup \{ e_{\alpha_1} ,\dots, e_{\alpha_j} \}$. One may write, 
from an order 2 Taylor formula, and since $F$ belongs to $S(M, \varepsilon)$:
%--
$$ (F \circ P_{E(I_j )}  -( F \circ P_{E(I_{j-1})})=\varphi_j + \psi_j$$
%--
with
%---
$$\varphi_j (x,\xi) = (\partial_{u_{\alpha_j}} F)(P_{E(I_{j-1})}(x,\xi)<e_{\alpha_j},x> 
+ (\partial_{v_{\alpha_j}} F) (P_{E(I_{j-1})} (x,\xi) < e_{\alpha_j},\xi>,$$
%---
$$ \psi_j(x,\xi)=a_j(x,\xi)< e_{\alpha_j},x>^2 +
 b_j(x,\xi) <e_{\alpha_j},x><e_{\alpha_j},\xi > +
 c_j(x,\xi)<e_{\alpha_j},\xi>^2,
$$
$$|( \partial _{u_{\alpha_j}} F) (P_{E(I_{j-1})} (x,\xi)| +
|( \partial_{v_{\alpha_j}}F)(P_{E(I_{j-1})}(x ,\xi)| \leq 2 M \varepsilon_j ,
$$
$$ |a_j(x,\xi)| + |b_j(x,\xi)| + |c_j(x,\xi)| \leq C M \varepsilon_j^2 .$$
%-----
We then deduce that, with a constant $C$ depending only on $h$:
%---
$$ \Vert \varphi_j  \Vert _{L^2 (B^2,\mu_{B^2,h})} +  \Vert \varphi_j
 \Vert _{L^2 (B^2,\nu_{n ,h})} \leq  C M  \varepsilon_{\alpha_j},  $$
%----
$$ \Vert \psi_j  \Vert _{L^2 (B^2,\mu_{B^2,h})} +  \Vert \psi_j
 \Vert _{L^2 (B^2,\nu_{n,h})} \leq  C M  \varepsilon_{\alpha_j}^2 .  $$
%----
Let us show that, when $j>k$, the functions $\varphi _j$ and $\varphi_k$ are 
orthogonal, both in $L^2(B^2,\mu_{B^2,h})$ and in $L^2(B^2,\nu_{n,h})$. Indeed,
 one may write, for $j>k$, with some suitable functions $U_{jk}$ and $V_{jk}$:
%----
$$ < \varphi_j,\varphi_k> _{L^2 (B^2,\mu_{B^2,h})} =
\int _{B^2} \Big [ U _{jk} ( P_{E_{j-1}} (X)) < u_{\alpha_j},X> +
V _{jk} ( P_{E_{j-1}} (X)) < v_{\alpha_j},X> \Big ] d\mu_{B^2 ,h}(X)$$
%---
and the same holds true for the measure $\nu_{nh}$. However, from (\ref{(1.8)}),
 the integral of the above cylinder function  is vanishing for
the measure $\mu_{B^2,h}$ and also for the measure $\nu _{n, h}$.
Thus, $\varphi_j$ and $\varphi_k$ are orthogonal if $j>k$. Then, one has:
%--
$$ ( F \circ P_{E(\Lambda _m} )  -( F \circ P_{E(\Lambda_n )} )
= \sum _{j=1}^p  (\varphi_j + \psi_j)$$
%---
and, from the orthogonality property:
%---
$$ \Vert ( F \circ P_{E_m} )  -( F \circ P_{E_n} ) \Vert _{L^2 (B^2,\mu_{B^2,h})}
\leq \left [\sum_{j=1}^p \Vert \varphi_j\Vert_{L^2(B^2,\mu_{B^2,h})}^2\right ]^{1/2}
  + \sum _{j= 1}^p  \Vert \psi_j\Vert _{L^2 (B^2 ,\mu_{B^2,h})}   $$
%---
and this is also similarly valid for the measure $\nu_{n, h}$. We
then deduce (\ref{(6.4)}). The proof of the Lemma is finished.\hfill $\square$ 

\bigskip

{\it End of the proof of Proposition \ref{3.3}.}  By Proposition \ref{4.10},
(with $m= 0$), we have, for all $f$ and $g$ in
 ${\cal D}_D$ ($D\in {\cal F} (B')$):
%---
$$  \left | Q_h^{Weyl}(\widetilde H_{E_n^{\perp} ,h/2}\widetilde F) (f,g) -
Q_h^{Weyl}(\widetilde F)(f,g)\right | \leq \ I_{D ,0,  h}(f) \ I_{D ,0,  h} (g)
\  N_0 (H_{E_n^{\perp},h/2}\widetilde F - \widetilde F). $$
 %---
For all $Y$ in $H^2$ and for every function $F$ on $B^2$ (resp. on $H^2$), let
us denote by  $\tau_Y$  the function $X\rightarrow F(X + Y)$ on $B^2$
 (resp. on $H^2$). According to (\ref{(1.12)}):
\begin{equation} \label{(6.XXX)} 
\begin{array}{llll}
\displaystyle
 N_0 (H_{E_n^{\perp},h/2}\widetilde F -\widetilde F) &\displaystyle= \sup _{Y\in H^2 }
\Vert \tau_Y (H_{E_n^{\perp},h/2}\widetilde F-\widetilde F)\Vert_{L^1(B ^2,\mu_{B^2,h/2})}
\\& =  \sup _{Y\in H^2 } \Vert H_{E_n^{\perp} ,h/2}( \tau_Y \widetilde F )
- ( \tau_Y \widetilde F )  \Vert _{L^1 ( B ^2,\mu_{B^2 ,h/2})}.
\end{array}
\end{equation}
One then applies Proposition \ref{6.1}, replacing $F$ with $\tau_Y  F$, for all
 $Y$ in $H^2$. Since $F$ is in $ S_2(M, \varepsilon)$ too, since
the sequence  $(\varepsilon_j)$ is square summable, inequality (\ref{(6.1)})
  proves that the right term of
 (\ref{(6.XXX)}) converges to $0$ when $n$ goes to infinity. This
proves Proposition \ref{3.3}.  \hfill $\square$

%%%%%%%%%%%%%%%%%%%%%%%%%%%%%%%%%%%%%%%%%%%%%%%%%%%%%%%%%%%%%%%%%
%%%%%%%%%%%%%%%%%%%%%%%%%%%%%%%%%%%%%%%%%%%%%%%%%%%%%%%%%%%%%%%%%
\section{Wick symbol of a Weyl operator }\label{7}

For all $X$ in $H^2$, let $\Psi_{Xh}$ be the coherent state defined by 
(\ref{(4.D.3)})  as a function on $H$ and $\widetilde \Psi_{X h}$ its
stochastic extension, defined almost everywhere on $B$ by 
(\ref{(4.D.4)}). For every operator $A$, bounded on  $L^2(B,\mu_{B,h/2})$,
one can define the Wick symbol of $A$ by 
 \begin{equation}
\sigma_h^{Wick} (A) (X) = < A \widetilde \Psi_{X h},\widetilde \Psi_{X h} >.
\label{(7.1)}
 \end{equation}
See Berezin \cite{Be} or Folland \cite{F}.

\begin{theo}\label{7.1}
 Let  $F$  be a function in $S_2(M, \varepsilon)$,  admitting a stochastic
 extension $\widetilde F$ (Definition  \ref{4.4}, second point, $p=2$).
One denotes by $Q_h^{Weyl} (\widetilde F)$ the quadratic form associated with 
$\widetilde F$ by Definition \ref{1.2} and $Op_h^{Weyl}(F)$  the bounded 
extension whose existence results from Theorem \ref{1.4}.
Then one has, for all $X$ in $H^2$:
\begin{equation}
 \sigma_h^{Wick} (Op_h^{Weyl}(F)) (X) = \int _{B^2} \widetilde F (X + Y)
 d\mu _{B^2,h/2} (Y). \label{(7.2)} 
\end{equation}
\end{theo}

According to Gross (\cite{G-4}, Prop.9) or Kuo (\cite{K}, (Th.6.2, Chapter 2))
the integral (\ref{(7.2)}) defines a continuous function on  $H^2$.

\bigskip

{\it  Proof.} We begin by proving that, for all $X$ and $Y$, in $H^2$:
\begin{equation} < Op_h^{Weyl}(F) (\widetilde \Psi_{X,h}),\widetilde \Psi_{Y,h}>
= \int_{B^2} \widetilde F(Z)
e^{-\frac{1}{4h} (|X|^2+ |Y|^2)} \widetilde K_h^{Weyl} (X,Y, Z)
d\mu_{B^2,h/2}(Z) . \label{(7.3)}
\end{equation}

Let $(E_n)$  be an increasing sequence in s ${\cal F}(B')$, whose union is
 dense in $H$. According to Definition \ref{1.2}, one has, for all $n\geq 1$
and all   $X$ and $Y$ in $H$:
%----
$$<  Op_h^{Weyl}(\widetilde F) (\Psi_{X,h}\circ P_{E_n}),\Psi_{Y,h}\circ P_{E_n} >
= \int _{B^2} \widetilde F (Z) H^{Gauss }_h (\Psi_{X,h}\circ  P_{E_n} ,
 \Psi_{Y,h}\circ  P_{E_n} ) ( Z )d\mu_{B^2,h/2} (Z). $$
%----
When $n$ goes to infinity one has, by definition of the stochastic extension:
%-----
$$ \lim _{n\rightarrow +\infty} \Psi_{X,h} \circ P_{E_n}= \widetilde \Psi_{X,h}$$
%---
The convergence is in $L^2(B,\mu_{B, h/2})$. Since $Op_h^{Weyl}(F)$ is a
 bounded extension, one deduces that:
%---
$$< Op_h^{Weyl}(F) (\widetilde \Psi_{X,h}) ,\widetilde \Psi_{Y,h}> =
\lim _{n\rightarrow +\infty}< Op_h^{Weyl}(\widetilde F) (\Psi_{X,h}\circ P_{E_n})  ,
 \Psi_{Y,h}\circ  P_{E_n} > .$$
 %------
 According to Proposition  \ref{4.9}, one has
 %----
$$
\begin{array}{lll}
\displaystyle
 \lim _{n\rightarrow +\infty}\int_{B^2} \widetilde F(Z)
 H^{Gauss }_h(\Psi_{X,h}\circ P_{E_n},\Psi_{Y,h}\circ P_{E_n} ) (Z) 
d\mu_{B^2,h/2} (Z)  =\\
\displaystyle
 =  \int _{B^2} \widetilde F (Z)  e^{-\frac{1}{4h} (|X|^2+ |Y|^2)}
 \widetilde K_h^{Weyl} (X,Y, Z)   d\mu_{B^2,h/2} (Z), 
\end{array}
$$
%----
which gives   (\ref{(7.3)}). Restricting (\ref{(7.3)}) to  $Y=X$ and using
 (\ref{(4.E.7)}), one gets:
%---
$$ \sigma_h^{Wick} (Op_h^{Weyl}(F)) (X) =
 \int _{B^2} \widetilde F (Z)  e^{ \frac{2}{h} \ell _{X} (Z) -\frac{1}{h} |X|^2}
 d\mu_{B^2,h/2} (Z).  $$
%---
According to Proposition  \ref{4.2}, where $B$ is replaced with par $B^2$, $h$ 
with $h/2$ and $a$ with $X$, one indeed obtains (\ref{(7.2)}).\hfill $\square$ 

%%%%%%%%%%%%%%%%%%%%%%%%%%%%%%%%%%%%%%%%%%%%%%%%%%%%%%%%%%%%%%%%%
%%%%%%%%%%%%%%%%%%%%%%%%%%%%%%%%%%%%%%%%%%%%%%%%%%%%%%%%%%%%%%%%%
\section{Examples}\label{8}

We shall give examples of Wiener spaces, of stochastic extensions and,
finally, of functions $F$ in $S_m(M, \varepsilon)$.

%%%%%%%%%%%%%%%%%%%%%%%%%%%%%%%%%%%%%%%%%%%%%%%%%%%%%%%%%%%%%%%%%
%%%%%%%%%%%%%%%%%%%%%%%%%%%%%%%%%%%%%%%%%%%%%%%%%%%%%%%%%%%%%%%%%
\subsection{Examples of Wiener spaces}\label{8.A}

We begin by recalling how, starting from a real separable infinite-dimensional
Hilbert space $H$, one can always construct a Banach space $B$ and an injection
 $i$ such that $(i, H, B)$ is an abstract Wiener space. This construction is due
to Gross (\cite{G-3}, example 2 p. 92)). Then we recall two other classical
examples of Wiener spaces, related to the Brownian motion or to the
field theory.  To conclude, we give a more original example, related to
interacting lattices.
 
\bigskip

{\it Example 8.1.} Let $H$ be a real separable  Hilbert space, with norm
 $|\cdot |$ and let  $D$  be an injective Hilbert-Schmidt operator in $H$. Let
 $B$ be the completion of $H$ with respect to the norm  $x \rightarrow |Dx|$
and $i$ be the canonical injection from $H$ into $B$. Then  $(i, H ,B)$ is an
 abstract Wiener space.

\bigskip

{\it Example 8.2. The classical Wiener space. } In this case, $H$ is the
 Cameron-Martin space, that is, the space of all real valued functions in
 $H^1(0, 1)$ vanishing at the origin. This space has the scalar product
\begin{equation}
< u,v> = \int_0^1 u'(t) v'(t) dt. \label{(8.1)}
\end{equation}
One can choose, for $B$, the space of all continuous functions on $[0, 1]$ 
vanishing at the origin. This space is equipped with the norm
\begin{equation}
 \Vert f \Vert  _B = \sup _{t\in [0, 1]} |f(t)|.\label{(8.2)}
\end{equation}
One proves (cf Kuo \cite{K}, Chapter 1, Section 5) that, if $i$ is the natural
 injection from $H$ into $B$, the triple $(i, H,B)$ is an  abstract Wiener 
space. The measure $\mu_{B,h}$ on $B$ is the classical Wiener measure.

\bigskip

The following example is concerned with the free euclidean field of mass $m$,
 or free boson field, or free Markoff field, cf \cite{SI}. The space $H_m$
below appears in \cite{SI}. Nevertheless, \cite{SI} does not use the Wiener
spaces of Definition \ref{4.1}, but a Gaussian measure supported by
 ${\cal S}'(\R^2)$. One finds in Reed-Rosen \cite{R-R} the following example, 
using the Wiener spaces of Definition \ref{4.1}.

\bigskip

{\it Example 8.3.} Let $H$ be the space of tempered distributions $f$ on $\R^2$
whose Fourier Transform  $\widehat f$ is in $L^2$ locally and such that, for
 all $m>0$:
\begin{equation}
 \Vert f \Vert _{H_m}^2 = 2 \int _{\R^2}
\frac{|\widehat f (\xi)|^2}{ m^2 +|\xi|^2 } d\lambda (\xi) < + \infty .
\label{(8.3)} 
\end{equation}
An example of a Banach space $B$ such that  $(i, H, B)$ is a Wiener space
 is developed in \cite{R-R}. Let $B$ be the completion of  ${\cal S}(\R^2)$
with respect to the following norm, where $\alpha >0$ and $\beta> 1/2$:
%---
$$\Vert f \Vert_B = \Vert P_1^{-\alpha}Q_{\beta}^{-1} P^{-1} f \Vert_{L^2},$$
%---
$$ P = (1-\Delta)^{1/2} \qquad  P_1 = ( 1-\partial_{x_1}^2 )^{1/2}
 \qquad  Q_{\beta} (x) = (1 +|x|^2) ^{1/2} [\ln ((1 +|x|^2)^{1/2}) ]^{\beta}. $$
%---
Let  $i$ be the natural injection from  $H$ in  $B$. Then $(i, H,B)$ is an
abstract Wiener space (cf \cite{R-R}). Indeed, $P$ is an isomorphism from
 $L^2(\R^2)$ on $H$. Reed and Rosen prove (Lemma 1) that
 $ P^{-1}P_1^{-\alpha}Q_{\beta}^{-1}$ is a Hilbert-Schmidt operator in $L^2(\R^2)$ 
 if $\alpha>0$ and $\beta > 1/2$. Hence $P(P^{-1}P_1^{-\alpha}Q_{\beta}^{-1})P^{-1}$
is Hilbert-Schmidt in $H$. The statement of  \cite{R-R} is, then,
 a consequence of example 8.1.

\bigskip

The following example shows that our results apply to interacting lattices. 
 The Wiener space $(i, H, B)$ which is defined there exceeds the frame of
 Example 8.1.

\begin{prop}\label{8.4}
Let  $H = \ell^2 (\Gamma,\R)$, where $\Gamma$ is a countable set, for example a
lattice in $\Z^d$. Choose a family  $b =(b_j)_{(j\in \Gamma)}$  of real positive
numbers such that, for all $\varepsilon >0$, the family of positive real numbers
\begin{equation} 
R_j (b_j,\varepsilon ) = \int_{\varepsilon b_j}^{+\infty} e^{-\frac{x^2}{ 2}} dx
\hskip 2cm j\in \Gamma \label{(8.4)}  
\end{equation}
is summable.
 Let $B$ denote the space  of all families $(x_j) _{(j\in \Gamma)}$ such that
  $\left ( \frac{|x_j |}{ b_j} \right )_{(j\in \Gamma )}$ is bounded and converges
 to zero when $j$ goes to infinity. This space has the norm
\begin{equation}
 \Vert x \Vert _{B} = \sup _{j\in \Gamma} \frac{|x_j|}{ b_j} .\label{(8.5)}  
\end{equation}
In this case, the space $H=\ell^2(\Gamma,\R)$ is densely embedded in $B$. 
Moreover, for all  $h >0$, the restriction to $H$ of the norm of $B$ is
 $H-$measurable.
\end{prop}

{\it  Proof.} A seminorm $N$ on a Hilbert space $H$ is said to be tame if
 it has the form $N(x)=\widetilde N (\pi_E(x))$, where $E$ is a subspace in
 ${\cal F}(H)$, $\pi_E:H \rightarrow E$ is the orthogonal projections on $E$
 and  $\widetilde N$ is a norm on $E$. One chooses an increasing sequence of
 finite subsets $\Gamma_p$ of $\Gamma $, whose union is $\Gamma $. One
 defines an increasing sequence $N_p$ of tame semi-norms on 
 $H=\ell^2(\Gamma,\R)$ by setting:
\begin{equation} 
 N_p(x)  = \sup _{j\in \Gamma _p}\frac{|x_{j}|}{ b_j} \ . \label{(8.6)}
\end{equation}
It is well known (cf. L. Gross \cite{G-3}, Theorem 1, p. 95), that, if $N_p$ is
 an increasing sequence of tame semi-norms on $H$ and if, for all  $h>0$ and
 $\varepsilon >0$,
\begin{equation} \lim _{p\rightarrow \infty}\mu_{ H,h  }\left ( \{  x \in H,\ \ \ \
N_p(x)  \leq \varepsilon \} \right ) > 0 \ ,\label{(8.7)}
\end{equation}
then  $\lim _{p \rightarrow \infty } N_p(x)  $ exists for all $x\in H$ and the
 limit defines a  $H-$measurable semi-norm (see Definition \ref{4.1}).
For all  $\varepsilon >0$, the set
%---
$$ C_p = \{ x \in H,\ \ \ \ N_p(x)  \leq \varepsilon \}$$
%---
is a cylinder set of $H$. It can be written $C_p=\pi_E^{-1}(\Omega)$, with
 $E=\R^{\Gamma_p}$ and
 $\Omega=\prod_{j\in\Gamma_p}[-\varepsilon b_j,\varepsilon b_j]$.  Its
 $\mu_{ H, h}$- measure, defined in  (\ref{(4.1)}),  is therefore
%----
$$ \mu_{H,h} (C_p) =  \prod _{j\in \Gamma_p} (2\pi h)^{-1/2}
\int _{-\varepsilon b_j}^{\varepsilon b_j} e^{-\frac{x^2}{ 2h}} dx\ .
$$
%---------
The sequence $ \mu_{ H,h} (C_p)$ decreases to a nonnegative limit. One has
%------
$$ \mu_{ H, h} (C_p) =  \prod _{j\in \Gamma_p}
\left [1-2(2\pi h)^{-1/2}\int_{\varepsilon b_j}^{+\infty} e^{-\frac{x^2}{2h}} dx\right ]
 =\prod _{j\in \Gamma _p}\left [ 1 - 2 (2\pi )^{-1/2}
 R_j( b_j ,\frac{ \varepsilon}{ \sqrt {h}} ) \right ]\,$$
%---
where $R_j (., . )$ is defined by (\ref{(8.4)}). Since the factors of this
product are positive, the limit is positive provided the family 
 $R_j (b_j ,\frac{\varepsilon}{\sqrt{ h}} )$  is summable. Hence the
sequence  $( N_p )$ satisfies (\ref{(8.7)}) and  the limit of this sequence,
which is the restriction of the norm of $B$ to $H$, is $H-$ measurable.
\hfill  $\square$ 

\bigskip

For example if $\Gamma = \Z^d$ and if $\gamma >0$, the family of real numbers
$b_j = (1+|j|)^{\gamma}$  is such that the family $R_j(b_j,\varepsilon)$ defined
in (\ref{(8.4)}) is summable for all $\varepsilon >0$. Hence this family $(b_j)$
can be used in the definition of the Banach space before (\ref{(8.5)}).

%%%%%%%%%%%%%%%%%%%%%%%%%%%%%%%%%%%%%%%%%%%%%%%%%%%%%%%%%%%%%%%%%
%%%%%%%%%%%%%%%%%%%%%%%%%%%%%%%%%%%%%%%%%%%%%%%%%%%%%%%%%%%%%%%%%
\subsection{Examples of stochastic extensions.}\label{8.B}

The most natural example of a stochastic extension is the continuous extension,
using density arguments. The following proposition is not very different 
from  \cite{K}, Chapter 1, Theorem 6.3.

\begin{prop}\label{8.5}
 If $f$ is bounded on $H$ and uniformly continuous with respect to the 
restriction, to $H$, of the norm of $B$, then the stochastic extension and 
the extension of $f$ obtained by density are almost everywhere the same.
\end{prop}

 Kree, R\c aczka and B. Lascar sometimes supposed the symbols were Fourier
 transforms of bounded measures on  $H^2$. The following proposition shows
that such functions admit stochastic extensions.

\begin{prop}\label{8.6}
Let  $\nu$ be a  complex bounded measure on the Borel $\sigma$-algebra of
a real separable Hilbert space $H$. Let $(i, H, B)$ be an abstract Wiener
space. For every $x$ in $H$, set 
\begin{equation}
 \widehat \nu (x) = \int_H e^{i u \cdot x} d\nu (u). \label{(8.B.1)}
\end{equation}
Then for every positive $h$, the function $\widehat \nu $ has a stochastic
 extension for the measure $\mu_{B, h}$, in the sense of Definition  \ref{4.4}.
\end{prop}

{\it Proof.} Let $(E_n)$  be an increasing sequence of ${\cal F}(H)$, whose
 union is dense in $H$. Let $\pi_n: H \rightarrow E_n$  be the orthogonal
 projection on $E_n$ and  $\widetilde \pi_n: B\rightarrow E_n $  its stochastic
 extension,  defined by (\ref{(4.B.1)}). If $m<n$, one has, using Fubini's
theorem and  Jensen's inequality:
%-----
$$ \Vert \widehat \nu \circ \widetilde \pi_n -
\widehat \nu \circ \widetilde \pi _m \Vert _{L^2(B,\mu _{B,h})}^2
\leq  |\nu| (H)  \int_{H\times B} \left | e^{ i u \cdot \widetilde \pi_n (x)} -
 e^{ i u \cdot \widetilde \pi_m (x)} \right |^2 d|\nu| (u) d\mu _{B,h} (x). $$
%---
According to (\ref{(4.7)}) and (\ref{(4.B.2)}), for all $u$ in $H$:
%-----
$$ \int_{ B} \left| e^{ iu\cdot \widetilde \pi_m (x)}- e^{i u\cdot \widetilde \pi_n (x)}
 \right |^2 d\mu _{B,h} (x) \leq 
\inf (2,  h | \pi_m (u) - \pi_n (u)|^2). $$
 %-----

Hence
 %-----
$$ \Vert \widehat \nu \circ \widetilde \pi _n -
\widehat \nu \circ \widetilde \pi _m \Vert _{L^2(B,\mu _{B,h})}^2
\leq |\nu|(H) \int_{H}\inf (2,h | \pi_n (u) - \pi_m (u)|^2) d|\nu|  (u) .$$
 %----
According the dominated convergence Theorem, the sequence of functions 
$\widehat \nu \circ \widetilde\pi_n$ is a Cauchy sequence in $L^2(B,\mu_{B,h})$.
The same argument shows that the limit does not depend on the sequence
of subspaces.  Therefore, the function $\widehat \nu$  admits a stochastic
 extension in the sense of $L^2(B,\mu_{B,h})$ (Definition  \ref{4.4}).

\bigskip

We now prove that, if the sequence $(\varepsilon_j)$ is summable, Hypothesis
 (H1) of Theorem \ref{1.4} implies Hypothesis (H2), that is the existence
of a stochastic extension.

\begin{prop}\label{8.7}
Let $F$ be a function in $S_1(M,\varepsilon)$, with respect to a Hilbert basis
$(e_j)_{(j\in \Gamma )}$, where the sequence $(\varepsilon_j) _{(j\in \Gamma)}$ 
is summable. Then, for every positive $h$, $F$ admits a stochastic
extension in $L^1(B^2,\mu _{B^2,h})$.
\end{prop}

{\it  Proof.} Up to a change of notations, one can suppose that $\Gamma$ is the
set of integers greater than $1$. For all $j\geq 1$, set
$u_j = (e_j,0)$ and $v_j =(0,e_j)$. Let $F_p$  be the subspace spanned by the
 $u_j$ and  $v_j$ ($j\leq p$). Let $(E_n)$ be an increasing sequence of
 ${\cal F}(H^2)$, whose union is dense in $H^2$. For all $m$ and $n$ such
that $m<n$, let $S_{mn}$  be the orthogonal supplement of  $E_m$ in  $E_n$.
We now prove that:
\begin{equation}
\int_{B^2}|F(\widetilde \pi_{E_m}(X))-F(\widetilde \pi_{E_n}(X))| d\mu_{B^2,h}(X)
\leq M \sum _{j=1}^{\infty} \varepsilon _j
 \Big (|\pi_{S_{mn}}(u_j)| +|\pi_{S_{mn}}(v_j) |\Big ). \label{(8.B.2)} 
\end{equation}

 For all $p\geq 1$, for all  $X$ and $Y$ in $H^2$ set:
%----
$$ F_p (X,Y) = F( \pi_{F_p} (Y) +  \pi_{(F_p)^{\perp} }(X)) - F(X).$$
%-----
Since $F$ is continuous, one has:
 \begin{equation}
\lim _{p\rightarrow \infty} F_p (X,Y) = F(Y) - F(X) \hskip 2cm (X,Y)\in H^4.
\label{(8.B.3)}
 \end{equation}
One can write
%---
$$ F_p(X,Y) = \sum_{j=1}^p \varphi_j(X,Y) \hskip 2cm \varphi_j (X,Y)=
F(\pi_{F_j}(Y) +\pi_{(F_j)^{\perp}}(X))- F(\pi_{F_{j-1}}(Y) +\pi_{(F_{j-1})^{\perp}}(X)).  $$
%----
According to Taylor's formula, there exist functions $A_j(X,Y)$ and
$B_j (X,Y)$ such that
%----
$$ \varphi_j (X,Y) = A_j(X,Y) u_j \cdot (X-Y) + B_j(X,Y) v_j \cdot (X-Y)
\hskip 2cm |A_j(X ,Y)|\leq M\varepsilon_j 
\qquad |B_j (X,Y)|\leq M \varepsilon_j.
$$
%----
Therefore,
%---
$$ |F_p(X,Y)|\leq M \sum_{j=1}^p\varepsilon_j (|u_j\cdot(X-Y)| +
|v_j\cdot(X-Y)|). $$
%------
It follows that, if $m<n$, for almost every $X$ in $B^2$:
%----
$$ |F_p(\widetilde \pi _{E_m}  (X),\widetilde \pi _{E_n}  (X) ) |
\leq M  \sum_{j=1}^p \varepsilon _j ( | u_j \cdot \widetilde \pi _{S_{mn}}  (X) | +
| v_j \cdot \widetilde \pi _{S_{mn}}  (X) | ) .$$
%------
According to (\ref{(4.B.2)}):
%----
$$ |F_p(\widetilde \pi _{E_m}  (X),\widetilde \pi _{E_n}  (X) ) |
\leq M  \sum_{j=1}^p \varepsilon _j ( |\ell _{\pi _{S_{mn}}(  u_j) }     (X) | +
|\ell _{\pi _{S_{mn}}(  v_j) }     (X) |  ). $$
%------
According to (\ref{(4.3)}), in $L^1(B^2, \mu _{B^2,h})$:
%----
$$\int_{B^2} |F_p(\widetilde \pi _{E_m}(X),\widetilde \pi_{E_n}(X))|d\mu_{B^2,h}(X)
\leq M  \sqrt {2h/\pi }\sum_{j=1}^p \varepsilon_j (|\pi_{S_{mn}}(u_j)| +
|\pi_{S_{mn}}(v_j)|) .$$
%------
Letting $p$ grow to infinity and using (\ref{(8.B.3)}), one obtains 
(\ref{(8.B.2)}). Since the sequence $(\varepsilon_j)_{(j\in \Gamma)}$ is summable,
the right term of (\ref{(8.B.2)}) converges to $0$ when $m$ goes to infinity, 
according to the dominated convergence Theorem, which proves the proposition.

\bigskip

We now give examples of functions belonging to $S_1(1,\varepsilon)$, where the
sequence $(\varepsilon_j)_{(j\in \Gamma )}$ is only square summable and not
summable and who, nevertheless, admit stochastic extensions. Indeed, one will 
prove in Section \ref{8.C} that, in the proposition below, if the sequence
 $(g_j)$ is bounded and if the function $V$ is $C^2$ with bounded derivatives
of order $\leq 2$, then the function $A_t$ ($t>0$) is in $S_1(1,\varepsilon)$,
 with $\varepsilon_j = C g_j$.

\begin{prop}\label{8.8}
Set $\Gamma = \Z^d$, $H= \ell^2(\Z^d)$. Let $B$ be a Banach space such that
$(i, H, B)$ is an abstract Wiener space. Let $V:\R \rightarrow \overline{\R}_+$
be a lipschitzian function. Let $(g_j)_{(j\in\Gamma)}$  be a square summable
 sequence of positive 
real numbers such that, if $|j-k|\leq 1$, the quotient 
$g_j/g_k$ is smaller than a constant $C>0$. For all $x $ in $H$, set:
\begin{equation}
p(x)=\sum_{ \begin{array}{c} (j,k)\in\Gamma\times \Gamma \\ |j-k|=1
  \end{array}}g_jg_k V( x_j - x_k)
\hskip 2cm A_t(x) = e^{-t p(x)} \hskip 2cm (t>0).\label{(8.B.4)}
\end{equation}
Then the function $A$ admits a stochastic extension $\widetilde A_t$ in $B$
in the sense of Definition \ref{4.4}.
\end{prop}

{\it  Proof.} Let $(E_n)$  be an increasing sequence in ${\cal F}(H)$, whose
reunion is dense in  $H$. Let $\widetilde \pi_{E_n}$ be the application  
defined in Proposition \ref{1.1}. Let us prove that the sequence
$(p \circ\widetilde \pi_{E_n})$ is a Cauchy sequence in $L^1(B,\mu_{B,h})$. 
Denoting by  $(e_j)_{(j\in \Gamma)}$ the canonical basis of $H$, one has, for all
 $m$ and $n$ such that $m<n$  and for all $x\in B$:
%------
$$ |p\circ\widetilde \pi_{E_n}(x) -p \circ \widetilde \pi_{E_m}(x)| \leq 
\sum_{ \begin{array}{c} (j,k)\in\Gamma\times \Gamma \\ |j-k|=1  \end{array}}g_jg_k
 |V(\ell_{\pi_{E_n}(e_j -e_k)} (x))-V( \ell _{\pi_{E_m}(e_j - e_k)} (x))|$$
%---
where $\pi_{E_n} : H \rightarrow E_n$ is the orthogonal projection. If $K$ is
a constant such that $|V(x) - V(y)|\leq K|x-y|$  for all real numbers
$x$ et $y$, one has:
%---
$$\int_B |p \circ\widetilde \pi_{E_n}(x) -p\circ \widetilde \pi_{E_m}(x)|
 d\mu_{B,h}(x)\leq K 
\sum_{ \begin{array}{c}(j, k)\in\Gamma\times \Gamma \\ |j-k|=1  \end{array}}g_jg_k 
\int_B \ell_{\pi_{E_n}(e_j-e_k)}(x)- \ell_{\pi_{E_m}(e_j-e_k)}(x))|d\mu_{B,h}(x). $$
%---
According to (\ref{(4.3)}), one deduces:
%---
$$\int_B  |p\circ \widetilde \pi_{E_n}(x) - p\circ \widetilde \pi_{E_m}(x)|
 d\mu_{B,h}(x) \leq K \sqrt {2h/\pi}
\sum_{\begin{array}{c} (j, k)\in \Gamma\times \Gamma \\ |j-k|=1 \end{array}}g_jg_k |
\pi_{E_n}(e_j - e_k) -\pi_{E_m}(e_j - e_k) |.$$
%---
Under our hypotheses, the family $(g_j)_{(j\in \Gamma)}$ is square summable and
the sequence of the projections  $\pi _{E_n}$ converges strongly to the
identity. Hence the right term converges to $0$ when $m= \inf(m,n)$ goes 
to infinity. The sequence of functions $( p \circ \widetilde \pi_{E_n})$ is thus
 a Cauchy sequence in $L^1(B,\mu _{B,h})$. The same arguments prove that the
limit does not depend on the sequence $(E_n)$. One concludes that the sequence
$(A_t\circ \widetilde \pi_{E_n})$  converges in probability, which proves 
the proposition. \hfill  $\square$ 

\bigskip

The following examples of stochastic extensions are concerned with linear,
quadratic or harmonic functions.

\begin{prop}\label{8.9}
Let $(i, H,B)$ be an abstract Wiener space satisfying (\ref{(1.4)}). For all
$(a,b)$ in the complexified $H_{\bf C}$, one defines two functions $\varphi $ 
and $F$ on $H^2$ by:
\begin{equation}
\varphi (x,\xi) = a\cdot x +b \cdot \xi \hskip 2cm
F(x,\xi) = e^{\varphi (x,\xi)}. \label{(8.B.5)}
\end{equation}
These functions have stochastic extensions $\widetilde \varphi $
and $\widetilde F$ (in the sense of Definition  \ref{4.4}, with $p=2$).
\end{prop}

{\it  Proof.}  It suffices to set 
%----
$$\widetilde \varphi (x,\xi ) = \ell _a (x) + \ell _b (\xi) \hskip 2cm
\widetilde F (X) = e^{\widetilde \varphi (X)}. $$
%----
Let $(E_n)$ be an increasing sequence in ${\cal F}(H)$, whose union is dense
 in $H$. According to (\ref{(4.B.2)}) and (\ref{(4.3)}) (with $p=2$), one has::
%----
 $$ \Vert \widetilde \varphi-
\varphi \circ\widetilde \pi_{E_n} \Vert_{L^2(B^2,\mu_{B^2, h/2})}^2
 = h (| a - \pi _{E_n}a|^2 + | b - \pi _{E_n}b|^2).$$
 %----
Hence $\widetilde \varphi $ is really a stochastic extension of $\varphi$ in
 $L^2 (B^2, \mu _{B^2, h/2})$ in the sense of Definition  \ref{4.4}.
Moreover, still according to (\ref{(4.B.2)}):
%----
 $$ \Vert \widetilde F -F \circ\widetilde \pi_{E_n} \Vert_{L^2(B^2,\mu _{B^2, h/2})}^2|
 = \int _{B^2}\left | e^{\ell_a (x) + \ell_b (\xi)}
 - e^{\ell_{ \pi _{E_n}a } (x) +\ell_{ \pi _{E_n}b } (\xi)}\right |^2 d\mu _{B^2,h/2} (x,\xi). $$
  %----
Using (\ref{(4.7)}) proves that the right term converges to $0$, which proves
 the proposition. 

\begin{prop}\label{8.10}
Let$(i, H,B)$ be an abstract Wiener space satisfying (\ref{(1.4)}). Let $T$
be a positive, self-adjoint, trace-class operator in $H ^2$. One defines a
function  $\varphi$ and, for all positive $t$ a function $F_t$ on $H^2$ by:
\begin{equation}
 \varphi (X ) = < TX,X> \hskip 2cm F_t(X) = e^{-t \varphi (X)} \hskip 2cm
X\in H^2. \label{(8.B.6)} 
\end{equation}
Then the function $\varphi$ has a stochastic extension $\widetilde \varphi $
in $L^2( B^2,\mu _{B^2,h/2})$ (in the sense of Definition  \ref{4.4}). The norm
$N_2 (\widetilde \varphi)$  defined in (\ref{(1.12)}) is finite. The function
$F_t$ has a stochastic extension $\widetilde F_t$, (in the sense of Definition 
 \ref{4.4}). For every bounded continuous function $G$, the function 
 $X \rightarrow G(\varphi (X))$ admits a stochastic extension   in the sense of
Definition  \ref{4.4}.
\end{prop}
 
{\it  Proof.} Let  $(u_j)_{(j\geq 0)}$ be a Hilbert basis of $H^2$ consisting
in eigenvectors of $T$ and let $(\lambda_j)$  be the corresponding eigenvalues.
Since the family  $(\lambda_j)$ is summable, the sequence of functions
%----
$$ S_n (X) = \sum _{j= 0}^n \lambda _j \ell _{u_j}(X)^2 $$
%----
is a Cauchy sequence in  $L^2(B^2,\mu _{B^2,h/2})$. This is a consequence of
the Wick Theorem (Theorem \ref{4.3}). Indeed, one has for all $u$ and $v$ in
 $H^2$, (with $2p=4$),
\begin{equation}
 \int _{B^2} \ell _u (X)^2 \ell _v (X)^2 d\mu _{B^2,h/2} (X) =
\frac{h^2}{ 4} \Big [ |u|^2 \ |v|^2 + 2 ( u\cdot v )^2 \Big ]. \label{(8.B.7)}
\end{equation}
The limit of this sequence of functions is denoted by
%----
$$ \widetilde \varphi (X) =\sum _{j\geq 0} \lambda _j \ell _{u_j}(X)^2.
$$
%---
One defines a measurable bounded function by:
%----
$$ \widetilde  F_t(X) = e^{ -t \widetilde \varphi (X) }.$$
%---------
Let $(E_n)$ be an increasing sequence in ${\cal F}(H)$,
 whose union is dense in  $H$. According to (\ref{(4.B.2)}):
%---
$$ \varphi \circ \widetilde \pi _{E_n} (X) = \sum _{j\geq 0}
\lambda _j (\ell _{\pi_{E_n} (u_j)} (X))^2. $$
%----
Hence
%-----
$$ \Vert \widetilde \varphi -
\varphi \circ \widetilde \pi_{E_n}\Vert_{L^2(B^2,\mu_{B^2,h/2})}^2
\leq \sum |A_{jk}^{(n)}|  |\lambda _j \lambda_k| $$
%----
$$ A_{jk}^{(n)} = \int _{B^2} \Big [\ell _{u_j}(X)^2 -\ell_{\pi_{E_n}(u_j)} (X)^2 \Big ]
 \Big [ \ell _{u_k} (X)^2 -\ell_{\pi_{E_n}(u_k)} (X)^2 \Big ] d\mu _{B^2,h/2} (X).$$
%----
Using the particular case (\ref{(8.B.7)}) of Wick's theorem one gets:
$$
 A_{jk}^{(n)} = \frac{h^2}{ 4} \Big [ 1 + 2 (u_j \cdot u_k)^2 +
|\pi_{E_n}(u_j) |^2\  |\pi_{E_n}(u_k) |^2 - 2 ( u_j \cdot \pi_{E_n}(u_k) )^2 
 - |\pi_{E_n}(u_j) |^2 - |\pi_{E_n}(u_k) |^2 \Big ]
$$
The sequence $A_{jk}^{(n)}$ is thus uniformly bounded and converges to $0$ when
$n$ goes to infinity, for all fixed $j$ and $k$. Hence if  $T$ is trace-class,
%----
$$ \lim _{n\rightarrow +\infty}  \Vert \widetilde \varphi - \varphi \circ \widetilde
 \pi_{E_n} \Vert _{L^2 ( B^2,\mu _{B^2,h/2})}^2 = 0.$$
%----
Therefore $\widetilde \varphi$ is really a stochastic extension of $\varphi $
in $L^2 ( B^2,\mu _{B^2,h/2})$. One deduces easily that, if $\varphi \geq 0$
 the function $\widetilde F$ is a stochastic extension of
$F$ in the sense of Definition  \ref{4.4}.

\bigskip

A necessary condition for a function $f$, continuous on $H$, to admit a
 stochastic extension in the $L^2 (B,\mu _{B,h})$ sense (cf. Definition  
\ref{4.4}), is that there exists a positive $M$ such that for all subspace
 $E$  in ${\cal F}(H)$:
\begin{equation} \Vert f \circ \tilde{\pi}_E \Vert_{L^2(B,\mu _{B,h})} \leq M.
 \label{(8.B.8)}
\end{equation}
According to  (\ref{(1.7)}), it is equivalent to the fact that, for every 
subspace $E$ in  ${\cal F}(H)$, one has:
\begin{equation}
 \Vert f  \Big | _E \Vert _{L^2 (E, \mu _{E,h})} \leq M. \label{(8.B.9)}
\end{equation}
But this condition is not sufficient.

\begin{theo}\label{8.11}
 Let $(i, H,B)$ be an  abstract Wiener space. Let $f$  be a continuous
 function on $H$. Suppose that the restrictions of $f$ to  the subspaces in
 ${\cal F}(H)$ are harmonic functions. Suppose there exist  $h>0$ and  $M>0$
such that, for every subset $E$ in ${\cal F}(H)$, (\ref{(8.B.8)}) holds.
Then  $f$ admits a stochastic extension in the  $L^2 (B ,\mu _{B ,h})$ sense.
\end{theo}

{\it  Proof.} For every subspace $E$ in ${\cal F}(H)$, set
 $f_E =f \circ\tilde{\pi}_{E}$. One notices that, for all $E$ and $F$ in 
${\cal F}(H)$,  such that $E \subset F$ and for every bounded and measurable 
function  $g: B \rightarrow {\bf C}$, one has
\begin{equation}
 \int _{B} [ f_{F}(x) - f_{E}(x) ]\  g\circ \tilde{\pi}_{E} (x) d\mu_{B, h}(x)=0.
  \label{(8.B.10)} 
\end{equation}
In order to prove it, let us take an orthonormal basis $\{ e_1, ...,e_n \}$
of $F$, such that $\{ e_1,...,e_m\}$ is a basis of $E$ $(m< n)$. Set
 $x' =(x_1,x_m)$, $x'' = (x_{m+1} ,... x_n)$,  $x = (x',x'')$ and:
%---
$$ \varphi_1 (x' ) = f( x_1 e_1 + ... + x_m e_m ) \hskip 2cm
\varphi_2 (x  ) = f( x_1 e_1 + ... + x_n e_n ),$$
%---
$$ \psi (x' ) = g(x_1 e_1 + ... + x_m e_m ).$$
%----
According to (\ref{(1.7)}), the left term of (\ref{(8.B.10)}) can be written as:
%---
$$ I( \varphi _1 ,\varphi_2,\psi)= \int _{\R^{n}} [ \varphi _2 (x) -
 \varphi _1 (x') ] \psi (x' )d\mu_{\R^{n},h} (x ).$$
%---
Since the function  $\varphi_2$ is harmonic, one has for all $x'$,
 according to the mean-value property:
%-----
$$ \int _{\R^{ ( n-m)}} \varphi_2(x',x" ) d\mu
_{\R^{(n-m)}, h} (x" ) = \varphi _2( x',0 ) = \varphi_1(x' ).$$
%----
Hence $ I(\varphi_1,\varphi_2,\psi)=0$ and the property is satisfied. Since
 $f_{E}$ is $\tilde{\pi}_{E }$- measurable, it amounts to the fact that $f_{E}$ 
is the conditional expectancy of $f_{F}$ with respect to $\tilde{\pi}_{E }$:
\begin{equation}
 f_{E} = E( f_{F} |\tilde{\pi}_{E}) \label{(8.B.11)}.
\end{equation}
 Let $(E_n)$ be an increasing sequence in ${\cal F}(H)$, whose union is dense
 in $H$. Then the sequence of $\sigma$-algebras $\sigma (\tilde{\pi}_{E_n })$ is
a filtration. The
 equality (\ref{(8.B.11)}), with  $E$ and $F$ replaced with $E_n$ and $E_{n+1}$,
proves that the sequence of functions $f_{E_n}$ is a martingale for this
 filtration. Moreover, according to (\ref{(8.B.8)}), it is bounded in $L^2$.
There exists a function $\widetilde f$ in  $L^2(B, \mu _{B,h})$ (depending, at
 first sight, on the sequence $(E_n)_{n\in \N}$) such that, for all  $n$:
\begin{equation}
 f_{E_n} = E( \widetilde f |\tilde{\pi}_{E_n })
\hskip 2cm \lim _{n\rightarrow +\infty} \Vert f_{E_n} - \widetilde f
\Vert _{L^2(B,\mu _{B,h}) } = 0. \label{(8.B.12)} 
\end{equation}

The second point shows that the conditions of Definition \ref{4.4} are 
satisfied. 

 Now let us check that the function $\widetilde f$ constructed above
 does not depend on the sequence $(E_n)_{n\in \N}$. Consider two increasing
sequences $(E_n)$ and $(E'_n)$ in ${\cal F}(H)$, whose union is dense in $H$
and let $\widetilde f$ and $\widetilde f'$ in $L^2(B,\mu_{B,h})$ be the
 functions given by the preceding construction. We have to prove that 
 $\widetilde f = \widetilde f'$. Set $E''_n = E_n+ E'_n$. Let  $\widetilde f''$
be the function associated with the sequence $(E''_n)$, according to the first
 step. We first prove that, for every $n$, 
\begin{equation}
{E}( \widetilde f''- \widetilde f |\tilde{\pi}_{E_n }) = 0. \label{(8.B.13)} 
\end{equation}
According to the point (\ref{(8.B.12)}) of the first step, one has, for all $n$,
${E}(\widetilde f''|\tilde{\pi}_{E_n''})=f_{E''_n}$. According to 
(\ref{(8.B.11)}),
one deduces that, since $E_n \subset E''_n$:
%----
$${E}( \widetilde f'' |\tilde{\pi}_{E_n }) = {E}( f_{E''_n} |\tilde{\pi}_{E_n})
  = f_{E_n}.  $$
%-----
But the limit $\widetilde f$ of (\ref{(8.B.12)}) is the function associated
 with the sequence $(E_n)$ by the first step. Hence
%---
$${E}( \widetilde f |\tilde{\pi}_{E_n }) =  f_{E_n}.$$
%---
Therefore, for all $n$, (\ref{(8.B.13)}) is satisfied. With the notations of the
 first part of the proof, one obtains that, for every bounded function 
$\varphi : B \rightarrow \R$, measurable with respect to the 
sub-$\sigma$-algebra $ \sigma (\tilde{\pi}_{E_n })$, one can write:
\begin{equation}
\int_{B} (\widetilde f''- \widetilde f )\varphi (x)d\mu_{B,h} (x) = 0 .
 \label{(8.B.14)} 
\end{equation}
Set ${\cal C}=\cup_{n=1}^{\infty}\sigma (\tilde{\pi}_{E_n })$ and let ${\cal T}$ be
the $\sigma$-algebra generated by ${\cal C}$. Equation (\ref{(8.B.14)}) still
 holds for all ${\cal T}$-measurable and bounded functions, according to 
an consequence of the   $\pi \lambda$ theorem of Dynkin. Moreover, one can see
that ${\cal T}^*$, the $\sigma$-algebra generated by ${\cal T}$  and the
 negligible Borel sets, is indeed the Borel $\sigma$-algebra.  This establishes
 (\ref{(8.B.14)}) for all bounded Borel functions  $\varphi$ and implies that
 $\widetilde f''=  \widetilde f $.

%%%%%%%%%%%%%%%%%%%%%%%%%%%%%%%%%%%%%%%%%%%%%%%%%%%%%%%%%%%%%%%%%
%%%%%%%%%%%%%%%%%%%%%%%%%%%%%%%%%%%%%%%%%%%%%%%%%%%%%%%%%%%%%%%%%
\subsection{Examples of symbols in $S_m(M, \varepsilon)$.}\label{8.C}

Let us now give examples of symbols belonging to the sets $S_m(M,\varepsilon)$.
The example below shows that Definition \ref{1.1} is coherent with former 
articles.
\begin{prop}\label{8.12}
Let $(i,H,B)$ be an abstract Wiener space satisfying (\ref{(1.4)}). For all $a$
 and $b$ in $H$, let $F_{a,b}$ be the function on $H^2$ defined by:
\begin{equation}
F _{a,b} (x,\xi) = e^{i(a\cdot x + b \cdot \xi)} . \label{(8.C.1)}  
\end{equation}
If $(e_j)_{(j\geq 1)}$ is an arbitrary Hilbert basis of $H$ and if $m$ is
a positive integer, the function $F_{a,b}$ is in the set  $S_m(1,\varepsilon)$
associated with this basis, with
\begin{equation}
 \varepsilon _j =  \max ( | e_j \cdot a |,|e_j \cdot b|)\label{(8.C.2)}
\end{equation}
and this sequence is square summable. Let $Op_h^{Weyl}(F_{a,b})$ be the operator,
bounded on $L^2(B,\mu_{B,h/2})$, associated with  $F_{a, b}$ by Theorem \ref{1.4} 
(this function admits a stochastic extension according to Proposition 
\ref{8.9}). Then the operator $Op_h^{Weyl}(F_{a,b})$ is equal to the operator 
 $U_{a,b,h}$ defined by
\begin{equation}
 (U_{a, b,h}f)(u)= e^{-\frac{h}{2}|b|^2+{i\frac{h}{2} a\cdot b+i\ell_{a+ib}(u)}}
 f ( u + h b) \hskip 2cm a.e. u\in B .\label{(8.C.3)}
\end{equation}
\end{prop}

{\it  Proof.}  The first claim, that $F_{a,b}$ is in $S_m(1,\varepsilon)$, is
 easy to verify. The fact that the right term of (\ref{(8.C.3)}) defines a
 bounded operator in $L^2(B,\mu_{B,h/2})$ is a consequence of Proposition 
\ref{4.2}. Let us prove that $U_{a,b,h}$ is equal to the operator associated
 with the function  $F_{a,b}$ by Theorem \ref{1.4}. Since $F_{a,b}$ is in
$S_2(1,\varepsilon)$, one can apply (\ref{(7.3)}). Hence it suffices to prove,
denoting by $\widetilde F_{a,b}$ the stochastic extension of  $F_{a,b}$, that,
for all $X$ and $Y$ in $H$,
\begin{equation}
< U_{a,b,h}(\widetilde \Psi_{X,h}),\widetilde \Psi_{Y,h}>= \int_{B^2}
\widetilde F_{a,b}(Z) H^{Gauss }_h(\widetilde \Psi_{X,h},\widetilde \Psi_{Y,h})(Z)
d\mu_{B^2,h/2} (Z) \label{(8.C.4)} 
\end{equation}
Indeed, the set of functions $\widetilde \Psi_{X,h}$ ($X\in H$) is total in
 $L^2(B,\mu _{B,h/2})$ (this is, for example, a consequence of Janson \cite{J},
 Theorem D.10). One easily verifies that:
%-----
$$  U_{a,b,h} (\widetilde \Psi_{X,h}) = e^{ \frac{i}{ 2}(a \cdot x + b\cdot \xi)}
\widetilde \Psi_{ x - h b,\xi + h a,h}.$$
%-----
Consequently, using  (\ref{(4.D.5)}):
%----
$$< U_{a,b,h} (\widetilde \Psi_{X,h}),  \widetilde \Psi_{Y,h}>
= e^{ -\frac{1}{4h}  (|x-y-h b|^2 +|\xi - \eta + h a|^2) +
\frac{i}{2h}(y\cdot\xi -x\cdot \eta +h a \cdot (x+y) +h b\cdot(\xi  \eta))}. 
$$
%----
One then uses the expression of the stochastic extension of the Wigner function
 of the coherent states, given by  (\ref{(4.E.10)}) and (\ref{(4.E.7)}), as
well as (\ref{(4.2)}). This yields
$$ 
\begin{array}{lll}
\displaystyle \int_{B^2} \widetilde F_{a,b}(Z) 
\widetilde H^{Gauss }_h(\widetilde \Psi_{X,h},\widetilde \Psi_{Y,h})(Z)
d\mu_{B^2,h/2}(Z) = \\ 
\displaystyle = e^{-\frac{1}{4h} (|x-y-hb|^2 +|\xi-\eta +ha|^2) +
\frac{i}{2h}(y\cdot \xi -x\cdot \eta + ha\cdot (x+y)+ hb\cdot (\xi+\eta))}.
\end{array}
$$
%-----
One has proved (\ref{(8.C.4)}). Therefore the operator associated with the
 function $F_{a,b}$ by Theorem \ref{1.4} is equal to the operator
$U_{a,b,h}$ defined by (\ref{(8.C.3)}).

\bigskip

\begin{prop}\label{8.13}
\begin{enumerate}
\item  For all  $a$ et $b$ in $H$, the quadratic Weyl form associated (in the
sense of Definition \ref{1.2}), with the stochastic extension
$\widetilde \varphi_{a,b} $ of the function
 $\varphi_{a,b} (x,\xi) = a \cdot x + b \cdot \xi$ satisfies
\begin{equation}
 Q_h^{Weyl}(\widetilde \varphi_{a,b})(f,g)= 
<\ell_{a+ib} f+ \frac{h}{i}\ b\cdot \frac{\partial f}{ \partial u},g>
_{L^2(B,\mu_{B,h/2})}
  \hskip 2cm f, g \in {\cal D},  \label{(8.C.5)}
\end{equation}
where the function $u \rightarrow \ell _{a+i b}(u)$  is introduced in 
Proposition \ref{1.1}.
\item Let $F$  be a function on $H^2$ which can be expressed as:
\begin{equation}
 F(t,\tau ) = \int _{H^2} e^{ i (   x\cdot  t +  \xi \cdot  \tau )} d\nu (x,\xi)
\label{(8.C.6)} 
\end{equation}
where $\nu$  is a complex measure, bounded on $H^2$. Let $\widetilde F$ be the
stochastic extension of $F$, which exists according to Proposition \ref{8.6}.
Let $Q_h^{Weyl} (\widetilde F)$ be the quadratic Weyl form associated with this
stochastic extension according to Definition \ref{1.2} and
 $Q_h^{AW} (\widetilde F)$  the quadratic anti-Wick form associated with it
according to (\ref{(2.4)}). Then for all $f$ and $g$ in ${\cal D}$:
\begin{equation}
  Q_h^{Weyl} (\widetilde F) (f,g ) = \int _{H^2} < U _{x,\xi, h}f,g>
d\nu (x,\xi), \label{(8.C.7)}
\end{equation}
\begin{equation}
Q_h^{AW} (\widetilde F) (f,g )= \int_{H^2}e^{-{\frac{h}{4}}(|x|^2+ |\xi|^2)} 
< U_{x,\xi, h}f,g>d\nu (x,\xi), \label{(8.C.8)} 
\end{equation}
where $U _{x,\xi, h}$ is defined by (\ref{(8.C.3)}). There exists an operator
$Op_h^{Weyl}(F)$ bounded on $L^2(B,\mu_{B,h/2})$, whose norm is smaller than
$|\nu|(H^2)$, such that 
$Q_h^{Weyl}(\widetilde F)(f,g)= <Op_h^{Weyl}(F)f,g>$ for all 
$f$ and $g$ in ${\cal D}$.
\end{enumerate}
\end{prop}

{\it  Proof. Point 1).}  Let us first suppose that 
$a$ and $b$ belong to $B'$. Let $E$ be in ${\cal F}(B')$, let $\widehat f$ and
$\widehat g$ be functions on $E$ such that $f =\widehat f\circ\widetilde \pi_E$
and $g =\widehat g \circ\widetilde \pi_E$.  One can assume that $a\in E$ and
 $b\in E$. Let $\widehat \varphi$ be the restriction of $\varphi_{a,b}$
 to $E^2$. 
According to (\ref{(1.14)}), one has:
%----
$$Q_h^{Weyl}(\widetilde \varphi_{a,b})(f,g)= 
<Op_h^{Weyl,Leb}(\widehat{\varphi }) \gamma_{E,h/2} \widehat{f}, 
 \gamma_{E,h/2} \widehat{g}>_{L^2 (E,\lambda)}. $$
%----
Using the usual Weyl calculus in a finite dimensional setting  gives:
%----
$$Op_h^{Weyl, Leb} (\widehat{\varphi}) \gamma_{E,h/2} \widehat{f}(u) =
a\cdot u\gamma_{E,h/2}\widehat{f}(u)
+\frac{h}{i} b\cdot \frac{\partial}{\partial u}\gamma_{E,h/2}\widehat{f}(u)=
 \gamma_{E,h/2} \widehat{F} (u) $$
%-------
where we set:
%----
$$ \widehat{F}(u) = (a+i b) \cdot u \widehat f (u) + \frac{h}{i} b\cdot 
\frac{\partial \widehat f}{\partial u}. $$
%---
Since $\gamma_{E,h/2}$ is an isometry, one can write:
%----
$$Q_h^{Weyl}(\widetilde \varphi_{a,b})(f,g)= 
<\widehat F,\widehat g>_{L^2( E,\mu _{E,h/2})} =
<\widehat F \circ \widetilde \pi_E,\widehat g
\circ \widetilde \pi_E >_{L^2(B,\mu_{B,h/2})}, $$
%------
which implies (\ref{(8.C.5)})  in the case when $(a,b)\in B'^2$.\\
If $a,b\in H$, let $(a_n)_n, (b_n)_n$ be sequences of $B'$ converging 
to $a,b$ in $H$. One can then check that
$N_1(\tilde{\varphi}_{a_n,b_n}-\tilde{\varphi}_{a,b})\leq C(|a-a_n|+|b-b_n|)$
for a real constant $C$. Then Proposition \ref{4.10}
implies that $ Q_h^{Weyl}(\widetilde \varphi_{a_n,b_n})(f,g)$
converges to the left term of  (\ref{(8.C.5)}). To treat the right term,
one  proves that there exists a constant $C(f,g)$ depending on $f$
 and $g$ such that
$$
|<\ell_{a+ib} f+ \frac{h}{i}\ b\cdot \frac{\partial f}{ \partial u},g>
_{L^2(B,\mu_{B,h/2})}|\leq ((|a-a_n|+|b-b_n|) C(f,g),
$$
using (\ref{(1.7)}) and the  fact that $\gamma_{E,h/2} \hat{f},
\gamma_{E,h/2} \hat{g}$ are rapidly decreasing. 
 This proves the convergence for the right term and establishes 
(\ref{(8.C.5)}) in the general case.

\smallskip

{\it Point 2).} Let $E$ be the space of all functions $F$ in
 $L^1(B^2,\mu_{B^2,h/2})$ such that, for all $Y$ in $H^2$, the function
 $F(Y + \cdot)$ is in $L^1(B^2,\mu_{B^2,h/2})$ and such that $N_1(F)$, defined in
 (\ref{(1.12)}), is finite. This space, equipped with the norm $N_1(\cdot)$, is
a Banach space. For all $(x,\xi)$ in $H^2$, the function $\widetilde F_{x,\xi}$,
 stochastic extension of the function $F_{x,\xi}$  defined in (\ref{(8.C.1)}), is
in $E$ and the application which associates  
$\widetilde F_{x,\xi}$ with every $(x,\xi)$ in $H^2$, is continuous from $H^2$
in $(E, N_1)$.  More precisely, one can prove that, for $U,V\in H^2$
and $Y\in H^2$,
$$
\int_{B^2} \vert e^{i \ell_U(X+Y)}-e^{i \ell_V(X+Y)}  \vert\ d\mu_{B^2, h/2}(X)
= \int_{\R} \vert e^{i (\sqrt{h/2} |U-V|u + Y\cdot(U-V))}\vert \ d\mu_{\R,1}(u) 
$$
and derive from that the continuity. This allows us to define
$
\int_{H^2} \widetilde F_{x,\xi} \ d\nu(x,\xi)
$
as a Bochner integral, following \cite{Y}, chapter V. One can notice that
$
\int_{H^2} \widetilde F_{x,\xi} \ d\nu(x,\xi)
$
is the stochastic extension of 
$
(z,\zeta) \mapsto \int_{H^2} e^{i(x\cdot z + \xi\cdot \zeta)} \ d\nu(x,\xi)
$
in the sense of Definition \ref{4.4} (for any $\mu_{B^2,t} $ and $p=2$).

Finally, for all $f$ and $g$ in ${\cal D}$, the application which
associates $Q_h^{Weyl}(F) (f,g)$ with every $F$ in $E$ is continuous on $E$,
according to Proposition \ref{4.10}.
 This yields 
$$
Q_h^{Weyl}\left(\int_{H^2} \widetilde F_{x,\xi} \ d\nu(x,\xi)  \right)(f,g)
=\int_{H^2}Q_h^{Weyl}\left( \widetilde F_{x,\xi}\right)(f,g) \ d\nu(x,\xi).
$$ 
Since $Op_h^{Weyl}( F_{x,\xi})$ is
the operator $U_{x,\xi, h}$ of (\ref{(8.C.3)}), one obtains (\ref{(8.C.7)}). 
 The upper bound comes from the fact that the
operator $U _{x,\xi, h}$ is unitary in  $L^2(B, \mu _{B, h/2})$. 

 Now let us prove (\ref{(8.C.8)}). According to 
(\ref{(2.3)}) and (\ref{(4.2)}), one has the following equality 
 for $(x,\xi)\in (B')^2$ :
%-----
$$ \widetilde H_{h/2} \widetilde F_{x,\xi }  = e^{-{\frac{h}{4}} (|x|^2
+ |\xi|^2)}  \widetilde F_{x,\xi } . $$
%-------
 Since both sides
are continuous functions on  $H^2$ ( thanks to  (\ref{(4.7)}) 
 and Proposition \ref{4.5})
and  valued in $L^2(B^2,\mu_{B^2,h/2})$, it
 holds true on $H^2$.
Hence
%-----
$$
Q_h^{AW}(\widetilde F_{x,\xi })(f,g)=  
e^{-{\frac{h}{4}} (|x|^2+ |\xi|^2)}  <U_{x,\xi,h }f,g>.
$$
%----
Similar arguments allow to integrate on $H^2$, which yields
 (\ref{(8.C.8)}). \hfill  $\square$ 

\bigskip

In the articles on Fock spaces, the operator on the right side of
 (\ref{(8.C.5)}) is usually called Segal field and its exponential $U_{a,b,h}$,
Weyl operator (see \cite{R-S} vol. II, sect. X, p.209-210 
or \cite{B-R} vol.II, sect.5.2.1.2, p.12-13).

\bigskip

In B. Lascar's survey  \cite{LA-1} (Definition 10 and Proposition 1.7), one 
finds the expression (\ref{(8.C.8)}) of the anti-Wick operator for $h=2$ in 
the case when $F$ is of the form (\ref{(8.C.6)}). The operator denoted by 
$\tau_{y'+i y''}$ in Proposition 1.7 of \cite{LA-1} is equal to
 $e^{-(1/2)( |y'|^2 + |y''|^2)}U_{y'',y', 2}$. To give the analogous expression of the
Weyl operator (for $h=2$), it would have been enough, in \cite{LA-1}, to 
suppress the factor $e^{-(1/2)( |y'|^2 +|y''|^2)}$  in the integral.

\bigskip

One finds in \cite{K-R} an analogous expression, using the creation and
 annihilation operators. For all $a$ and $b$ in $H$, one denotes by $a(-b+ia)$ 
and $a^{\star}(-b+ia))$ the operators whose Weyl symbols are respectively the
functions on $H^2$ $(-b-ia)\cdot(x+i \xi)/(\sqrt{2h})$ and
$(-b+ia)(x-i \xi)/(\sqrt{2h})$ and which are not bounded on $L^2(B,\mu_{B h/2})$.
With these notations, one verifies that
\begin{equation}
U_{a,b,h} = e^{ \sqrt \frac{h}{2}  a^{\star} (-b+ia)) } e^{ -  \sqrt \frac{h}{2}  a(-b+ia ) }
e^{-{\frac{h}{4}} (|a|^2 + |b|^2)}. \label{(8.C.9)} 
\end{equation}
One finds in \cite{K-R} (formula (\ref{(6.4)}), maybe up to some insignificant
modifications), when $F$ is of the form (\ref{(8.C.6)}), an expression
comparable to (\ref{(8.C.7)}) and using the factorization (\ref{(8.C.9)}) above.

\bigskip

The last example is inspired by interacting lattices. It reminds of a lattice
 of harmonic oscillators with a coupling between nearest neighbours.

\begin{prop}\label{8.14}
Set $\Gamma = \Z^d$ and  $H= \ell^2(\Z^d)$. Let 
$V:\R \rightarrow \overline{\R}_+$ be a smooth function whose derivatives of
order greater than $1$ are bounded. Let $(g_j)_{(j\in \Gamma)}$  be a sequence of
positive real numbers such that, if $|j-k|\leq 1$, the quotient $g_j/g_k$ is
 smaller than a constant $K_0>1$. For all $(x,\xi)$ in $H^2$ and all $t>0$ set:
\begin{equation}
 f(x,\xi) = \sum _{j\in \Gamma}g_j^2  \xi_j^2  +
\sum _{\begin{array}{c}(j, k)\in\Gamma\times \Gamma \\|j-k|=1  \end{array}}
g_jg_k V( x_j - x_k)
\hskip 2cm F_t (x,\xi) = e^{-tf(x, \xi)},\label{(8.C.10)}
\end{equation}
where $|\cdot |$  refers to the norm $\ell ^{\infty}$ on $\R^d$. \\
For all $m\geq 1$, the function $F_t$ is in $S_m(1,\varepsilon^{(m)})$, with 
\begin{equation}\varepsilon_j^{(m)}  \leq  C_m \max (  g_j^2,g_j^{1/m}),
  \label{(8.C.11)}
\end{equation}
 where $C_m$ depends only on $m$, $t$, $d$, $K_0$ and on the bounds for the
derivatives of $V$, up to order $2m$.
\end{prop}

{\it  Proof.}  Let $p$ and $A_t$ be the functions defined by (\ref{(8.B.4)})
 (Proposition  \ref{8.8}), identified with functions on  $H^2$, independent on
 $\xi$. With the notations of Definition \ref{1.4}, we have, for each finitely 
supported multi-index $\alpha$ such that $\alpha \not= 0$ and $\alpha_j \leq m$
for all $j$ and for all $x\in \R^n$:
\begin{equation}
|\partial_u^{\alpha} p(x)| \leq \prod_{j\in S(\alpha)}(\lambda_j^{(m)})^{\alpha _j}
 \hskip 2cm  \lambda_j ^{(m)} = 2 \ 3^d K_0 \max ( g_j^2,g_j^{1/m})
 \max _{1\leq k \leq 2m} \left( \Vert V^{(k)}\Vert _{\infty}\right )^{1/k}
 \label{(8.C.12)} 
\end{equation}
where $S(\alpha)$ is the set of sites $j\in\Gamma$ such that $\alpha_j \not=0$.
We have used the fact that the number of elements of $\Z^d$ which are in the 
unit ball for the $\ell ^{\infty}$ norm is lesser than $3^d$. We shall use a
multidimensional version of the Faa di Bruno formula, due to Constantine Savits
 \cite{C-S}. For every multi-index $\alpha$, let $F(\alpha)$ be the set of
all applications $\varphi$ from the set of multiindices $\beta \not = 0$ 
such that $\beta \leq \alpha$ into the nonnegative integers such that 
\begin{equation}
\sum _{0 \not= \beta\leq \alpha } \varphi (\beta) \beta = \alpha . \label{(8.C.13)} 
\end{equation}
 If $\alpha \not= 0$, Constantine-Savits' formula  can be written:
\begin{equation}
\partial_u^{\alpha}e^{-tp(x)} = \alpha ! e^{-tp(x)}  \sum _{\varphi \in F(\alpha)}
\prod _{0\not= \beta \leq \alpha} \frac{1}{ \varphi (\beta)!}
 \left [\frac{\partial_u ^{\beta }(-tp(x))}{\beta! }\right ]^{\varphi (\beta) } .
 \label{(8.C.14)}
\end{equation}
In the above sum, we may replace $F(\alpha)$ with the subset $E(\alpha)$ of all
$\varphi $ in $F(\alpha)$ such that $\varphi (\beta) = 0$ for all multi-index
$\beta$ such that $\partial _u^{\beta } p(x)$ vanishes identically. By  
(\ref{(8.C.12)}) and (\ref{(8.C.13)}), we have, for each $\varphi\in F(\alpha)$:
\begin{equation}
 \prod_{0\not= \beta \leq \alpha} 
\left|\partial^{\beta}(-tp(x))\right|^{\varphi (\beta) }\leq
 t^{r(\varphi)} \prod _{j\in S(\alpha) } ( \lambda_j^{(m)} )^{\alpha_j}\hskip 2cm 
r(\varphi) = \sum_{0 \not= \beta\leq \alpha} \varphi(\beta) \label{(8.C.15)}
\end{equation}
If $\alpha_j \leq m$ for all $j$, we have $r(\varphi) \leq m |S(\alpha)|$ for
 each $\varphi \in F(\alpha)$.  We have also $\alpha !\leq (m !)^{|S(\alpha)|}$.
 For each $\beta \not= 0$ such that $\beta \leq \alpha$ and for each $\varphi$ 
in $F(\alpha)$, there exists $j$ such that $\beta_j \not = 0$, and therefore :
\begin{equation}
 \varphi (\beta) \leq \varphi (\beta)\beta_j \leq \sum_{0\not = \gamma \leq \alpha }
\varphi (\gamma ) \gamma_j = \alpha_j \leq m . \label{(8.C.16)} 
\end{equation}
Let $I(\alpha)$ be the set of all multi-indices $\beta \not =0$ such that
$\beta \leq \alpha$ and such that $\partial _u^{\beta } p$ does not vanish
identically. By (\ref{(8.C.16)}), the number of elements of $E(\alpha)$ is not
 greater than the number of applications of $I(\alpha)$ to the set
 $\{ 0, ...,m\}$, that is to say $(m+1)^{|I(\alpha)|}$. By the form 
(\ref{(8.B.4)}) of $p$, for each point $j$, the number of multi-indices
 $\beta$ such that $\beta _j \not= 0$, $\beta_k \leq m$ for all $k\in \Gamma$
 and such that $\partial ^{\beta}p(x)$ does not vanish identically, is smaller
 than $3^d m^2$.  Therefore, if $\alpha_k \leq m $ for all $k$, the number of
 elements of $I(\alpha)$ is not greater than $3^d m^2 |S(\alpha)|$. Hence
 the number of elements of $E(\alpha)$ is not greater than
 $(m+1)^{3^d m^2 |S(\alpha)|}$.
Therefore, if $t<1$, and if $\alpha_k \leq m$ for all $k$:
%----
$$ \left | \partial_u^{\alpha}e^{-tp (x)}  \right | \leq
\left [m !(m+1)^{3^d m^2}\right]^{|S(\alpha)|}
\prod_{j\in S(\alpha)}(\lambda_j^{(m)})^{\alpha_j}. $$
%-----
If $t>1$, we have:
%----
$$ \left | \partial_u^{\alpha}e^{-tp(x)}  \right | \leq t^{m|S(\alpha)|}
 \left [ m! (m+1)^{3^d m^2} \right ] ^{|S(\alpha)|}
 \prod _{j\in S(\alpha)} (\lambda _j^{(m)})^{\alpha_j} . $$
%-----
Therefore, in both cases:
\begin{equation}
|\partial_u^{\alpha}e^{-t p (x)}| \leq \prod_{j\in S(\alpha)}
(\varepsilon'_j(t))^{\alpha _j} \hskip 2cm 
\varepsilon'_j(t) = m!(m+1)^{3^d m^2}\max(1,t^m)\lambda_j^{(m)}.\label{(8.C.17)} 
\end{equation}
Now set:
%----
$$ q( \xi) = \sum _{j\in \Gamma}g_j^2 \xi_j^2 \hskip 2cm B_t(\xi) =e^{-tq(\xi) }.$$
%------
One easily checks that:
\begin{equation}
|\partial_v^{\beta}e^{-t q (\xi )}|\leq 
\prod _{j\in S(\beta)}(\varepsilon''_j(t) )^{\beta _j} \label{(8.C.18)} 
\end{equation}
with
$$\varepsilon''_j(t) = C_m g_j \sqrt {t} \hskip 2cm C_m = \max _{\alpha \leq m}
 \sup _{\xi \in \R} \left | \partial^{\alpha} e^{-\xi^2}\right |^{1/|\alpha|}.  $$
%-----
The proposition follows, setting 
$\varepsilon_j^{(m)} (t) =\max (\varepsilon'_j(t),\varepsilon''_j(t))$. 
\hfill  $\square$ 

\bigskip

If $H$ is the Hilbert space of Proposition \ref{8.14}, let $B$ be a Banach
 space containing $H$  and such that $(i, H,B)$ is an abstract Wiener space.
 The space $B$ may be the one of Proposition  \ref{8.4}. If the family
 $(g_j)_{(j\in \Gamma)}$ is summable, according to Propositions \ref{8.8} and
 \ref{8.10}, the function $ F_t$ defined on $H^2$ in (\ref{(8.C.10)}) admits
 a stochastic extension $ \widetilde F_t$ in  $B^2$. Moreover, for  $m=2$, the
 family  $(\varepsilon^{(2)}_j)_{(j\in \Gamma)}$ satisfying (\ref{(8.C.11)}) is
square summable. According to Theorem \ref{1.4}, one can associate with $ F_t$
a Weyl operator, bounded in $L^2(B,\mu _{B,h/2})$.

\end{document}